\documentclass[a4paper]{amsart} 
\usepackage[ascii]{inputenc} 

\usepackage{amssymb}
\usepackage{mathrsfs}
\usepackage[hidelinks]{hyperref}

\usepackage{xcolor}

\newcommand{\red}[1]{\textcolor{red}{#1}}

\newcommand{\magenta}[1]{\textcolor{magenta}{#1}}

\usepackage[shortlabels]{enumitem}
\setlist[enumerate,1]{label={(\Alph*)}}
\setlist[enumerate,2]{label={(\alph*)}}
\setlist[enumerate,3]{label={$\bullet_{\arabic*}$}}

\newenvironment{PROOF}[2][\proofname.]
   {\begin{proof}[#1]\renewcommand{\qedsymbol}{\textsquare$_{\mathrm{#2}}$}}
   {\end{proof}}

\newtheorem{theorem}{Theorem}[section]

\newtheorem{claim}[theorem]{Claim}
\newtheorem{conclusion}[theorem]{Conclusion}

\newtheorem{lemma}[theorem]{Lemma}
\newtheorem{observation}[theorem]{Observation}

\theoremstyle{definition}

\newtheorem{definition}[theorem]{Definition}

\newtheorem{problem}[theorem]{Problem}

\theoremstyle{remark}

\newtheorem{notation}[theorem]{Notation}
\newtheorem{question}[theorem]{Question}
\newtheorem{remark}[theorem]{Remark}
\newtheorem{note}[theorem]{Note}
\newtheorem{thesis}[theorem]{Thesis}

\newcommand{\acc}{\operatorname{acc}}
\newcommand{\Hom}{\operatorname{Hom}}

\newcommand{\rk}{\operatorname{rk}}
\newcommand{\tlim}{\operatorname{tlim}}

\newcommand{\Ens}{\mathrm{Ens}}
\newcommand{\BB}{\mathrm{BB}}

\newcommand{\cer}{\mathrm{cr}}

\newcommand{\CH}{\mathrm{CH}}
\newcommand{\Ch}{\mathrm{Ch}}
\newcommand{\ch}{\mathrm{ch}}

\newcommand{\ged}{\mathrm{gd}}
\newcommand{\Fil}{\mathrm{Fil}}
\newcommand{\NPT}{\mathrm{NPT}}
\newcommand{\pp}{\mathrm{pp}}
\newcommand{\Reg}{\mathrm{Reg}}
\newcommand{\reg}{\mathrm{reg}}
\newcommand{\rs}{\mathrm{rs}}
\newcommand{\stat}{\mathrm{stat}}

\newcommand{\bd}{\mathrm{bd}}

\newcommand{\cf}{\mathrm{cf}}

\newcommand{\CON}{\mathrm{CON}}

\newcommand{\cov}{\mathrm{cov}}

\newcommand{\Dom}{\mathrm{Dom}}

\newcommand{\GCH}{\mathrm{GCH}}

\newcommand{\id}{\mathrm{id}}

\newcommand{\nacc}{\mathrm{nacc}}

\newcommand{\nor}{\operatorname{nor}}

\newcommand{\Ord}{\mathrm{Ord}}
\newcommand{\otp}{\mathrm{otp}}
\newcommand{\pcf}{\mathrm{pcf}}

\newcommand{\Rang}{\mathrm{Rang}}

\newcommand{\tcf}{\mathrm{tcf}}

\newcommand{\tr}{\mathrm{tr}}

\newcommand{\ZFC}{\ensuremath{\mathrm{ZFC}}}

\newcommand{\wilog}{\ensuremath{\text{without loss of generality}}}

\newcommand{\bfC}{\mathbf{C}}
\newcommand{\bfc}{\mathbf{c}}

\newcommand{\bfj}{\mathbf{j}}
\newcommand{\bfK}{\mathbf{K}}

\newcommand{\bfL}{\mathbf{L}}
\newcommand{\bfM}{\mathbf{M}}

\newcommand{\bfU}{\mathbf{U}}

\newcommand{\bfV}{\mathbf{V}}

\newcommand{\bbL}{\mathbb{L}}

\newcommand{\bbR}{\mathbb{R}}

\newcommand{\bbZ}{\mathbb{Z}}

\newcommand{\Iff}{\underline{iff}}
\newcommand{\If}{\underline{if}}

\newcommand{\then}{\underline{then}}

\newcommand{\mn}{\medskip\noindent}
\newcommand{\sn}{\smallskip\noindent}

\newcommand{\cA}{\mathscr{A}}

\newcommand{\cC}{\mathscr{C}}
\newcommand{\cD}{\mathscr{D}}

\newcommand{\cF}{\mathscr{F}}

\newcommand{\cP}{\mathscr{P}}
\newcommand{\cS}{\mathscr{S}}

\newcommand{\clF}{\mathcal{F}}
\newcommand{\clH}{\mathcal{H}}
\newcommand{\clI}{\mathcal{I}}

\newcommand{\clP}{\mathcal{P}}

\newcommand{\gB}{\mathfrak{B}}

\newcommand{\gK}{\mathfrak{K}}

\newcommand{\ga}{\mathfrak{a}}
\newcommand{\gb}{\mathfrak{b}}
\newcommand{\gc}{\mathfrak{c}}
\newcommand{\gd}{\mathfrak{d}}

\newcommand{\eps}{\varepsilon}

\newcommand{\rest}{\restriction}

\newcommand{\caret}{{\char 94}}
\newcommand{\LL}{\langle}
\newcommand{\RR}{\rangle}

\usepackage[hidelinks]{hyperref}
\begin{document}
\makeatletter\def\shfiuwefootnote{\gdef\@thefnmark{}\@footnotetext}\makeatother\shfiuwefootnote{Version 2022-09-24\_3. See \url{https://shelah.logic.at/papers/E12/} for possible updates.}

\title [Analytical Guide]{Analytical Guide and updates for 
Cardinal Arithmetic \\
 E-12}

\author {Saharon Shelah}
\address{Einstein Institute of Mathematics\\
Edmond J. Safra Campus, Givat Ram\\
The Hebrew University of Jerusalem\\
Jerusalem, 91904, Israel\\
 and \\
 Department of Mathematics\\
 Hill Center - Busch Campus \\ 
 Rutgers, The State University of New Jersey \\
 110 Frelinghuysen Road \\
 Piscataway, NJ 08854-8019 USA}
\email{shelah@math.huji.ac.il}
\urladdr{http://shelah.logic.at}
\thanks{We thank Todd Eisworth and Tanmay Inamda for their helpful comments. The author thanks Alice Leonhardt for the beautiful typing up to 2019. The author would like to thank the ISF (Israel Science Foundation) for grant 1838(19) and earlier grants by the BSF, ISF, and the NSF.
The author thanks an individual who wishes to remain anonymous for generously funding typing services, and thanks Matt Grimes for the beautiful and careful typing.
  The reader should note that the version in my website is usually
  more updated than the one in the mathematical archive.  References
  like \cite[Th0.2=Ly5]{Sh:950} means the label of Th.0.2 is y5. 
 First Typed at Rutgers - 97/Mar/17}
 


\makeatletter
\@namedef{subjclassname@2020}{\textup{2020} Mathematics Subject Classification}
\makeatother
\subjclass[2020]{Primary: 03E04}
\keywords{set theory, pcf, cardinal arithmetic}
\date{September 23, 2022}
\begin{abstract}
\underline{Part A}:  A revised version of the guide in \cite{Sh:g}, with
corrections and expanded to include later works. 

\noindent
\underline{Part B}:  Corrections to \cite{Sh:g}. 

\noindent
\underline{Part C}:  Contains some revised proof and improved
theorems. 

\noindent
\underline{Part D}:  Contains a list of relevant references.  

\underline{Recent} (July 2022) \underline{additions}
\begin{itemize}
    \item \S14 = \ref{12a.2} on: no choice
    
    
    \item On inner models \ref{12.29}, see \cite{Sh:805} 
    
    \item On Black Boxes and abelian groups \ref{14.44} - \ref{14.53}, see \cite{Sh:750} and \cite{Sh:898}
    
    \item On somewhat free scales \ref{2.16}, see \cite{Sh:1008}, 
    
    \item On $n$-dimensional Black Boxes, quite free abelian groups such that $\Hom(G,\bbZ) = \{0\}$, \ref{14.56}, see \cite{Sh:1028}, 
    
    \item Survey on the existence of universal models; in particular, abelian groups \ref{14.32}, see \cite{Sh:1151}.
\end{itemize}

\end{abstract}

\maketitle
\numberwithin{equation}{section}
\setcounter{section}{-1}

\newpage

\centerline {Contents}

\noindent
Part A: Analytic Guide, \hfill pg.\pageref{A} 
\smallskip

\noindent
0. \quad $\check I[\lambda]$ and partial squares, \hfill pg.\pageref{A0}
\smallskip

\noindent
1. \quad Guessing clubs, \hfill pg.\pageref{A1}
\smallskip

\noindent
2. \quad Existence of l.u.b.,\hfill pg.\pageref{A2}
\smallskip

\noindent
3. \quad Uncountable cofinality,\hfill pg.\pageref{A3} 
\smallskip

\noindent
4. \quad Products, $T_D(f),\bfU(f)$,\hfill pg.\pageref{A4}
\smallskip

\noindent
5. \quad pcf Theory,\hfill pg.\pageref{A5}
\smallskip

\noindent
6. \quad Representation and pp,\hfill pg.\pageref{A6} 
\smallskip

\noindent
7. \quad Covering numbers,\hfill pg.\pageref{A7} 
\smallskip

\noindent
8. \quad Bounds in cardinal arithmetic,\hfill pg.\pageref{A8} 
\smallskip

\noindent
9. \quad J\'onsson algebras,\hfill pg.\pageref{A9} 
\smallskip

\noindent
10. \quad Coloring = negative partition relations,\hfill pg.\pageref{A10}
\smallskip

\noindent
11. \quad Trees and linear orders, \hfill pg.\pageref{A11}
\smallskip

\noindent
12. \quad Boolean Algebras and general topology,\hfill pg.\pageref{A12}
\smallskip

\noindent
13. \quad Strong covering, Forcing, and Partition Calculus,\hfill pg.\pageref{A13}
\smallskip

\noindent
14.  \quad Axiom of Choice, weak versions,\hfill pg.\pageref{A14}
\smallskip

\noindent
15.  \quad Transversals and freeness,\hfill pg.\pageref{A15}
\smallskip

\noindent
16. \quad Model Theory, Algebra, and Black Boxes\hfill pg.\pageref{A16} 
\smallskip

\noindent
17. \quad Discussion,\hfill pg.\pageref{A17}
\smallskip


\noindent
Part B: Corrections to \cite{Sh:g},\hfill pg.\pageref{B} 
\bigskip

\noindent
Part C: Details on changes and additions,\hfill pg.\pageref{C}
\smallskip

\noindent
17. \quad Short expansions, (label k),\hfill pg.\pageref{C1} 
\smallskip

\noindent
18. \quad More on II, \S(3.5), \ref{3.5}:\\ \indent \quad \quad Trees with
$\kappa$-branches ordered by $<_{J^{\bd}_\kappa}$, (label k),\hfill pg.\pageref{C2}
\smallskip

\noindent
19. \quad More on III, \S(4.10), \ref{4.10}:\\ \indent \quad \quad  Densely running away from
colours (label ac),\hfill pg.\pageref{C3}
\smallskip

\noindent
20. \quad Guessing clubs by countable $C$'s, (label 19),\hfill pg.\pageref{C4} 

\newpage

\begin{notation}:  $=^+$ appears in the 
following context: $\mu =^+ \sup\{\ldots\}$ means 
``both sides are equal, and if in the right side the $\sup$ is not
obtained, then it is singular."

For a set $C$ of ordinals, $\acc(C) = \{\alpha \in C:\alpha = 
\sup(\alpha \cap C)\}$, $\nacc(C) = C \setminus \acc(C)$.

The aim of this guide is to help the reader find out what is said in 
\cite{Sh:g} and related works of the author, what are the theorems and
definitions or where to look for them.

Let $[A]^\kappa = \{a \subseteq A:|a| = \kappa\}$, similarly
$[A]^{< \kappa}$ and $[A]^{\le \kappa}$. We denote 
$[A]^{\le \kappa}$ also as ${\clP}_{\le \kappa}[A]$.
\end{notation}

\newpage
 \section*{Part A - Analytic guide} \label{A}

\section {$I[\lambda]$ and partial squares} \label{A0}
\bigskip

See \cite{Sh:108}, \cite{Sh:88a}, \cite[2.3(5)]{Sh:345a}, equivalent forms
\cite[1.2]{Sh:420}, preservation of stationary subsets by $\mu$-complete
forcing \cite[21]{Sh:108}, \cite[10]{Sh:88a}.

\begin{definition}
\label{0.1}  
Let $\lambda = \cf(\lambda) > \aleph_0$.  For $S \subseteq \lambda$
we have: $S \in \check I[\lambda]$ \Iff \, for some club $E$ of $\lambda$ and
$\langle C_\alpha:\alpha < \lambda \rangle$ we have: $C_\alpha$ is a closed
subset of $\alpha$, otp$(C_\alpha) < \alpha$,

\[
[\beta \in \nacc(C_\alpha) \Rightarrow C_\beta = \beta \cap
C_\alpha] \text{ and }
[\alpha \in E \cap S \Rightarrow \alpha = \sup(C_\alpha)]
\]

\mn
(and every $\beta \in \nacc(C_\alpha)$ is a successor ordinal); note
$E \cap S$ has no inaccessible cardinal as a member.  Note that
\cite[1.2]{Sh:420} says that the definition just given is equivalent to 
those used in \cite{Sh:108}, \cite{Sh:88a}.

We can demand further $\alpha \in E \cap S \Rightarrow \otp(C_\alpha)
= \cf(\alpha)$.  But we can demand less: for each $\alpha$ we are
given $< \lambda$ candidates for $C_\alpha$, and for $C$ a candidate for
$\alpha$ and $\beta < \alpha,C \cap (\beta + 1)$ is a candidate for some
$\gamma < \alpha$.  $\check I[\lambda]$ is a 
normal ideal, and in many cases of the form ``non-stationary ideal 
$+ S$" (see \cite{Sh:108}; \cite{Sh:88a}).
\end{definition}

\begin{remark}
\label{0.1A} 
$\check I[\lambda]$ is a normal ideal but many times it has the form 
$\{A \subseteq \lambda:A \cap S \text{ non-stationary}\}$ and then $S$ is the
``bad" set of $\lambda$.  This holds for 
$\check I[\lambda] \restriction \{ \delta < \lambda:\cf(\delta) 
= \kappa\}$ if $\lambda = \lambda^{< \kappa}$ or less 
(see \cite{Sh:108}, \cite{Sh:88a}).
\end{remark}

\begin{claim}
\label{0.2}
1) If $\lambda$ is regular, then $S = S^{\lambda^+}_{< \lambda} = 
\{\delta < \lambda^+:\cf(\delta) < \lambda\}$ is the union of $\lambda$
sets on each of which we have a square (see below) hence belongs to 
$\check I[\lambda]$, see \cite[4.1]{Sh:351}.

\noindent
2) If $\lambda = \lambda^{< \kappa}$, then $\{\delta < \lambda^+:
\cf(\delta) < \kappa\}$ is the union of $\lambda$ sets on
each of which we have a square (see \cite{Sh:237e}), hence the set belongs
to $\check I[\lambda]$.

\noindent
3) Moreover, if $\lambda > \aleph_0$ is regular and 
$\alpha < \lambda \Rightarrow |\alpha|^{< \kappa} < \lambda$ or just
$\alpha < \lambda \Rightarrow \cov(|\alpha|,\kappa,\kappa,2) <
\lambda$ then $\{\delta < \lambda:\cf(\delta) < \kappa\} \in 
\check I[\lambda]$ (see \cite[2.8]{Sh:420}, the case $\kappa =
\aleph_0$ is trivial).   By Dzamonja, Shelah \cite{Sh:562} 
the same assumption gives $\{\delta < \lambda^+:\cf(\delta) <
\cf(\lambda)\}$ is the union of $\le \lambda$ sets on each of which
we have square.  Also in \cite{Sh:562} there are results on 
getting squares with $\lambda$ singular and results 
with an inaccessible instead of $\lambda^+$.
\end{claim}

\begin{definition}
\label{0.2A}
$S \subseteq \mu$ has a square if we
have $S^+,S \subseteq S^+ \subseteq \mu$ and $\langle C_\alpha:
\alpha \in S^+ \rangle$ such that: $C_\alpha$ is a closed subset of 
$\alpha$ of order type $< \alpha$, and $\alpha \in C_\beta \Rightarrow 
C_\alpha = \alpha \cap C_\beta$ and [$\alpha$ is a limit ordinal iff 
$\alpha = \sup(C_\alpha)]$ for $\alpha \in S$; also if $\alpha 
\in S \Rightarrow \cf(\alpha) \le \kappa(< \kappa)$, we can add 
``$\otp(C_\alpha) \le \kappa(< \kappa)$".
\end{definition}

\begin{remark}
\label{0.3}
Related ideals \cite[2.3,2.4]{Sh:345a}
\cite[2.3,2.4,2.5,5.1,5.1A,5.2]{Sh:371}.
\end{remark}

\begin{claim}
\label{0.4}
If $\kappa^+ < \lambda = \cf(\lambda)$, then we 
can find a stationary

\[
S \subseteq \{\delta < \lambda:\cf(\delta) = \cf(\kappa)\},\ 
S \in \check I[\lambda]
\]

\mn
\cite[1.5]{Sh:420} (somewhat more \cite[1.4]{Sh:420}).
\end{claim}

\begin{note}
\label{0.5}
Negative consistency results:
\cite{Sh:108}, (``GCH + the bad set for $\aleph_{\omega +1}$ is stationary")
Magidor, Shelah \cite{Sh:204}, Hajnal, Juhasz, Shelah \cite{Sh:249}, 
consistency of $\check I[\lambda]$ large but 
stationary sets reflect \cite{Sh:351}.
\end{note}

\begin{note}
\label{0.6}
On killing stationary sets by forcing \cite{Sh:108}, 
\cite[18,19]{Sh:88a}, \cite[2.4]{Sh:371}.
\end{note}

\begin{note}
\label{0.7}
On consequences of pcf structure (\cite{Sh:108},
\cite[Ch.VIII,\S5?]{Sh:g}, \cite[5.17,5.18]{Sh:589}), e.g.
(GCH) the bad stationary subsets of $\aleph_{\omega+1}$ do 
not reflect (\cite{Sh:108} or \cite{Sh:88a}).
\end{note}
\newpage

\section {Guessing clubs} \label{A1}

\begin{definition}
\label{1.1}
Definition of ideals \cite[1.3,1.5,3.1]{Sh:365}: definition of $g \ell$ 
\cite[2.1]{Sh:365}: also \cite[1.8]{Sh:380}.
\end{definition}

\noindent
For example
\begin{definition}
For $\bar C = \langle C_\delta:\delta \in S \rangle$,
$S \subseteq \lambda = \cf(\lambda) > \aleph_0$, $C_\delta$ a club of $\delta$:

\[
\id^b(\bar C) = \{A \subseteq \lambda:\,\text{for some club } E \text{ of }
\lambda, \text{ for no } \delta \in S \cap A \cap E, \text{ is }
C_\delta \subseteq E \}
\]

\begin{equation*}
\begin{array}{clcr}
\id^a(\bar C) = \{A \subseteq \lambda:\,&\text{for some club } E 
\text{ of } \lambda, \text{ for no } \delta \in S \cap A \cap E, \\
  &\,\text{ is } \sup(C_\delta \setminus E) < \sup C_\delta \}
\end{array}
\end{equation*}

\begin{equation*}
\begin{array}{clcr}
\id_p(\bar C) = \big\{A \subseteq \lambda:\,&\text{for some club } E \text{ of }
\lambda, \text{ for no } \delta \in S \cap A \cap E, \\
  &\,\text{is } \delta = \sup\big(E \cap \text{ nacc}(C_\delta) \big)\big\}.
\end{array}
\end{equation*}
\end{definition}

\begin{note}
\label{1.2}
1) Easy facts \cite[1.4,1.6]{Sh:365}.

\noindent
2) For $\lambda,S \subseteq \lambda$ stationary, concerning 
the existence of $\bar C = \langle C_\delta:\delta \in
S \rangle$ ``guessing clubs of $\lambda$" \cite[\S2]{Sh:365} (and
\cite[Ch.III,7.8A-G]{Sh:e}).
\end{note}

\begin{claim}
\label{1.3}
The following items give sufficient conditions for the properness
of the above ideals for $\lambda$ regular uncountable:
\mn
\begin{enumerate}
\item[$(a)$]   If $\delta \in S \Rightarrow \cf(\delta) < \mu$
for some $\mu < \lambda$, \underline{then}  we can find clubs $C_\delta$ for
$\delta \in S$ such that $\id^b(\langle C_\delta:\delta \in S \rangle)$ is
a proper ideal (i.e. for every club $E$ of $\lambda$ for some $\delta,
C_\delta \subseteq E$) by \cite[2.3(2)]{Sh:365}.
\sn
\item[$(b)$]   If $\lambda = \mu^+,\mu$ regular,
\then
\begin{enumerate}
    \item[$(\alpha)$]
There is a sequence
$\langle (S_\eps, \bar C_\eps) : \eps < \lambda \rangle $
such that:
\begin{enumerate}[$\bullet_1$]
    \item each $ S_\eps$ is a subset of $\lambda$ and 
    $$\bigcup \{ S_\eps : \eps< \lambda^+ \} = \{\delta < \lambda ^+ : 
    \cf(\delta) < \lambda \}$$
    
    \item for each $\eps < \lambda$, $ \bar C_\eps$
    has the form $\langle C_{\eps, \alpha} : \alpha \in S_\eps\rangle$
    and is a partial square, which means that:
\sn
    $C_{\eps, \alpha} $ is a closed subset of $\alpha$,
    is unbounded if $ \alpha $ is a limit ordinal, is included in $ S_\eps$
    and $$\beta \in C_{\eps, \alpha}\ \Rightarrow\ C_{\eps, \beta} = 
    C_{\eps, \alpha } \cap \beta$$
    
    \item also, $C_{\eps, \alpha}$ is of cardinality $< \lambda$
\end{enumerate}

\sn
    \item[$(\beta)$]
    For every limit ordinal $\delta(*) < \lambda$, for some partial square 
    $ \bar C = \langle C_ \alpha : \alpha \in S_* \rangle$ of $\lambda$ we have 
    $\alpha \in S_* \Rightarrow {\rm otp}(C_ \alpha)\le \delta (*) $, and letting
    $$S^*= \{ \alpha \in S_* : \otp ( \alpha ) = \delta(*)\}$$ 
    we have that $\bar{C} \rest S^*$ guesses clubs;
    that is, $ \varnothing \notin \id^a(\bar{ C })$.
    
    \sn

\item[$(\gamma)$] 
For every limit ordinal $\delta(*) < \lambda^+$, for some club $E$ of $\lambda^+$,
an ordinal $\eps < \lambda$ and and limit ordinal $\Upsilon < \lambda$
divisible by $\delta(*)$ with the same cofinality, the sequence 
$\bar{C}_{\eps, E}$ is as required in $(\beta)$ replacing
$\delta (*)$ by $\Upsilon$ where 
(of course $\bar{C}_\eps$ is from sub-clause $(\alpha)$)
$$\bar{C}_{\eps, E} = \langle C_{\eps, \alpha } \cap E :
\alpha \in S_\eps \cap \alpha \rangle$$
\underline{Why?} clause $ ( \beta )$ follows
from clause $ ( \gamma ) $, clause $ ( \gamma ) $
follows from clause $ ( \alpha ) $  as usual
(trying enough times). 

  For clause $ ( \alpha ) $ see \cite[2.14(2)]{Sh:365} but the proof there is inaccurate, see \cite[4.4,pg.47]{Sh:351} or see part (B), \pageref{pg136}, see part B here.
\end{enumerate} 

\sn
\item[$(c)$]   If $\lambda = \mu^+$, $\mu$ regular, $S \subseteq
S^\lambda_\mu$ is stationary, \underline{then}  we can find 
$\bar C = \langle C_\delta : \delta \in S \rangle$, $C_\delta$ a club
of $\delta$, $\otp(C_\delta) = \mu$, $[\alpha \in \nacc(C_\delta)
\Rightarrow \cf(\alpha) = \mu]$ and $\id_p(\bar C)$ a proper ideal (i.e.
for every club $E$ of $\lambda$ for some $\delta$, 
$\delta = \sup(E \cap \nacc(C_\delta))$), \cite[2.3(1)]{Sh:365}, \cite{Sh:413}, 
\cite[\S3]{Sh:572}.
\sn
\item[$(d)$]   If $[\lambda = \mu^+$, $\mu$ singular, and 
$\delta \in S \Rightarrow \cf(\delta) = \cf(\mu) > \aleph_0]$ \underline{or}
[$\lambda$ inaccessible and $\delta \in S \Rightarrow \cf(\delta) \in (\aleph_0,\delta)$], \underline{then} for some $\bar C = \LL C_\alpha : \alpha \in S \RR$ we have: 
$\id^a(\bar C)$ is proper and for each $\delta \in S$ we
have: $\langle \cf(\alpha) : \alpha \in \nacc(C_\delta) \rangle$
converges to $|\delta|$ (and is strictly increasing) \cite[2.6,2.7]{Sh:365}.
\sn
\item[$(e)$]  If $S^* \subseteq \lambda$ is stationary and does not reflect
outside itself and $S \subseteq \lambda$ is stationary, \underline{then} for some
$\bar C = \langle C_\delta : \delta \in S \rangle$ we have 
$\nacc(C_\delta) \subseteq S^*$, and $\id_p(\bar C)$ is a proper ideal, \cite[2.13]{Sh:365}.
\sn
\item[$(f)$]   Similar theorems with ideals \cite[1.7,2.4]{Sh:380},
\cite[1.11,1.12]{Sh:413} other related ideals \cite[1.10]{Sh:380}.
\sn
\item[$(g)$]  More in the places above and \cite[2.6,2.8,2.9]{Sh:413} and
\cite{Sh:449}.
\sn
\item[$(h)$]   Assume $\lambda = \cf(\lambda)$, $S \subseteq \{\delta < 
\lambda^+:\cf(\delta) = \lambda\}$ is stationary and $\chi$
satisfies one of the following:  $\lambda = \chi^+$ \underline{or} $\chi =
\min\{\tau < \lambda:(\exists \theta \le \tau)\ \tau^\theta \ge 
\lambda\}$ \underline{or} $\lambda$ strongly inaccessible not Mahlo.
\underline{then}  we can find $\langle C_\delta,h_\delta:\delta \in S \rangle$ such
that: $C_\delta = \{\alpha_{\delta,\zeta}:\zeta < \lambda\}$ is a club of
$\delta,\alpha_{\delta,\zeta}$ increasing with $\zeta$, 
$h_\delta : C_\delta \rightarrow \chi$ and for every club $E$ of $\lambda^+$, 
for stationarily many $\delta \in S$, for each $i < \chi$,
\[
\{\zeta < \lambda:\alpha_{\delta,\zeta} \in E,\ \alpha_{\delta,\zeta +1} \in E, 
\text{ and } h_\delta(\alpha_{\delta,\zeta}) = i\}
\]
is a stationary subset of $\lambda$ (see \cite[\S3]{Sh:413}, 
\cite[\S3]{Sh:572}).
If $\lambda$ is a limit of inaccessibles, we can demand cf$(\alpha_{\delta,
\zeta+1}) > \zeta$.
\sn
\item[$(i)$]   If $\lambda,\bar C = \langle C_\delta:\delta \in S^+ \rangle$
is as in \ref{0.1}, $S \subseteq S^+, \sup\{|C_\alpha|^+:\alpha \in S\}
< \lambda$ then for some club $E$ of $\lambda$, 
$\bar C' = \langle g \ell(C_\delta,E):\delta \in S^+ \cap \acc(E) \rangle$ 
is as in \ref{0.1} and for every club $E_1 \subseteq E$ of $\lambda$, for 
stationarily many $\delta \in S$, we have
$\alpha \in C'_\delta \Rightarrow \sup(C'_\delta \cap \alpha) 
\le \sup(E \cap \alpha)$.
\sn
\item[$(j)$]  Assume $\lambda = \cf(\lambda)$ and $S 
\subseteq S^{\lambda^+}_\lambda = \{\delta < \lambda^+:\cf(\delta) 
= \lambda\}$ is stationary.
Then we can find an $S$-club system $\bar C = \langle C_\delta : \delta \in S \rangle$ 
and $h : S \rightarrow \lambda$ such that for any club $E$ of 
$\lambda^+$, for stationarily many $\delta \in S$, for every $i < \lambda$,
the set $\nacc(C_\delta) \cap h^{-1}(\{i\})$ is unbounded in $\delta$
(under reasonable assumption 
$\big|\{C_\delta \cap \alpha:\alpha \in \nacc(C_\delta)\} \big| \le \lambda$), 
see \cite[3.3]{Sh:413}.
\end{enumerate}
\end{claim}

\begin{note}
\label{1.4}
On $\otimes_{\bar C}$, $\otimes^\kappa_{\bar C}$ for some $S$-club system \cite[2.12,2.12A,4.10]{Sh:365} and a
colouring theorem \cite[4.9]{Sh:365} (see earlier \cite{Sh:276}).
Where $\lambda$ is a Mahlo cardinal,
\mn
\begin{enumerate}
\item[$\otimes_{\bar C}$]   $\bar C$ has the form $\LL C_\delta : \delta \in S \RR$, $S \subseteq \lambda$ a set of inaccessibles, $C_\delta$
a club of $\delta$ such that: for every club $E$ of $\lambda$ for stationary
many $\delta \in S$, $E \cap \delta \setminus C_\delta$ is unbounded in
$\delta$
\end{enumerate}
\mn
and for $\kappa < \lambda$:
\mn
\begin{enumerate}
\item[$\otimes^\kappa_{\bar C}$]   $\bar C$ has the form 
$\langle C_\delta : \delta \in S^\lambda_{\in} \rangle$, 
$S^\lambda_{\in} = \{ \mu < \lambda:\mu$ inaccessible$\}$, such that: for every club $E$ of
$\lambda$, for stationarily many $\delta \in S^\lambda_{\in} \cap 
\acc(E)$, for no $\zeta < \kappa$ and $\alpha_\eps \in 
S^\lambda_{\in} (\eps < \zeta)$ is $\nacc(E) \cap \delta \setminus 
\bigcup\limits_{\eps < \zeta} 
C_{\alpha_\eps}$ bounded in $\delta$.
\end{enumerate}
\mn
By \cite[4.9]{Sh:365} if $\kappa$ is a Mahlo cardinal and $\otimes^\kappa
_{\bar C}$, \underline{then}  for some 2-place function $c$ from $\kappa$ to
$\omega$, for every pairwise disjoint $w_i \subseteq \kappa$, $|w_i| < \kappa$ 
for $i < \kappa$, and $n$, for some $i < j$, $\Rang(c \rest w_i \times
w_j) \subseteq (n,\omega)$.  By \cite[4.10B]{Sh:365}, $\otimes^2_{\bar C} 
\Leftrightarrow \otimes^{\aleph_0}_{\bar C}$, also $\otimes^2_{\bar C}$ is a 
strengthened form of ``$\kappa$ not weakly compact", which fails under 
mild conditions (\cite[4.10A]{Sh:365}).
See more in \cite[4.13]{Sh:365}.
\end{note}

\begin{note}
\label{1.5}:  id$_p(\bar C,\bar I)$ is decomposable
\cite[3.2,3.3]{Sh:365}.
\end{note}

\begin{note}
\label{1.6}
If $\kappa^+ < \lambda$, we can find $\langle 
{\cP}_\alpha:\alpha < \lambda \rangle$ such that: 
\begin{enumerate} 
\item[(s)] ${\cP}_\alpha$ is a family of $< \lambda$ closed subsets of $\alpha$,
\item[(b)]
$ \beta \in \nacc(C) \text{ and } C \in {\cP}_\alpha \Rightarrow C \cap
\beta \in {\cP}_\beta $ 
\item[(c)]
for every club $E$ of $\lambda$ for stationarily many $\alpha < \lambda$,
there is $C \in {\cP}_\alpha$, $\kappa = \otp(C)$, $\alpha = \sup(C)$
and $C \subseteq E$ \cite[1.3]{Sh:420} (we can replace $\kappa$ by $\delta
(*)$, $|\delta(*)| = \kappa$).
\end{enumerate}

\end{note}

\begin{note}
\label{1.7}
More on \ref{1.3}(c) in \cite[\S3]{Sh:413} and better in \cite[\S3]{Sh:572}.
\end{note}

\begin{note}
\label{1.8}
If we want to preserve $$\alpha \in \nacc(C_\alpha) 
\cap \nacc(C_\beta) \cap \nacc(C_\gamma) \Rightarrow C_\beta 
\cap \alpha = C_\gamma \cap \alpha$$ we can weaken the
guessing to: $\forall$ clubs $E$, $\exists^{\stat} \delta$ such that $E$
is not disjoint to any interval of $C_\alpha$.  See the proof of
\cite[6.2]{Sh:430}, \cite{Sh:562}.
\end{note}

\begin{note}
\label{1.9}
On ideals related to J\'onsson algebras and guessing clubs:
\cite{Sh:380}, \cite[\S1]{Sh:413} (used in \S8 here).
\end{note}
\newpage

\section {Existence of lub} \label{A2}
\bigskip

We discuss here lub of $\bar f = \langle f_\alpha : \alpha < \delta
\rangle \mod I$, where $f_\alpha \in {}^\kappa\Ord$, $I$ an ideal
on $\kappa$, $\kappa^+ < \cf(\delta)$.  See \cite{Sh:68}, \cite{Sh:111}, 
\cite[\S14]{Sh:282} and better \cite[\S1]{Sh:355}.

\begin{definition}
\label{2.1}
We say ``$f$ is a lub of $\langle f_\alpha:
\alpha < \delta \rangle \mod I$" where $I$ is an ideal on 
$\Dom(I),f_\alpha: \Dom(I) \rightarrow$ ordinals, if 
$\bigwedge\limits_{\alpha < \delta} f_\alpha \le_I f$, and 
$$\bigwedge\limits_{\alpha < \delta} f_\alpha \le f' \Rightarrow f 
\le f' \mod I.$$

We say ``$f$ is an eub (exact upper bound) of 
$\langle f_\alpha : \alpha < \delta \rangle \mod I$" where $I$ 
is an ideal on $\Dom(I)$, $f_\alpha:\Dom(I) \rightarrow$ ordinals, if 
$\bigwedge\limits_{\alpha < \delta} f_\alpha \le_I f$ and if
$g <_I \max\{f,1\}$ \underline{then}  for some $\alpha < \delta$ we have 
$g \le_I f_\alpha$ (see \cite[1.4(4)]{Sh:345a}); 
usually $\alpha < \beta \Rightarrow f_\alpha \le_I 
f_\beta$; ``$f$ is an eub of $\langle f_\alpha:\alpha < \delta
\rangle$ mod $I$" says more than ``$f$ is a lub of $\langle 
f_\alpha:\alpha < \delta \rangle \mod I$".
\end{definition}

\begin{claim}
\label{2.2}
The trichotomy theorem on the existence
of eub \cite[1.2,1.6]{Sh:355} (slightly more \cite[6.1]{Sh:430}, on
eub $\ne$ lub, see example \cite[6.1A]{Sh:430}). 

For example for $I$ a maximal ideal on $\kappa$, $f_\alpha \in {}^\kappa\Ord$ 
for $\alpha < \delta$, $\cf(\delta) > \kappa^+$, 
$\bar f = \langle f_\alpha / I : \alpha < \delta \rangle$ increasing, 
\underline{either} $\bar f$ has a $<_I$-eub, \underline{or} 
for some sequence $\bar w = \langle w_i : i < \kappa \rangle$ 
of sets of ordinals, $|w_i| \le \kappa$ we have:

\[
\bigwedge\limits_{\alpha < \delta} \, \bigvee\limits_{\beta < \delta} \,
\Big(\exists g \in \prod\limits_{i < \kappa} w_i \Big)
\big[f_\alpha/I < g/I < f_\beta/I \big].
\]

\mn
The $\cf(\delta) > \kappa^+$ is necessary by \cite{Sh:673}.
\end{claim}

\begin{definition}
\label{2.3}
\cite[2.6]{Sh:345a}.  We define:

\begin{equation*}
\begin{array}{clcr}
\ged_I(\bar f) =: \big\{\alpha < \delta:\,&\cf(\alpha) > \kappa,
\text{ and there is an unbounded } A \subseteq \alpha \\
  &\text{ and members } s_i \text{ of } 
I \text{ for } i \in A \text{ such that:} \\
  &\,i \in A \text{ and } j \in A \text{ and } i < j \text{ and } \\
&\zeta \in \kappa \setminus (s_i \cup s_j) \Rightarrow 
f_i(\zeta) \le f_j(\zeta)\big\}.
\end{array}
\end{equation*}

\mn
Sufficient conditions for the existence of eub \cite[1.7]{Sh:355} is
that $\ged_I(\bar f)$ is a stationary subset of $\delta$.  
\end{definition}

\begin{definition}
\label{2.4}
Let ${\ga}$ be a set of regular cardinals and
$N \prec (\clH(\lambda),\in)$: 
we define $\Ch^{\ga}_N(\theta) = \sup(N \cap \theta)$
for $\theta \in {\ga}$ \cite[3.5]{Sh:345a}, 
\cite[3.4(stationary)]{Sh:355}, \cite[1.2,1.3,1.4]{Sh:371} 
more \cite[3.3A,5.1A]{Sh:400} and \cite[\S6]{Sh:430}.
\end{definition}

\begin{claim}
\label{2.5}
On the good/bad/chaotic division.  
For $\bar f$ a $<_I$-increasing sequence
of functions from $\kappa$ to ordinals, we have a natural division of
$\{\delta < \ell g(\bar f),\cf(\delta) > \kappa^+\}$ to there:
\mn
\begin{enumerate}
\item[$\bullet_1$]   to $\ged_I(\bar f)$ (see \ref{2.3} above),
\sn
\item[$\bullet_2$]  $\ch(\bar f) = 
\{\delta < \ell g(\bar f) : \text{for some ultrafilter }
D \text{ on } \ell g(\bar f) \text{ disjoint to } I \text{ and }
w_i \subseteq \text{ ordinals for } i \in \Dom(I),\ |w_i| \le |\Dom(I)| 
\text{ and } \bigwedge\limits_{i < \delta} \, \bigvee\limits_{j < \delta}
 \, (\exists g \in \prod w_i)[f_i \le_D g \le_D f_j]\}$ and
\sn
\item[$\bullet_3$] $\bd_I(\bar f) = \ell g(\bar f) \setminus (\ged_I(\bar f)
\cup \ch_I(\bar f))$.
\end{enumerate}
\end{claim}

\begin{note}
\label{2.5d}
In \ref{2.5}:
\mn
\begin{enumerate}  
\item[(a)]   for every $\delta < \ell g(\bar f)$ of uncountable
cofinality there is a club $C$ of $\delta$ such that $\delta 
\in \ged_I (\bar f) \wedge \alpha \in C \wedge \cf(\alpha) > 
\kappa \Rightarrow \alpha \in \ged_I(\bar f)$ and $\delta \in 
\ch_I(\bar f) \Rightarrow C \subseteq \ch_I(\bar f)$;
\sn
\item[(b)]   for bd$_I(\bar f)$ to be non-trivial, 
$\ell g(\bar f)$ should not be so small among the alephs. 
\end{enumerate}
\mn
There are connections to NPT (see \S12) and 
$\check I[\cf(\ell g(\bar f))]$ (see \S1) (and consistency of the 
existence of counterexamples; see
\cite{Sh:108}, \cite{Sh:204}, \cite[1.6]{Sh:355}, \cite{Sh:523}).
\end{note}

\begin{problem}
\label{2.6}
Is the following consistent: $\{\delta < \aleph_{\omega + 1}:
\cf(\delta) = \aleph_2\} \notin \check I[\aleph_{\omega + 1}]$ or 
$2^{\aleph_0} < \aleph_\omega$ and $\{\delta <
\aleph_{\omega +1}:\text{cf}(\delta) = (2^{\aleph_0})^+\} \notin
\check I[\aleph_{\omega +1}]$ (also for inaccessibles) 
or $\bar f = \langle f_\alpha:
\alpha < \aleph_{\omega +1} \rangle$, $f_\alpha \in 
\prod\limits_{n < \omega} \aleph_n$, $\ch_{J^{\bd}_\omega}(\bar f) \cap 
\{\delta < \aleph_{\omega +1}: \cf(\delta) = \aleph_2\}$ 
stationary or $(\forall S)[S \in \check I[\aleph_2] \text{ and } 
\bigwedge\limits_{\delta \in S} \cf(\delta) = \aleph_1 \Rightarrow S$
not stationary]?
\end{problem}

\begin{note}
\label{2.7}
More on \S2, see in \S12 (in universes without full choice).
\end{note}

\begin{note}
\label{2.8}
See more in \cite{Sh:506} for generalization to the case
$\cf(\delta) \le |\Dom I|$.  
On existence of eub see \cite[3.10]{Sh:506} and \cite[6.4]{Sh:589}.
\end{note}

\begin{claim}
\label{2.9}
Assume $\lambda = \cf(\lambda) \ge \mu > 2^\kappa,
f_\alpha \in {}^\kappa\Ord$ for $\alpha < \kappa$.  \underline{then}  for some
$\beta^*_i \, (i < \kappa)$ and $w \subseteq \kappa$ we have: $i \in w
\Rightarrow \cf(\beta^*_i) > 2^\kappa$ and for every $f \in
\prod\limits_{i \in w} \beta^*_i$ for unboundedly many $\alpha 
< \lambda$ we have $i \in w \Rightarrow f(i) < f_\alpha(i) < 
\beta^*_i$ and $i \in \kappa \setminus w \Rightarrow f_\alpha(i) 
= \beta^*_i$; \cite[6.6D]{Sh:430} (slightly more general); more 
detailed proof \cite[6.1]{Sh:513}, more variants \cite[\S7]{Sh:620}.
\end{claim}

\begin{note}
\label{2.10}
On decreasing sequences see \cite[6.1,6.2]{Sh:589}.
See also \cite{Sh:829}.
\end{note}

\begin{claim}
\label{2.13}
The restriction in \ref{2.5} to $\cf(\delta) > \kappa$ is necessary
(by Kojman-Shelah \cite{Sh:673}.
\end{claim}

\begin{claim}
\label{2.16}
See \cite{Sh:1008} for more; e.g. even if 
$\pp(\aleph_\omega) = \aleph_{\omega+1}$ then there is a $<_{J^\bd_\omega}$-increasing sequence $\LL f_\alpha : \alpha < \aleph_{\omega+1} \RR$ of members of $\prod\limits_n \omega_n$ which is $(\aleph_\omega,\aleph_4)$-free; i.e. for every $\alpha < \aleph_\omega^+$ there is an equivalence relation $E_\alpha$ on each class of cardinality $< \aleph_4$ and $h : \alpha \to \omega$ such that if $\beta,\gamma < \alpha$ are not $E_\alpha$-equivalent then 
$$\big\{ f_\beta(n) : n \geq h(\beta) \big\} \cap \big\{ f_\gamma(n) : n \geq h(\gamma) \big\} = \varnothing$$ 
(See more in \cite{Sh:1008}, \cite{Sh:1028}.)
\end{claim}
\newpage

\section {Uncountable cofinality and $\aleph_1$-complete filters and 
products: \cite{Sh:71}, \cite{Sh:111}, \cite{Sh:256}} \label{A3}
\bigskip

\begin{note}
\label{3.1}
Assume $\langle \lambda_i:i \le \kappa
\rangle$ is an increasing continuous sequence of singulars, $\aleph_0 <
\kappa = \cf(\kappa) < \lambda_0$.  
Let $\lambda = \lambda_\kappa$.  If $\{i < \kappa:\pp(\lambda_i) 
= \lambda^+_i\}$ is a stationary subset of
$\kappa$, \underline{then}  $\pp(\lambda) = \lambda^+$, \cite[2.4(1)]{Sh:355}.

Moreover, $\pp(\lambda_\kappa)$ is bounded by
$\lambda^{+\|h\|}_\kappa$ where $\pp(\lambda_i) = 
\lambda^{+h(i)}_i$ hence we have a bound on $\pp(\lambda)$ in
many cases \cite[2.4]{Sh:355}, \cite[1.10]{Sh:371}.
\end{note}

\begin{note}
\label{3.2}
Definition of various ranks and niceness of filters in 
\cite[1.1,1.2,1.4,3.12]{Sh:386} (more generally on pair
$(t,D)$ or for $D \in \Fil(e,y)$ see \cite[\S5]{Sh:410} 
and \cite[\S3,\S4,\S5]{Sh:420}).  For $\kappa = \cf(\kappa) > 
\aleph_0$, $D$ a normal filter on $\kappa$ and $f \in {}^\kappa\Ord$ 
let $\rk^2(f,D)$ be $\le \alpha$ iff for every $A \in D^+$
and $g <_{D+A}f$ for some normal filter, $D_1 \supseteq D+A$ we have
$\rk^2(g,D_1) \le \beta$ for some $\beta < \alpha$. 
$D$ is nice if $f \in {}^{\Dom(D)}\Ord \Rightarrow \rk^2 (f,D) <
\infty$.
\end{note}

\begin{note}
\label{3.3}
If for any $A \subseteq 2^{\aleph_1}$ in
$K[A]$, there are Ramsey cardinals (or suitable Erd\"os cardinals 
which occurs if cardinal arithmetic is not trivial, essentially by 
Dodd and Jensen \cite{DJ1}), \underline{then}  every normal 
filter on $\omega_1$ is nice \cite[1.7,1.13]{Sh:386}; more in 
\cite[\S1]{Sh:386}, \cite[\S3,\S4,\S5]{Sh:420}.
\end{note}

\begin{note}
\label{3.4}
\cite[2.2,2.2A,2.4,2.7]{Sh:386}, \cite[\S4]{Sh:420} $A_e(f)$ \cite[3.3]{Sh:386}.
\end{note}

\begin{note}
\label{3.5}
Rank, basic properties: 
  
\cite[2.3,2.4,2.8,2.9,2.10,2.11,2.12,2.14,2.21,3.4,3.8]{Sh:386}.
\end{note}

\begin{note}
\label{3.6}
Rank, connection to forcing:

\cite[Definition 2.6 $(E^t_p)$,2.6A,2.7A]{Sh:386},\cite[\S3]{Sh:420}.
\end{note}

\begin{note}
\label{3.7}
Rank, relation with $T_D$

\cite[2.15,2.16,2.17,2.18,2.19,2.20,2.22]{Sh:386}.
\end{note}

\begin{note}
\label{3.8}
Ranks-going down: ranks when we divide
$\omega_1$, \cite[3.2]{Sh:386} each $f$ successor \cite[3.6]{Sh:386}, each
$f$ limit \cite[3.7]{Sh:386}.
\end{note}

\begin{note}
\label{3.9}
Rank, getting $\kappa$-like reduced products \cite[3.10,3.11,3.11A]{Sh:386}.
\end{note}

\begin{note}
\label{3.10}
Generic ultrapower with all $\kappa > 
\beth_2(\aleph_1)$ represented: \cite[1.3]{Sh:333}, just for one
\cite[1.4]{Sh:333} (earlier \cite{Sh:111}).
\end{note}

\begin{note}
\label{3.11}
Ranks are $< \infty$ \cite[3.13-18]{Sh:386}.
\end{note}

\begin{note}
\label{3.12}
Preservative pairs (see \ref{3.15}),
definition and basic properties \cite[4.15]{Sh:386}.
\end{note}

\begin{note}
\label{3.13}
Specific functions are preservative: 

\cite[4.6]{Sh:386} ($H_s =$ successor), \cite[5.8]{Sh:386} 
($H^{ia}=$ next inaccessible), \cite[5.9]{Sh:386} 
($H^{\epsilon-m}=$ next $\epsilon$-Mahlo).
\end{note}

\begin{note}
\label{3.14}
The class of preservative pairs is closed under:
\mn
\begin{enumerate}
\item  $H^*(i)$ iterating $H \, i$ times \cite[4.7,4.8,4.9]{Sh:386}
\sn
\item  composition \cite[4.10]{Sh:386}
\sn
\item  $\sup_{n < \omega} H^n$ \cite[4.11]{Sh:386}
\sn
\item  iterating $\alpha$ times, $\alpha < \omega_1$ \cite[4.12]{Sh:386}
\sn
\item  more \cite[4.13]{Sh:386}
\sn
\item  induction \cite[\S2]{Sh:333}.
\end{enumerate}
\end{note}

\begin{note}
\label{3.15}
Preservative pairs are bounds on
cardinal exponentiation \cite[5.1,5.2,5.3]{Sh:386}.
\end{note}

\begin{note}
\label{3.16}
If $\rk^2_E(f) = \rk^3_E(f) =
\lambda$ inaccessible, \underline{then}  modulo (fil $E$) almost every $f(i)$ is
inaccessible \cite[5.7]{Sh:386}.
\end{note}

\begin{note}
\label{3.17}
Generalizing normal filters and then ranks 
\cite[\S5]{Sh:410}, \cite[\S3,\S4,\S5]{Sh:420}.
\end{note}

\begin{note}
\label{3.881}
Combinatorial theorem using ranks, \cite{Sh:881}, if
$\lambda > \cf(\lambda) > \aleph_0$ and $2^{\cf(\lambda)} < \lambda$
then $\lambda \rightarrow (\lambda,\omega +1)^2$.
\end{note}

\begin{note}
\label{3.18}
For set theory with weak choice much remains (see \cite{Sh:497}, here \S12).
\end{note}
\newpage

\section {Products, $T_D(f),\mathbf U$} \label{A4}

We deal with computing $T_D(f),\mathbf U_D(f)$ and reduced products
$\prod\limits_{i < \kappa} f(i)/D$ from pcf, mainly when $(\forall i)[f(i) >
2^\kappa]$ see \cite[\S3]{Sh:506}, \cite[\S1]{Sh:589}, \cite[\S4]{Sh:589} on
$T_D$ earlier, Galvin Hajnal \cite{GH}.

\begin{definition}
\label{4.1}
1) Define $$T_D(f) =: \min\Big\{|{\cF}|:{\cF} \subseteq \prod\limits_{i} 
(f(i)+1) \text{ and } f \ne g \in {\cF} \Rightarrow f \ne_D g \Big\}$$ 
(i.e. $\{i:f(i) \ne g(i)\} \in D$) and $\cF$ is maximal with respect to 
those properties. 

$T_\Gamma(f) = \sup\{T_D(f):D \in \Gamma\}$ for $\Gamma$ set of filters on
Dom$(f)$, similarly for $\Gamma$ set of ideals and naturally $T_\Gamma
(\lambda)$. 

\noindent
2) 
\begin{align*}
    \bfU_D(f,< \theta) = \min \Big\{ |{\cA}|&\ :  {\cA} \subseteq
    \prod[f(i)]^{< \theta},\ A \in \cA \Rightarrow
    |A| < \theta,\\ 
    & \text{ such that for every } g \in {}^\kappa \Ord \text{ with } g <_D f, \\
    & \text{ for some } \bar A \in {\cA} \text{ we have } \{i < \kappa:g(i) \in A_i\} \ne \varnothing \mod D \Big\}
\end{align*}

If $\theta = \kappa^+$ we may omit it, (note: if $\cf(\theta) > \kappa$ we
can replace $\bar A$ by $\bigcup\limits_{i < \kappa} A_i$.  

[See more: \cite{Sh:430}, \cite{Sh:552}.
\end{definition}

\begin{note}
\label{4.2}
If $\lambda > 2^{< \theta},\theta \ge \sigma = \cf
(\sigma) > \aleph_0$ and $\Gamma = \Gamma(\theta,\sigma)$ (the set of
$\sigma$-complete ideals on a cardinal $< \theta$) we have

\[
T_\Gamma(\lambda) = \cov(\lambda,\theta,\theta,\sigma)
\]

\mn
(the latter can be computed from case of $\pp_\Gamma$); 
\cite[5.9,p.94]{Sh:355}. 
If $\theta^\kappa < \mathbf U_D(\lambda)$, \underline{then}  $T_D(f) = \mathbf U_D (f)$.
\end{note}

\begin{note}
\label{4.3}
A pcf characterization when $\lambda \le T_D(f)$ holds,
under $2^{\Dom(D)} < \min\limits_i f(i)$ and
$(\forall \alpha)(\alpha < \lambda \Rightarrow |\alpha|^{\aleph_0} 
< \lambda)$; see \cite[3.15]{Sh:506}. (Note if $A_n \in D$, 
$\bigcap\limits_{n < \omega} A_n = \varnothing$, then 
$T_D(f) = T_D(f)^{\aleph_0}$.)

See more in \cite[\S3]{Sh:506}.
\end{note}

\begin{note}
\label{4.3a}
On sufficient conditions for $T_J(\bar\lambda) \ge
\lambda$ and $T_J(\bar\lambda) = \lambda$, see \cite{Sh:829}.
\end{note}

\begin{note}
\label{4.4}
Assume $D$ is a filter on $\kappa,\mu = \cf(\mu) >
2^\kappa,f \in {}^\kappa\Ord$ and: $D$ is $\aleph_1$-complete or
$(\forall \sigma < \mu)(\sigma^{\aleph_0} < \mu)$.  Then $(\exists A \in D^+)
T_{D+A}(f) \ge \mu$ \Iff \, for some $A \in D^+$ and 
$\langle \lambda_i:i < \kappa\rangle = \bar \lambda \le_{D+A} f$ 
we have $\prod\limits_{i < \kappa} \lambda_i/(D+A)$ has true cofinality
$\mu$ (for approximations see \cite[\S3]{Sh:506}, proof \cite[1.1]{Sh:589},
note $\Leftarrow$ is trivial).  This is connected to the problem of the depth
of products (e.g. ultraproducts) of Boolean Algebra.
\end{note}

\begin{note}
\label{4.5}
If $2^{2^\kappa} \le \mu <T_D(\bar \lambda)$ and
$\mu^{< \theta} = \mu$, \underline{then}  for some $\theta$-complete ideal $E
\subseteq D$ we have $\mu < T_E(\bar \lambda)$, \cite[3.20]{Sh:506}.
\end{note}

\begin{note}
\label{4.6}
On $\prod\limits_{i < \kappa} \lambda_i/D$ see
\cite[3.1-3.9B]{Sh:506}, essentially this gives full pcf characterization
when it is $> 2^\kappa$.  In particular for an ultrafilter $D$ on $\kappa$
with regularity $\theta$ (i.e. not $\theta$-regular but $\sigma$-regular
for $\sigma < \lambda$) and $\lambda_i > 2^\kappa$, we have

\[
\prod\limits_{i < \kappa} \lambda_i/D = \sup \Big(\tcf 
\prod\limits_{i < \kappa} \{ \lambda'_i/D : 2^\kappa < \lambda'_i = 
\cf(\lambda'_i) \le \lambda_i\}\Big)^{< \reg(D)}
\]

\mn
(see mainly \cite[3.9]{Sh:506}).
\end{note}

\begin{note}
\label{4.956}
Assume $\LL \lambda_i:i < \kappa \RR$ tends to
$\lambda$.  A full characterization of $\prod\limits_{i < \kappa}
\lambda_i/D = \lambda$ (via weak normal ultrafilters) appears in
\cite{Sh:956}. 
\end{note}

\begin{note}
\label{4.7}
Assume $D$ is an ultrafilter on $\kappa$ and $\theta$ is
the regularity of $D$ (i.e. minimal $\theta$ such that $D$ is
not $\theta$-regular).  \underline{Then}  every $\lambda = 
\lambda^\theta > 2^\kappa$ can be represented as 
$\prod\limits_{i < \kappa} \lambda_i/D$.  (Note $\lambda = 
\lambda^{< \theta}$ is necessary) (see \cite[\S6]{Sh:589}).
\end{note}

\begin{note}
\label{4.8}
Assume $\theta < \kappa$, $J_* = [\kappa]^{< \theta}$, and
$\lambda > \kappa^\theta$. Then
\begin{equation*}
\begin{array}{clcr}
T_{J_*}(\lambda) = \sup \Big\{\tcf\big(\prod\limits_{n < n_i \atop
i < \kappa} \lambda_{i,n}/J\big)& :\ 
n_i < \omega,\ \lambda_{i,n} \text{ regular} \in [\kappa^\theta,\lambda), \\
  & J \text{ is an ideal on } 
\bigcup\limits_{i < \kappa} \{i\} \times n_i,\\
  &A \subseteq \kappa,\  |A| \ge \theta \Rightarrow 
\bigcup\limits_{i \in A} \{i\} \times n_i \in J^+ \text{ and} \\
  &\prod\limits_{n < n_i \atop i < \kappa} \lambda_{i,n}/J
\text{ has true cofinality}\Big\}.
\end{array}
\end{equation*}

\mn
This is just a case of the ``$\theta$-almost disjoint family $\subseteq
[\lambda]^\kappa$" problem as clearly $T_J(\lambda) = \sup\big\{\cA : \cA
\subseteq [\lambda]^\kappa \text{ is } \theta$-almost disjoint; i.e. $A \ne
B \in {\cA} \Rightarrow |A \cap B| < \theta \big\}$. 

See \cite[\S6]{Sh:410}.
\end{note}

\begin{note}
\label{4.9}
If $\lambda \ge \kappa > \beth_\omega(\theta)$ then in
\ref{4.7}, $T_J(\lambda) = \lambda$. 

(See \cite{Sh:460}).
\end{note}

\begin{note}
\label{4.10}
(\cite[1.2]{Sh:430}).  Assume $\lambda > \mu = \cf
(\mu) > \theta > \aleph_0$ and $\cov(\theta,\aleph_1,\aleph_1,2) < \mu$.
Then the following are equal

\begin{equation*}
\begin{array}{clcr}
\lambda(0) = \min \big\{\kappa:&\text{if } {\ga} \subseteq \Reg \cap 
\lambda \setminus \mu,|{\ga}| \le \theta \text{ then we can partition } 
{\ga} \text{ to } \langle {\ga}_n:n < \omega \rangle \\
  &\text{ such that } {\gb} \subseteq {\ga}_n \text{ and }
|{\gb}| \le \aleph_0 \Rightarrow \max \pcf ({\gb}) \le \kappa \\
  &\text{ and } [{\ga}_n]^{\le \aleph_0} \text{ is included in the ideal
generated by} \\
  &\{{\gb}_\theta[{\ga}_n]:\theta \in {\gd}_n\} \text{ for some }
{\gd}_n \subseteq \kappa^+ \cap \pcf({\ga}_n) \text{ of cardinality} 
< \mu \big\}
\end{array}
\end{equation*}

\begin{equation*}
\begin{array}{clcr}
\lambda(1) = \min\big\{|{\cP}|:&{\cP} \subseteq [\lambda]^{<\mu}
 \text{ and for every } A \in [\lambda]^{\le \theta}  \\
  &\text{ for some partition } \langle A_n:n < \omega \rangle \text{ of } A
\text{ we have:} \\
  &\langle {\cP}_n:n < \omega \rangle,\ {\cP}_n \subseteq {\cP},\ 
|{\cP}_n| < \mu,\ \mu > \underset {B \in {\cP}_n} {\to} \sup(B) \\
  &\text{ and } n < \omega \text{ and } a \in [A_n]^{\aleph_0} \Rightarrow (\exists
A \in {\cP}_n)[a \subseteq A]\big\}.
\end{array}
\end{equation*}
\end{note}
\newpage

\section {pcf theory: \cite{Sh:68}, \cite[Ch.XIII,\S5,\S6]{Sh:b},
\cite{Sh:282}, \cite{Sh:345}} \label{A5}
\bigskip

${\ga}$ denotes a set of regulars, $\min({\ga}) > |{\ga}|$ (except for a
generalization in \cite[\S3]{Sh:371}).

For a partial order $P$ let $\cf(P) = \min\{|A|:A \subseteq P,
\bigwedge\limits_{p \in P} \, \bigvee\limits_{q \in A} p \le q\}$.  We
say that $P$ has true cofinality if it has a well ordered cofinal 
subset whose cofinality is called $\tcf(P)$ (equivalently - a linearly 
ordered cofinal subset).

\begin{note}
\label{5.1}
$J_{< \lambda}[{\ga}],J_{\le \lambda}[{\ga}]$ 
\cite[1.2(2),(3)]{Sh:345a}, also \cite[3.1]{Sh:345a},
\cite[3.1]{Sh:371}.  
For example we define:

\[
J_{< \lambda}[{\ga}] =: \{{\gb} \subseteq {\ga}:
\text{ for every ultrafilter } D \text{ on } {\gb},
\cf(\textstyle\prod {\gb}/D) < \lambda\}.
\]

\[
J_\lambda[{\ga}] =: \{{\gb} \subseteq {\ga}:
\text{ for every ultrafilter } D \text{ on } {\gb},
\tcf(\textstyle\prod {\gb}/D) \ne \lambda\}.
\]
\end{note}

\begin{note}
\label{5.2}
Definition of variants of $\pcf$ 
\cite[1.2(1),(2)]{Sh:345a}, \cite[3.1]{Sh:371} for example

\[
\pcf({\ga}) = \{\cf(\textstyle\prod {\ga}/D):D \text{ an ultrafilter on } {\ga}\}.
\]

\mn
For given cardinals $\theta > \sigma$ let

\begin{equation*}
\begin{array}{clcr}
\pcf_{\Gamma(\theta,\sigma)}({\ga}) = \{\tcf(\textstyle\prod  {\gb}/J):\,&{\gb} 
\subseteq {\ga},|{\gb}| < \theta,J
\text{ is a } \sigma \text{-complete ideal} \\
  &\,\text{ on } {\gb} \text{ and } \textstyle\prod  {\gb}/J
\text{ has true cofinality}\},
\end{array}
\end{equation*}

\mn
$\Gamma(\theta)$ means $\Gamma(\theta^+,\theta)$.
\end{note}

\begin{note}
\label{5.3}
Trivial properties \cite[1.3,1.4]{Sh:345a}.
\end{note}

\begin{note}
\label{5.4}
Basic properties \cite[1.5,1.8,2.6,2.8,2.10,2.12]{Sh:345a}.
\end{note}

\begin{note}
\label{5.5}
$|\pcf(\ga)| 
\le 2^{|{\ga}|}$ \cite[1.8(5)]{Sh:345a}, 
$\pcf({\ga})$ has a last member. \cite[1.9]{Sh:345a} Also, if 
$|{\ga} \cup {\gb}| < \min({\ga} \cup {\gb})$ then $(\pcf({\ga})) \cap
(\pcf({\gb}))$ has a last member; actually, $|{\ga}| < \min ({\ga})$, 
$|{\gb}| < \min({\gb})$ suffices (by \cite[6.4A]{Sh:430}, we
can take intersections of many ${\ga}_i$).
\end{note}

\begin{note}
\label{5.6}
If $D,D_i \, (i < \kappa < \min({\ga}))$ are filters 
on ${\ga}$, $E$ a filter on $\kappa$,

\[
D = \{{\gb} \subseteq {\ga}:\{i:{\gb} \in D_i\} \in E\}
\text{ and } \lambda_i = \tcf(\textstyle\prod {\ga}/D_i)
\]

\mn
(well defined) \underline{then}  $\tcf(\prod {\ga}/D_i)$, and 
$\tcf(\prod\limits_{i < \kappa} \lambda_i/E)$
are equal \cite[1.10]{Sh:345a}.  Moreover, $\bigwedge\limits_{i} \kappa <
\lambda_i$ is enough \cite[1.11]{Sh:345a}. 

(And see more in \cite[3.3,3.6]{Sh:410}, generalization
\cite[1.10]{Sh:506}).
\end{note}

\begin{note}
\label{5.7}
(Repeating \ref{2.3})  What is $\Ch^{\ga}_N$ (where
$N \prec (\clH(\chi),\in)$), $\Ch^{\ga}_N(\theta) = \sup(N \cap \theta)$
for $\theta \in {\ga}$) \cite[3.4]{Sh:345a}, \cite[3.4]{Sh:355},
``stationary $F \subseteq \prod{\ga}$" \cite[1.2,1.3,1.4]{Sh:371}, more
in \cite[\S6]{Sh:430}. 
\end{note}

\begin{note}
\label{5.8}
$\cf(\prod {\ga}) = \max \pcf({\ga})$
\cite[3.1,more 3.2]{Sh:355}, \cite[3.4]{Sh:345a}, other representation
\cite{Sh:506}.
\end{note}

\begin{note}
\label{5.9}
There is a generating sequence $\langle {\gb}_\theta
[{\ga}]:\theta \in \pcf({\ga}) \rangle$; i.e.
$J_{\le \theta}[{\ga}] = J_{< \theta}[{\ga}] + {\gb}_\theta
[{\ga}]$, so $J_{< \lambda}[{\ga}]$ is the ideal on ${\ga}$
generated by $\{{\gb}_\theta[{\ga}]:\theta < \lambda\}$ and $\prod
{\gb}_\theta[{\ga}]/J_{< \theta}[{\ga}]$ has true cofinality
$\theta$ and $J_\lambda[{\ga}]$ is the ideal on ${\ga}$ generated by
$\{{\gb}_\theta[{\ga}]:\theta < \lambda\} \cup \{{\ga} \setminus
{\gb}_\lambda[{\ga}]\}$; \cite[2.6]{Sh:371} also \cite[3.1]{Sh:345a} 
+ \cite[\S1]{Sh:420}, more in \cite[4.1A]{Sh:400}; nice good 
cofinal $\bar f$: \cite[\S3]{Sh:345a}, \cite[3.4A]{Sh:355}, 
\cite[1.2,1.3,1.4]{Sh:371}, \cite[4.1A(2)]{Sh:400}.
Another representation is included in \cite{Sh:506} (see \ref{2.7} on the
framework and \ref{5.23} below); it uses \ref{0.4} from \cite{Sh:420}.
\end{note}

\begin{note}
\label{5.1a}
If ${\gb} \subseteq {\ga}$ and ${\gc} =
\pcf({\gb})$, \underline{then}  for some finite ${\gd} \subseteq
{\gc},{\gb} \subseteq \bigcup\limits_{\theta \in {\gd}} {\gb}_\theta
[{\ga}]$ see \cite[3.2(5)]{Sh:345a}.
\end{note}

\begin{note}
\label{5.11}
Cofinality sequence \cite[3.3]{Sh:345a}, 
\cite[2.1]{Sh:371}, more in the proof of \cite[4.1]{Sh:400}.
\end{note}

\begin{note}
\label{5.12}
$\bar f$ is $x$-continuous (nice) \cite[3.3,3.5,3.8(1),(2)]{Sh:345a}.
\end{note}

\begin{note}
\label{5.13}
For a discussion of when ${\ga}$ has a generating
sequence which is smooth and/or closed \cite[3.6,3.8(3)]{Sh:345a},
\cite[4.1A(4)]{Sh:400}; {\it smooth} means $\mu \in {\gb}_\lambda
[{\ga}] \Rightarrow {\gb}_\mu[{\ga}] \subseteq {\gb}_\lambda
[{\ga}]$, {\it closed} means $\pcf({\gb}_\lambda[{\ga}]) =
{\gb}_\lambda[{\ga}]$.  If for example $|\pcf({\ga})| <
\min({\ga})$ we can have both \cite[3.8]{Sh:345a} and more, then
we can use the ``pcf calculus" style of proof.  Proofs in this style can 
generally be carried further but become a little complicated, 
as done in \cite[6.7-6.7E]{Sh:430} (particularly \cite[6.7C(3)]{Sh:430}.)
On a generalization, see \cite{Sh:506}.
\end{note}

\begin{note}
\label{5.14}
If $\lambda = \max \pcf({\ga})$, and $\mu =: \sup
(\lambda \cap \pcf({\ga}))$ is singular, \underline{then}  for ${\gc}
\subseteq \pcf({\ga})$ unbounded in $\mu$, $\tcf(\prod{\gc}/J^{\bd}_\mu) = \lambda$ 
\cite[3.7]{Sh:345a}, \cite[2.10(2)]{Sh:371}, where for $A$ a set of ordinals, 
$J^{\bd}_A = \{B \subseteq A : \sup(B) < \sup (A)\}$.
\end{note}

\begin{note}
\label{5.15}
If $\lambda \in \pcf({\ga})$, then for some
${\gb} \subseteq {\ga}$ we have: $\lambda = \max \pcf({\gb})$ and
$\lambda \cap \pcf({\gb})$ has no last element and $\lambda \notin
\pcf({\ga} \setminus {\gb})$; see \cite[,2.10(1)]{Sh:371}.
\end{note}

\begin{note}
\label{5.16}
If $(\forall \mu)[\mu < \lambda \Rightarrow \mu^{< \kappa}
< \lambda]$, then $J_{< \lambda}[{\ga}]$ is $\kappa$-complete
\cite[1.6(1)]{Sh:371}.
\end{note}

\begin{note}
\label{5.17}
Localization: if $\lambda \in \pcf({\gb})$,
${\gb} \subseteq \pcf({\ga})$ (and we assume just
$|{\gb}| < \min({\gb})$), \then , for some $\gc \subseteq \gb$ 
we have $|{\gc}| \le |{\ga}|$ and $\lambda \in \pcf({\gc})$, \cite[3.4]{Sh:371}.

Also if $\lambda \in \pcf_{\sigma\text{-complete}}({\gb})$,
${\gb} \subseteq \pcf({\ga})$ then for some $\gc \subseteq \gb$, we have 
$|{\gc}| \le |{\ga}|$ and $\lambda \in \pcf_{\sigma\text{-complete}}({\gc})$; see \cite[6.7F(4),(5)]{Sh:430}.
\end{note}

\begin{note}
\label{5.18}
\mn
\begin{enumerate}
\item[$(a)$]   $\pcf({\ga})$ cannot contain an interval of
Reg ($=$ the class of regulars) of cardinality $|{\ga}|^{+4}$.
\end{enumerate}
\mn
In fact:
\mn
\begin{enumerate}
\item[$(b)$]   for no ${\ga}$ and $\chi$ is $\{i < |{\ga}|^{+4}:
\chi^{+i+1} \in \pcf({\ga})\}$ unbounded in $|{\ga}|^{+4}$.
\end{enumerate}
\mn
[Why?  If so, there is $\lambda \in \pcf((\chi,\chi^{+|{\ga}|
^{+4}}) \cap \pcf({\ga}))$ such that $\lambda > \chi^{|{\ga}|^{+4}}$,
hence by localization for some ${\gc} \subseteq
(\chi,\chi^{+|{\ga}|^{+4}}) \cap \pcf({\ga})$ of cardinality
$\le |{\ga}|$ we have $\lambda \in \pcf({\gc})$, hence for
some limit ordinal $\delta < |{\ga}|^{+4}$, 
$\pp_{|{\ga}|}(\chi^{+ \delta}) \ge \lambda \ge \chi^{+|{\ga}|^{+4}}$ and we get 
a contradiction by \cite[\S4]{Sh:400}.]
\end{note}

\begin{note}
\label{5.19}
Defining $(\mu,\theta,\sigma)$-inaccessibility 
\cite[3.1,3.2]{Sh:410}.
\end{note}

\begin{note}
\label{5.20}
On $\pcf({\gb})$ for ${\gb} \subseteq \pcf({\ga}),|{\gb}| < \min({\gb})$ or
even with no inaccessible accumulation points, see \cite[1.12]{Sh:345a},
\cite[\S3]{Sh:371}, mainly: having $\mathfrak{b} 
^*_\lambda[{\ga}] \subseteq
\pcf[{\ga}]$.
\end{note}

\begin{note}
\label{5.21}
Uniqueness of $\bar f$ ($<_J$-increasing cofinal) 
\cite[2.7,2.10]{Sh:345a}.
\end{note}

\begin{note}
\label{5.22}
If $J = J_{< \lambda}[{\ga}]$, $\lambda = \tcf(\prod {\ga}/J)$, 
${\ga} = \bigcup\limits_{i < \alpha}{\ga}_i$, \underline{then}  
for some finite ${\gb}_i \subseteq \pcf({\ga}_i)$ (with $i < \alpha$) 
we have $\lambda = \max \pcf(\bigcup\limits_{i < \alpha} {\gb}_i)$, 
and for $w \subseteq \alpha$ we have
$$\max\pcf\Big(\textstyle\bigcup\limits_{i \in w} {\gb}_i\Big) < \lambda \Leftrightarrow
\Big(\bigcup\limits_{i \in w} {\ga}_i \Big) \in J$$ 
\cite[\S1]{Sh:371}, more in \cite[\S6]{Sh:430}.
\end{note}

\begin{note}
\label{5.23}
If $\bar \lambda = \langle \lambda_i:
i < \kappa \rangle,I^*$ a weakly $\theta$-saturated ideal on $\kappa$ (see
below) $\theta = \cf(\theta) < \lambda_i$, \underline{then}  the pcf
analysis, e.g. from \ref{5.9} holds for $\bar \lambda$ when we restrict
ourselves to ideals on $\kappa$ extending $I^*$ (see 
\cite[\S1,\S2]{Sh:506}). 

E.g. $\theta$ can play the role of $\kappa =
\Dom(I)$ if $I$ is weakly $\theta$-saturated, i.e.
\mn
\begin{enumerate}
\item[$(*)_{I,\theta}$]  there is no division of $\kappa$ to $\theta$ sets
none of which is in $I$.
\end{enumerate}
\end{note}

\begin{note}
\label{5.24}
If $|{\ga}| < \min({\ga}),
\aleph_0 \le \sigma = \cf(\sigma)$, then for some $\alpha < \sigma$
and $\lambda_\beta,{\ga}_\beta \, (\beta < \alpha)$ we have
\mn
\begin{enumerate}
\item[(i)]   ${\ga} = \bigcup\limits_{\beta < \alpha} {\ga}_\beta$,
\sn
\item[(ii)]  $\lambda_\beta = \max \pcf({\ga}_\beta)$
\sn
\item[(iii)]   $\lambda_\beta \notin \pcf({\ga} \setminus {\ga}_\beta)$ and
\sn
\item[(iv)]   $\lambda_\beta \in \pcf_{\sigma-\text{\rm com}}({\ga}_\beta)$.
\end{enumerate}
\mn
[Why?  We prove this by induction on $\max \pcf({\ga})$, hence by the
induction hypothesis we can ignore (iii) as we can regain it. Now let

\begin{equation*}
\begin{array}{clcr}
J = \big\{\gb \subseteq \ga:\,&\text{we can find } \alpha <
\sigma,\ \langle {\ga}_\beta : \beta < \alpha \rangle \\
  &\,\text{such that } {\gb} = \bigcup\limits_{\beta < \alpha} 
{\ga}_\beta \text{ and (ii), (iv) above}\big\}.
\end{array}
\end{equation*}

\mn
Clearly $J$ is a family of subsets of ${\ga}$, includes the
singletons, and is closed under subsets and under unions of 
$< \sigma$ members.  If ${\ga} \in J$ we are done.  
If not, choose ${\gc} \subseteq {\ga}$
such that ${\gc} \notin J$ and (under these restrictions)
$\lambda_{\gc} =: \max \pcf({\gc})$ is minimal.  Now by the
minimality of $\lambda_{\gc}$, $J_{< \lambda_{\gc}}[{\ga}] \subseteq J$, 
so ${\gb}_{\lambda_{\gc}}[{\ga}]$ satisfies the requirement
for ${\gb} \in J$ (with $\alpha = 1$). Contradiction].
\end{note}

\begin{note}
\label{5.25}
See more in \ref{7.18} and \cite{Sh:497} and particularly \cite{Sh:513}.
\end{note}

\begin{note}
\label{5.26}
If $\lambda = \max \pcf({\gb})$ and $\lambda \cap
\pcf({\gb})$ has no last element (see \ref{5.15}) and $\mu <
\sup(\lambda \cap \pcf({\gb}))$, \underline{then}  for some ${\gc}
\subseteq \text{ pcf}({\mathfrak b}) \setminus \mu$ of cardinality $\le
|{\gb}|$ we have $\lambda = \max \pcf({\gc})$ and $\theta \in
{\gc} \Rightarrow \max \pcf(\theta \cap {\gc}) < \theta$ (see
\cite[2.4A,2.4(2)]{Sh:413}, an ex.).
\end{note}

\begin{note}
\label{5.29}
Above, the demand $|\ga| < \min(\ga)$ was essential, but:
\begin{itemize}
    \item If $\ga \subseteq \pcf(\gb)$, $|\gb| < \min(\gb)$ \underline{then},
    almost always, the pcf theorem holds for $\ga$, fully by \cite{Sh:E11}. 
    
    \item In the other direction, see \ref{5.23}.
\end{itemize}
\end{note}

\newpage

\section {Representation and pp} \label{A6}
\bigskip

\begin{note}
\label{6.1}
Definition of pp and variants \cite[1.1]{Sh:355}.  For
$\lambda$ singular 

\begin{equation*}
\begin{array}{clcr}
\pp_\theta(\lambda) = \sup\{\tcf(\textstyle\prod{\ga}/J):\,&{\ga}
\text{ is a set of } \le \theta \text{ regular cardinals}, \\
  &\,\text{unbounded in } \lambda,\ J \text{ an ideal on } {\ga}
\text{ including } J^{\bd}_{\ga} \\
  &\,\text{and } \textstyle\prod {\ga}/J \text{ has true cofinality}\},
\end{array}
\end{equation*}

\[
\pp(\lambda) = \pp_{\cf(\lambda)}(\lambda),
\]

\[
\pp^+(\lambda) \text{ is the first regular without such a representation}
\]

\[
\pp_\Gamma(\lambda) \text{ means that we restrict ourselves to } J
\text{ satisfying } \Gamma,
\]

\[
\pp^*_I(\lambda) = \pp_{\{I\}}(\lambda)
\]

\mn
and

\[
\pp_I(\lambda) = \sup\{\pp^*_J(\lambda):J \text{ an ideal extending } I\},
\]

\begin{equation*}
\begin{array}{clcr}
\lambda =^+ &\pp(\lambda) \text{ means more than equality; the 
supremum in the right hand side} \\
  &\text{is obtained if it is regular}.
\end{array}
\end{equation*}
\end{note}

\begin{note}
\label{6.2}
Downward closure: 

If $\lambda = \tcf(\prod\limits_{i < \kappa} \lambda_i/I)$, 
$\lambda_i = \cf(\lambda_i) > \kappa$, and 
$\kappa < \lambda' = \cf(\lambda') < \lambda$, \underline{then} 
for some $\lambda'_i$ we have
$\kappa \le \lambda'_i = \cf(\lambda'_i) < \lambda_i$ and 
$\lambda' =\tcf(\prod\limits_{i < \kappa} \lambda'_i/I)$. Moreover, 
$$\tlim_I \lambda_i = \mu < \lambda' < \lambda \Rightarrow \tlim_I \lambda'_i = \mu$$ 
$\lambda' = \tcf(\prod \lambda'_i,\le_I)$ is exemplified by
$\mu^+$-free $\bar f$, which means: 
if $w \subseteq \lambda' \text{ and } |w| \le \mu$,
then for some $\langle s_\alpha:\alpha \in w \rangle$, $s_\alpha \in I$ and
for each $i < \kappa$, $\langle f_\alpha(i):\alpha \in w,i \notin s_\alpha
\rangle$ is without repetition; in fact, we get ``strictly increasing."
\cite[1.3,1.4,2.3]{Sh:355} more \cite[4.1]{Sh:400} a generalization 
\cite[3.12]{Sh:506}.
\end{note}

\begin{note}
\label{6.3}
If $\lambda > \kappa \ge \cf(\lambda)$, $I$ an ideal on
$\kappa$, $\kappa$ is an increasing union of $\cf(\lambda)$ members of $I$,
\underline{then}  $\big\{\tcf \big(\prod\limits_{i < \kappa} \lambda_i/I\big) : \tlim_I \lambda_i = \lambda$ and $\lambda_i = \cf(\lambda_i)\big\}$ is an initial
segment of $\Reg \setminus \lambda$, so the first member is $\lambda^+$,
\cite[1.5,2.3]{Sh:355}.
\end{note}

\begin{note}
\label{6.4}
If $\lambda > \cf(\lambda) > \aleph_0$, \underline{then} 
for some increasing continuous $\langle \lambda_i:i < \cf(\lambda)
\rangle$ with limit $\lambda,\prod\limits_{i < \cf(\lambda)}
\lambda^+_i/J^{\bd}_{\cf(\lambda)}$ has 
true cofinality $\lambda^+$, \cite[2.1]{Sh:355}.
\end{note}

\begin{note}
\label{6.5}
{\rm pp}$(\lambda) > \lambda^+$ 
contradicts ``large cardinal"
type assumptions, for example ``every $\mu$-free abelian group is free"
\cite[2.2,2.2B]{Sh:355}, for the parallel fact on cov see 
\cite[6.6]{Sh:355}.
\end{note}

\begin{note}
\label{6.6}
\mn
\begin{enumerate}
\item[$(a)$]  (inverse monotonicity)  If $\mu > \lambda > \kappa \ge
\cf(\lambda) + \cf(\mu)$ and $\pp^+_\kappa(\lambda) > \mu$,
\underline{then}  $\pp^+_\kappa(\lambda) \ge \pp^+_\kappa(\mu)$
\sn
\item[$(b)$]  so given $\kappa_0 \le \kappa_1 < \mu$ if $\lambda$ is minimal
such that $\lambda > \kappa_1 \ge \kappa_0 \ge \cf(\lambda)$, 
$\pp(\lambda) \ge \mu$, \then: 
${\ga} \subseteq \Reg \cap [\kappa_1,\lambda)$, $|{\ga}| \le \kappa_0$, 
$\sup({\ga}) < \lambda$ implies $\max \pcf({\ga}) < \lambda$;
equivalently, $\lambda' \in (\kappa_1,\lambda)$ and  
$\cf(\lambda') \le \kappa_0 \Rightarrow \pp_{\kappa_1}(\lambda') < \lambda$ 
\cite[2.3]{Sh:355} (with more)
\sn
\item[$(c)$]   assume $\kappa \le \chi < \mu$, and
\[
(\forall \lambda)[\lambda \in (\chi,\mu) \text{ and } \cf(\lambda) \le \kappa
\Rightarrow \pp(\lambda) < \mu]
\]
\underline{then}  for every ${\ga} \subseteq (\chi,\mu)$ of cardinality
$\le \kappa$, $\sup({\ga}) < \mu$ we have $\max \pcf({\ga}) < \mu$ 
[by $(d)$ below and \ref{6.9} below]
\sn
\item[$(d)$]   $\max \pcf({\ga}) \le \sup\{\pp_{|{\ga} 
\cap \mu|}(\mu):\mu \notin {\ga},\ \mu = \sup({\ga} \cap \mu)\}$ 
[by the definition].
\end{enumerate}
\mn
Similar assertion holds for pp$_\Gamma,\Gamma$ is ``nice" enough.
\end{note}

\begin{note}
\label{6.7}
\mn
\begin{enumerate}
\item[$(A)$]   If $\lambda$ is singular, $\mu < \lambda$, \underline{then}  for
some $\delta \le \cf(\lambda)$ and increasing sequence $\langle
\lambda_i:i < \delta \rangle$ of regular cardinals in $(\mu,\lambda)$ 
and $\delta = \cf(\delta) \vee \delta < \omega_1$ we have: 
$\lambda_i > \max \pcf\{\lambda_j:j < i\}$ and $\lambda^+ = 
\tcf(\prod \lambda_i/J^{\bd}_\delta$), \cite[3.3]{Sh:355}
\sn
\item[$(B)$]  If $\lambda$ is singular and $\aleph_0 < \cf(\lambda)
= \kappa$ and $\bigwedge\limits_{\mu < \lambda} \mu^\kappa < \lambda$ and
$\lambda < \theta = \cf(\theta) \le \lambda^\kappa$, \underline{then}  for
some increasing sequence $\langle \lambda_i:i < \kappa \rangle$ of regulars
$< \lambda,\lambda = \sum\limits_{i < \kappa} \lambda_i$ and 
$\prod\limits_{i < \kappa} \lambda_i/J^{\bd}_\kappa$ has 
true cofinality $\theta$ (see \cite[1.6(2)]{Sh:371}).  
Moreover, we can demand $i < \kappa \Rightarrow \max \pcf\{\lambda_j:
j < i\} < \lambda_i$.  We can weaken the hypothesis
to $\aleph_0 < \kappa = \cf(\lambda) < \lambda_0 < \lambda$ and
$(\forall \mu)[\lambda_0 < \mu < \lambda $ 
and 
$ \cf(\mu) \le \kappa
\Rightarrow \pp(\mu) < \lambda]$ (see \cite[1.6(2)]{Sh:371}.  
If we allow $\cf(\lambda) = \kappa = \aleph_0$ we still get this, but 
for possibly larger $J$, see \cite[6.5]{Sh:430}.
\end{enumerate}
\end{note}

\begin{note}
\label{6.8}
$\pp_{\Gamma(\theta,\sigma)}$ can be reduced to finitely many
$\pp_{\Gamma(\theta)}$, see \cite[5.8]{Sh:355}.
\end{note}

\begin{note}
\label{6.9}
If $\mu > \theta \ge \cf(\mu)$ and for every large
enough $\mu' < \mu$:

\[
[\cf(\mu') \le \theta \Rightarrow \pp_\theta(\mu') < \mu]
\]

\mn
\then

\[
\pp(\mu) =^+ \pp_\theta(\mu) =^+ \pp_{\Gamma(\cf(\mu))}(\mu)
\]

\mn
\cite[1.6(3)(5),1.6(2)(4)(6),1.6A]{Sh:371}.
\end{note}

\begin{note}
\label{6.10}
If $\langle {\gb}_\zeta:\zeta < \kappa \rangle$ is
increasing, $\lambda \in \pcf\big(\bigcup\limits_{\zeta < \kappa} {\gb}_\zeta\big)
\setminus \bigcup\limits_{\zeta < \kappa} \pcf({\gb}_\zeta)$, \underline{then} :
\mn
\begin{enumerate}
\item  for some ${\gc} \subseteq \bigcup\limits_\zeta \pcf({\gb}_\zeta),
|{\gc}| \le \kappa$, we have $\lambda \in \pcf({\gc})$
\sn
\item  if $\kappa = \cf(\kappa) > \aleph_0$, \underline{then} for some club
$C \subseteq \kappa$ and $\lambda_\zeta \in \pcf\big(\bigcup\limits_{\xi < \zeta}
{\gb}_\xi\big)$ for $\zeta \in C$, we have 
$\lambda = \prod\limits_{\zeta \in C} \lambda_\zeta/J^{\bd}_\kappa$, 
$\lambda_\zeta\ (\zeta \in C)$ is increasing and 
$$\zeta \in C \Rightarrow \lambda_{\zeta +1} > 
\max \pcf\{\lambda_\xi : \xi \in C \text{ and } 
\xi \le \zeta\}$$ \cite[1.5]{Sh:371}.
\end{enumerate}
\end{note}

\begin{note}\label{6.11}
If $\mu > \theta \ge \cf(\mu) \ge \sigma =
\cf(\sigma)$, and for every large enough $\mu' < \mu$:

\[
[\sigma \le \cf(\mu') \le \theta \Rightarrow \pp_{\Gamma
(\theta^+,\sigma)}(\mu') < \mu]
\]

\mn
\then

\[
\pp_{\Gamma(\cf(\mu)^+,\sigma)}(\mu) = \pp_{\Gamma(\cf(\mu))}(\mu)
\]

\mn
\cite[1.2]{Sh:410} and more there.
\end{note}

\begin{note}\label{6.12}:  If $\mu > \kappa = \cf(\mu) > \aleph_0$ and for
every large enough $\mu' < \mu$

\[
(\mu')^\kappa < \mu \text{ or just } [\cf(\mu') \le 
\cf(\mu) \Rightarrow \pp_\kappa(\mu') < \mu],
\]

\mn
\underline{then}  $\pp^+(\mu) = \pp^+_{J^{\bd}_\kappa}(\mu)$ and we can 
get the conclusion in \ref{6.7}(B) above 
\cite[1.8]{Sh:371}.  
Generalization for $\Gamma(\theta,\sigma)$ in \cite[1.2]{Sh:410}.
\end{note}

\begin{note}\label{6.13}:  If $\lambda > \kappa = \cf(\lambda) > \aleph_0,
\lambda > \theta$ \underline{then}  for some increasing continuous sequence $\langle
\lambda_i:i < \kappa \rangle$ with limit $\lambda$:
\mn
\begin{enumerate}
\item[$(a)$]   for every $i < \kappa,\lambda_i < \mu < \lambda_{i+1} \text{ and }
\cf(\mu) \le \theta \Rightarrow \pp_\theta(\mu) < \lambda_{i+1}$
\end{enumerate}
\mn
\underline{or}
\mn
\begin{enumerate}
\item[$(b)$]   for every $i < \kappa,\pp_{\theta + \cf(i)}
(\lambda_i)\ge \lambda$ \cite[1.9; more 1.9A]{Sh:371}.
\end{enumerate}
\end{note}

\begin{note}
\label{6.14}
If $\sigma \le \cf(\mu) \le \theta < \kappa < \mu$ then:

\[
\pp_\theta(\mu) < \mu^{+ \theta^+} \Rightarrow \pp_\kappa(\mu) = 
\pp_\theta(\mu);
\]

\mn
and

\[
\pp_{\Gamma(\theta^+,\sigma)}(\mu) < \mu^{+ \theta^+} \Rightarrow 
\pp_{\Gamma(\kappa^+,\sigma)}(\mu) = \pp_\theta(\mu)
\]

\mn
\cite[3.6; more 3.7,3.8]{Sh:371}.
\end{note}

\begin{note}
\label{6.15}
If $\langle \mu_i:i \le \kappa \rangle$ is increasing
continuous, $\mu_0 > \kappa^{\aleph_0} > \kappa = \cf(\kappa) >
\aleph_0$ and $\cov(\mu_i,\mu_i,\kappa^+,2) < \mu_{i+1}$, then for
some club $E$ of $\kappa$ we have 
$$\delta \in E \cup \{\kappa\} \Rightarrow
\pp_{J^{\bd}_{\cf(\delta)}}(\mu_\delta) = \cov(\mu_\delta,
\mu_\delta, \kappa^+,2)$$ 
so e.g. for most limit $\delta < \omega_1$,
pp$_{J^{\bd}_\omega}(\beth_\delta) =^+ \beth_{\delta +1}$ (see
\cite[part C, remark to X, \S5, p.412]{Sh:E12}.
\end{note}

\begin{note}
\label{6.16}
If $\pp^+_\sigma(\mu) > \lambda = \cf(\lambda)$,
(so $\cf(\mu) \le \sigma$) then
\mn
\begin{enumerate}
\item[$(a)$]   for some ${\ga}$ an unbounded subset of $\mu$, 
$|{\ga}| \le \sigma$, $\lambda = \tcf(\prod {\ga}/J^{\bd}_\ga) =
\max \pcf({\ga})$ \underline{or}
\sn
\item[$(b)$]    for some ${\ga} \subseteq (\mu,\lambda)$ of cardinality
$\le \sigma$, $\lambda = \max \pcf({\ga})$ and $\theta \in {\ga}
\Rightarrow \max \pcf(\theta \cap {\ga}) < \theta$ (see
\cite[2.4A,2.4(2)]{Sh:413}).
\end{enumerate}
\end{note}

\begin{note}\label{6.50}
For singular $\mu$ we like to have, for $\lambda \in (\mu, \pp^+(\mu)) \cap \Reg$, 
an increasing sequence $\bar \lambda$ of regulars of length $\cf(\mu)$ converging 
to $\mu$ such that $\lambda = \tcf (\prod \bar \lambda, <_{J^{\bd}_{\cf\mu}})$. 
We had such results for $\lambda = \mu^+$ and when 
$\alpha < \mu \Rightarrow |\alpha|^{\cf(\mu)} < \mu$; see \cite[\S1]{Sh:371} 
(and more in \cite{Sh:371}). 
\begin{enumerate}
    \item[$\otimes_1$] If $\lambda > \kappa \geq \cf(\lambda) > \aleph_0$, 
    and $\mu < \lambda$ of cofinality $\leq \kappa$ (satisfying 
    $\pp_\kappa(\mu) < \lambda$) is large enough, then
    \begin{enumerate}
        \item for every regular $\theta \in (\lambda, \pp_\kappa^+(\mu))$ there is
        an increasing sequence\\ $\LL \lambda_i : i < \cf(\lambda) \RR$ of 
        regular cardinals $< \lambda$, with limit $\lambda$, such that $\theta$ is the
        true cofinality of $\prod \lambda_i / J_{\cf\lambda}^{\bd}$. This means that
        
        \item $\pp_\kappa(\lambda) = \pp(\lambda) = \pp_{ J_{\cf\lambda}^{\bd} }^*(\lambda)$
    \end{enumerate}
\end{enumerate}
\end{note}

\begin{note}\label{6.53}
Now from a paper by Gitik and the author \cite{Sh:1013}, we quote directly:
\begin{enumerate}
    \item Assume that $\kappa > \aleph_0$ is a weakly compact cardinal. Let $\mu > 2^\kappa$ be a singular cardinal of cofinality $\kappa$. Then for every regular 
    $\lambda < \pp^+_{\Gamma(\kappa)}(\mu)$ there is an increasing sequence 
    $\LL \lambda_i : i < \kappa \RR$ of regular cardinals converging to $\mu$ such that $\lambda = \tcf(\prod\limits_{i < \kappa} \lambda_i, <_{J^{\bd}_\kappa})$.
    
    \item Let $\mu$ be a strong limit cardinal and $\theta$ a cardinal above $\mu$. Suppose that at least one of them has an uncountable cofinality. Then there is $\sigma_* < \mu$ such that for every $\chi < \theta$ the following holds: 
    $$\theta > \sup \big\{ \sup \pcf_{\sigma_*\text{-complete}}(\ga) : \ga \subseteq \Reg \cap (\mu^+, \chi),\ |\ga| < \mu \big\}.$$
    
    \medskip As an application we show that:\medskip
    
    \item if $\kappa$ is a measurable cardinal and $\bfj : \bfV \to \bfM$ is the elementary embedding by a $\kappa$-complete ultrafilter over $\kappa$, then for every $\tau$ the following holds:
    \begin{enumerate}
        \item if $\bfj(\tau)$ is a cardinal then $\bfj(\tau) = \tau$;
    
        \item $|\bfj(\tau)| = |\bfj(\bfj(\tau))|$;
    
        \item for any $\kappa$-complete ultrafilter $W$ on $\kappa$, $|\bfj(\tau)| = |\bfj_W(\tau)|$.
    \end{enumerate}
\end{enumerate}
\end{note}

\newpage

\section {Covering number} \label{A7}
\bigskip

\begin{definition}
\label{7.1}
\cite[5.1]{Sh:355}

\begin{equation*}
\begin{array}{clcr}
\cov(\lambda,\mu,\theta,\sigma) = \min \big\{|{\cP}|\ : 
& {\cP} \subseteq [\lambda]^{<\mu} \text{ such that } \forall a \subseteq \lambda \text{ with } |a| < \theta, \\
  &\exists \alpha < \sigma \text{ and } \exists A_i \in {\cP}\ 
(i < \alpha) \text{ such that } a \subseteq \bigcup\limits_{i < \alpha} A_i \big\}.
\end{array} 
\end{equation*}

\mn
So $\cov(\lambda,\kappa^+,\kappa^+,2) = \cf([\lambda]^{\le \kappa},
\subseteq)$.
\end{definition}

\begin{note}
\label{7.2}
Basic properties \cite[5.2,5.3]{Sh:355} see also
\cite[3.6]{Sh:355}; for example if $\lambda > \theta > \cf(\lambda) \ge
\sigma$, \underline{then}  for some $\mu < \lambda$ we have

\[
\cov(\lambda,\lambda,\theta,\sigma) = \cov(\lambda,\mu,\theta,\sigma).
\]
\end{note}

\begin{note}
\label{7.3}
$\cov$ and cardinal arithmetic and $T_\Gamma(\lambda)$ see e.g.
\cite[5.10]{Sh:355}, \cite[5.6,5.7,5.9,5.10, Definition of 
$T_\Gamma$]{Sh:355}. For example,

\[
\lambda^\kappa = \cov(\lambda,\kappa^+,\kappa^+,2) + 2^\kappa.
\]

\mn
By this and \ref{7.4}, \ref{7.5} below we shall use assumptions on
cases of $\pp$ rather than conventional cardinal arithmetic.
\end{note}

\begin{note}
\label{7.4}
On $\cov = \pp$: if $\lambda \ge \mu \ge \theta > \sigma = \cf(\sigma) > \aleph_0$, 
$\lambda > \mu \vee \cf(\mu) \in [\sigma,\theta)$, \underline{then} $\cov(\lambda,\mu,\theta,\sigma) = \sup
\{\pp_{\Gamma(\theta,\sigma)}(\chi): \chi \in [\mu,\lambda],\cf(\chi) \in
[\sigma,\theta)\}$, we have $=^+$ if $\mu = \theta$; \cite[5.4]{Sh:355}.

Assuming for simplicity $\lambda = \mu$, if $=^+$ fails, then for some
${\ga} \subseteq \Reg \cap \mu$ we have $|{\ga}| < \mu,\sup ({\ga}) 
= \mu$ and

\begin{equation*}
\begin{array}{clcr}
\cov(\lambda,\mu,\theta,\sigma) = \sup\{\tcf(\textstyle\prod{\gb}/J):\,
&{\gb} \subseteq {\ga},|{\gb}| < \theta,\mu = \sup({\gb}), \\
  &J \text{ is an ideal on } {\gb} \text{ extending } 
J^{\bd}_{\gb}\};
\end{array}
\end{equation*}

\mn
see \cite[6.12]{Sh:513}.
\end{note}

\begin{note}
\label{7.5}
The parallel of \ref{7.4} for $\sigma = \aleph_0$ ``usually holds", i.e.:
\mn
\begin{enumerate}
\item[$(a)$]   for $\lambda$ singular, $\cov(\lambda,\lambda,
\cf(\lambda)^+,2) = \pp(\lambda)$ if for every singular $\chi <
\lambda$, $\pp(\chi) = \chi^+$; \cite[\S1]{Sh:400} (and weaker assumptions and 
intermediate stages there)
\sn
\item[$(b)$]   if $\cf(\lambda) = \aleph_0$, 
$\bigwedge\limits_{\mu < \lambda}\mu^{\aleph_0} < \lambda$ and 
$\pp(\lambda) < \cov(\lambda, \lambda, \aleph_1, 2)$, \underline{then}  
$$\{\mu:\lambda < \mu = \aleph_\mu < \pp(\lambda)\}$$
is uncountable \cite[5.9]{Sh:400}, more in 
\cite[6.4]{Sh:430}. If $\lambda$ is a strong limit, then the set has 
cardinality $> \lambda$;
\sn
\item[$(c)$]   few exceptions: if $\langle \lambda_i:i \le \kappa \rangle$
is increasing continuous and $\kappa = \cf(\kappa) > \aleph_0$,
$\bigwedge\limits_{i < \kappa} \cov(\lambda_i,\lambda_i,\kappa^+,2) <
\lambda_\kappa$, \underline{then}  for some club $C$ of $\kappa$, 
$\delta \in C \cup \{\kappa\}$ implies equality, i.e.
$\cov(\lambda_\delta,\lambda_\delta,\aleph_1,2) = \pp(\lambda_\delta)$.
\cite[5.10]{Sh:400}
\sn
\item[$(d)$]   For example, for a club of $\delta < \omega_1$,
$2^{\beth_\delta}$ suffices for the parallel of \ref{6.4}, 
\cite[5.13]{Sh:400}
\sn
\item[$(e)$]   if on $\lambda$ there is a $\aleph_1$-saturated 
$\lambda$-complete ideal extending $J^{\bd}_\lambda$ 
(for example, $\lambda$ a real valued measurable) \underline{then}  
$\cov(\lambda,\aleph_1,\aleph_1,2) \le \lambda$ \cite[\S3]{Sh:430} 
and more there
\sn
\item[$(f)$]   in clause (c), if $\kappa^{\aleph_0} < \lambda_0$ we can
add $\pp(\lambda_\delta) =^+ \pp_{J^{\bd}_\omega}(\lambda_\delta)$;
of course, there $\cov(\lambda_\delta,\lambda_\delta,\aleph_1,2) =
\cov(\lambda_\delta,\lambda_\delta,\kappa^+,2)$.
\end{enumerate}
\end{note}

\begin{note}
\label{7.6}
$\cov =$ minimal cardinality of a stationary $S$
\cite[3.6,5.12]{Sh:355}, \cite[3.6,3.8,3.8A,5.11,5.2A]{Sh:400}, 
\cite[Ch.VII,\S1,\S4]{Sh:g}, \cite[2.6(using 2.2),3.7]{Sh:410}, 
finally \cite[3.6]{Sh:420}; for example

\[
\cf({\cS}_{\le \kappa}(\lambda),\subseteq) = \min\{|S|:S \subseteq 
[\lambda]^{\le \kappa} \text{ is stationary}\}.
\]

\mn
Moreover, we got a measure one set of this cardinality for an appropriate 
filter; for another filter see \cite{Sh:580}.
\end{note}

\begin{note}
\label{7.7}
Covering by normal filters ({\rm prc}), \cite[\S4]{Sh:371},
\cite[\S1]{Sh:410}, generalization \cite[\S5]{Sh:410}, essentially
\cite[proof of 4.2 second case]{Sh:430}.  To quote \cite[\S1]{Sh:410}. 
\end{note}

\begin{note}
\label{7.8}
On $\cf_J(\prod {\ga},<_I)$, a generalization, 
see \cite[3.1]{Sh:400}.
\end{note}

\begin{note}
\label{7.9}
Computing $\cf^\sigma_{< \theta}(\prod {\ga})$ 
\cite[3.2]{Sh:400}; computing from it $\pp(\lambda)$ for non-fixed point
$\lambda$ by it \cite[3.3]{Sh:400}.
\end{note}

\begin{note}
\label{7.10}
$\cov$ is $\cf^\sigma_{< \theta} (\prod(\Reg \cap
\lambda),<_{J^{\bd}_\lambda})$ is $\pp_{\Gamma(\theta,\sigma)}$, 
when $\cf(\sigma) > \aleph_0$ \cite[3.3,3.4,3.5]{Sh:400}.
\end{note}

\begin{note}
\label{7.11}
$\cov(\lambda,\lambda,\cf(\lambda)^+,2) =^+
\pp(\lambda)$ when $\lambda$ is singular non-fixed point 
\cite[3.7(1), and more 3.7(1)-(5),3.8]{Sh:400}.
\end{note}

\begin{note}
\label{7.12}
Computing $\cov(\lambda,\theta,\theta,2)$ by using
$\cf_{< \theta}$ when $\theta > \cf(\lambda) = \aleph_0$, 
see
\cite[5.1,5.2,5.3,5.4,5.4A,5.5]{Sh:400},restriction to subset of $\lambda \cap
  \Reg$ is 
  \cite[5.5A]{Sh:400}.
   See more in \cite[6.5]{Sh:400}; here, in \ref{5.7}, \ref{5.8}.
\end{note}

\begin{note}
\label{7.13}
Finding a family ${\cP}$ of subsets of $\lambda$ covering
many of the countable subsets of $\lambda$, for example, if $a \in 
[\lambda]^{\aleph_1}$ we can find $H:a \rightarrow \omega$ such that each
countable subset of $H^{-1}(\{0,\dotsc,n\})$ is included in a member of
${\cP}$.  I.e. we characterize the minimal cardinality of such ${\cP}$
by $\pcf$ \cite[2.1-2.4]{Sh:410}, \cite[1.2]{Sh:430} more in \cite{Sh:513}.
\end{note}

\begin{note}
\label{7.14}
Characterizing the existence of ${\cP} \subseteq [\lambda]^{\aleph_1}$, 
$|{\cP}| > \lambda$ with pairwise finite intersection
\cite[\S6]{Sh:410} more in \cite[1.2]{Sh:430}, \cite{Sh:513}.
\end{note}

\begin{note}
\label{7.15}
If $\lambda \ge \mu > \sigma = \cf(\sigma) > \aleph_0$, \underline{then}  $\{\cov(\lambda,\mu,\theta,\sigma) : \mu \ge \theta > \sigma\}$ is finite \cite{Sh:412}.
\end{note}

\begin{note}
\label{7.16}
Let $\lambda > \kappa > \aleph_0$ be regular, then:
$\bigwedge\limits_{\mu < \lambda} \cov(\mu,\kappa,\kappa,2) < \lambda$
\underline{iff} for every $\mu < \lambda$ and 
$\langle a_\alpha : \alpha < \lambda \rangle$, 
$a_\alpha \subseteq \mu$, $|a_\alpha| < \kappa$ for some unbounded 
$s \subseteq \lambda$, $|\bigcup\limits_{a \in s} a_\alpha| < \kappa$ 
(a problem of Rubin-Shelah \cite{Sh:117}, see \cite[6.1]{Sh:371}, 
\cite[3.1]{Sh:430}).  For $\lambda$ successor of regular, a stronger 
theorem: see \cite[\S6]{Sh:371}; more \cite[6.13,6.14]{Sh:513}.
\end{note}

\begin{note}
\label{7.17}
If $\mu > \lambda \ge \kappa$, $\theta = \cov(\mu,\lambda^+,\lambda^+,\kappa)$ 
and $\cov(\lambda,\kappa,\kappa,2) \le \mu$ (or at least $\le \theta$), 
\underline{then} $\cov(\mu,\lambda^+,\lambda^+,2) = \cov(\theta,\kappa,\kappa,2)$, 
\cite[2.1]{Sh:430}.
\end{note}

\begin{note}
\label{7.18}
1) If $\lambda \ge \beth_\omega$, \underline{then} for some $\kappa < \beth_\omega$, $\cov(\lambda,\beth^+_\omega,\beth^+_\omega,\kappa) = \lambda$, \cite[1.1]{Sh:460}; 
any strong limit singular can serve instead of $\beth_\omega$.

\noindent
2) For a singular limit cardinal $\mu$ (for example $\mu =
\aleph_\omega$) sufficient conditions (for replacing $\beth_\omega$ by
$\mu$) are given in \cite[2.1,4.1]{Sh:460}.  For example such a condition is 
\mn
\begin{enumerate}
\item[$(*)_{\kappa,\mu}$]   ${\ga} \subseteq \Reg \setminus
\mu \wedge |{\ga}| < \mu \Rightarrow |\pcf_{\kappa-\mathrm{complete}}
({\ga})| < \mu$.
\end{enumerate}
\mn
3) So for every $\lambda \ge \beth_\omega$ for some $n$ and ${\cP} \subseteq
[\lambda]^{< \beth_\omega}$ of cardinality $\lambda$, 
every $X \in [\lambda]^{\le \beth_\omega}$ is
the union of $\le \beth_n$ sets from ${\cP}$; (\cite[2.5]{Sh:460}) and
the inverse \cite[4.2]{Sh:460} (see \cite{Sh:513}). 

Also if the statement above holds for e.g. $\aleph_\omega$ then
$(*)_{\aleph_n,\aleph_{\omega +1}}$ holds (by \cite[2.6]{Sh:460}).
\end{note}

\newpage

\section {Bounds in cardinal arithmetic} \label{A8}
\bigskip

\begin{note}
\label{8.1}
If $\langle \lambda_i : i \le \kappa \rangle$ is increasing continuous, 
$J$ a normal ideal on $\kappa$ and $\pp_J(\lambda_i) \le \lambda^{+h(i)}_i$, 
\underline{then} $\pp_J(\lambda_\kappa) \le \lambda^{\|h\|}_\kappa$ 
\cite[2.4]{Sh:355}, \cite[1.10,1.11]{Sh:371} where $\|h\|$ is Galvin Hajnal rank, i.e.

\[
\|h\| = \sup\big\{\|f\| + 1:f <_{D_\kappa} h \big\},
\]

\mn
$D_\kappa$ the club filter on $\kappa$.
\end{note}

\begin{note}
\label{8.2}
Let $C_0$ be the class of infinite cardinals and
define by induction:

\[
C_\zeta =: \{\lambda:\text{ for every } \xi < \zeta,\ \lambda \text{ is a fixed
point of } C_\zeta \}
\]

\mn
(i.e. $\lambda = \otp(C_\zeta \cap\lambda)$), \underline{then} for example

\[
\pp \big(\omega_1 \text{-th member of } C_1 \setminus \beth_2(\aleph_1) \big) <
\beth_2(\aleph_1)^+ \text{-th member of } C_1 \setminus \beth_2(\aleph_1)
\]

\mn
\cite[5.6]{Sh:386}.
\end{note}

\begin{note}
\label{8.3}
For $\zeta < \omega_1$ we have

\[
\pp_{\nor}(\aleph^\zeta_\omega(\beth_2(\aleph_1))) <
\aleph^\zeta_{(\beth_2(\aleph_1))^+}(\beth_2(\aleph_1))
\]

\mn
and more on $\aleph^\zeta_\delta$, see \cite[5.4,5.5]{Sh:386}, where

\[
\aleph^0_\alpha(\lambda) = \lambda^{+ \alpha},\ \aleph^{\zeta +1}_0(\lambda)
= \lambda,\ \aleph^{\zeta +1}_{\alpha +1}(\lambda) = \aleph^i_\zeta(\aleph_0)
\]

\mn
where

\[
\zeta = \aleph^{i+1}_\alpha(\lambda) + 1 \text{ and } 
\aleph^{\zeta +1}_\delta(\lambda) = \bigcup\limits_{\alpha < \delta}
\aleph^{i+1}_\alpha(\lambda),
\]

\mn
and for $i$ limit,

\[
\aleph^i_0(\lambda) = \lambda,\ \aleph^i_{\alpha +1}(\lambda) = 
\bigcup\limits_{j<i} \aleph^j_{\alpha +1}(\aleph^i_\alpha(\lambda)) 
\text{ and } \aleph^i_\delta(\lambda) = 
\bigcup\limits_{\alpha < \delta} \aleph^i_\alpha(\lambda).
\]
\end{note}

\begin{note}
\label{8.4}
If there are no [there are $\le \aleph_1$] inaccessibles
below $\lambda$, $\lambda > 2^{\aleph_1}$, $\cf(\lambda) = \aleph_1$,
 \underline{then}  there are no [there are $\le 2^{\aleph_1}$] inaccessibles below
$\pp(\lambda)$ \cite[5.10]{Sh:386}, similarly for Mahlo,
 $\epsilon$-Mahlo.
\end{note}

\begin{note}
\label{8.5}
If $\bigwedge\limits_{\delta < \omega_1} 
\pp(\aleph_\delta) < \aleph_{\omega_1},\pp(\aleph_{\omega_1}) 
= \aleph_{\alpha^*}$, \underline{then}  there are $|\alpha^*|$ subsets of 
$\omega_1$ with pairwise countable intersection 
\cite[1.7(1), more(2)]{Sh:371} getting Kurepa trees
\cite[2.8.2.9]{Sh:371}.
\end{note}

\begin{note}
\label{8.6}
The minimal counterexample to Tarski statement is simple,
Jech-Shelah \cite{Sh:385}. 

In \cite{Ta1} Tarski showed that for every limit ordinal $\beta,
\prod\limits_{\xi < \beta} \aleph_\xi = \aleph^{|\beta|}_\beta$, and
conjectured that

\[
\prod\limits_{\xi < \beta} \aleph_{\sigma_\xi} = \aleph^{|\beta|}_\alpha
\]

\mn
holds for every ordinal $\beta$ and every increasing sequence 
$\{\sigma_\xi\}_{\xi < \beta}$ such that $\lim_{\xi < \beta} 
\sigma_\xi = \aleph_\alpha$. 

Now: if a counterexample exists, then there exists one of length
$\omega_1 + \omega$ (Jech and Shelah \cite{Sh:385}).
\end{note}

\begin{note}
\label{8.7}
$\pp(\aleph_{\alpha + \delta}) < \aleph_{\alpha + 
|\delta|^{+4}}$ \cite[2.1,2.2, more 2.3-2.8]{Sh:400}.
\end{note}

\begin{note}
\label{8.8}
If $\delta < \aleph_4,\cf(\delta) = \aleph_0$ then
$\pp(\aleph_\delta) < \aleph_{\omega_4}$.   If $|\delta| + 
\cf(\delta)^{+3} < \kappa$, then $\pp(\aleph_{\alpha + \delta}) <
\aleph_{\alpha + \kappa}$ \cite[4.2,4.3,4.4]{Sh:400}, more
\cite[3.3-3.6]{Sh:410}.
\end{note}

\begin{note}
\label{8.9}
More on the number of inaccessibles: \cite[\S4]{Sh:430}.

E.g. \cite[4.4]{Sh:430}. For transparency, assume that for no core model $\bfK[A]$,
with $A$ a set of ordinals, do we have covering (here the SCH holds). Then
\begin{enumerate}
    \item Assume $\mu > \cf(\mu) = \aleph_1$, $\mu_0 < \mu$, 
    $$\sigma \geq \big|\{ \lambda \in (\mu_0,\mu) : \lambda \text{ inaccessible} \}\big| < \mu.$$
    Then
    $$ \sigma^{+4} > \big| \{ \lambda : \mu < \lambda < 
    \pp_{\Gamma(\sigma,\aleph_1)}(\mu),\ \lambda \text{ inaccessible} \} \big| $$
    
    \item The parallel of \cite[4.3]{Sh:400}.
\end{enumerate}
\end{note}

\begin{note}
\label{8.10}
By Gitik and Shelah \cite{Sh:412}:
\mn
\begin{enumerate}
\item[$(a)$]  If $\mu$ is a J\'onsson limit cardinal not strong limit,
\underline{then}  $\langle 2^\sigma:\sigma < \mu \rangle$ is eventually constant.
\sn
\item[$(b)$]  If $\mu$ is a limit cardinal, $\mu_0 < \mu$ and
$\bigwedge\limits_{\theta \in (\mu_0,\mu)} \mu \rightarrow 
[\theta]^{< \omega}_{\theta,\mu_0}$, \underline{then}  $\langle 2^\theta:
\mu_0 < \theta < \mu \rangle$ has finitely many values.
\sn
\item[$(c)$]   If on $\mu$ there is a $\mu^+_0$-saturated, uniform
$\mu$-complete ideal for example $\mu$ a real value measurable
$\le 2^{\aleph_0}$, \underline{then}  the assumption of (b) holds, hence its
conclusion.
\end{enumerate}
\end{note}
\newpage

\section {J\'onsson algebras} \label{A9}
\bigskip

\begin{note}
\label{9.1}
Definition and previously known results: 
\cite[4.3,4.4]{Sh:355}.  A J\'onsson algebra is one with no proper 
subalgebra with the same cardinality.  A J\'onsson cardinal is $\lambda$ 
such that there is no J\'onsson
algebra with countable vocabulary and cardinality $\lambda$.
\end{note}

\begin{note}
\label{9.2}
Definition of $\id_j(\bar C)$, $\id^j_\theta(\bar C)$ see
\cite[1.8]{Sh:380}, $\id^j_J(\bar C)$ see \cite[1.16]{Sh:380} (also with $k$
instead of $j$).
\end{note}

\begin{note}
\label{9.3}
J\'onsson games: Definition \cite[2.1]{Sh:380}, connection to
\cite[2.3]{Sh:380} (for example $\lambda = \aleph_{\omega +1}$).
\end{note}

\begin{note}
\label{9.4}
$\lambda^+$ (for a singular $\lambda$) is not a J\'onsson cardinal when:
\mn
\begin{enumerate}
\item[$(a)$]   $\lambda$ is not an accumulation point of inaccessible
J\'onsson cardinals \cite[4.5 more 4.6]{Sh:355}
\sn
\item[$(b)$]   weaker hypothesis (for $\lambda^+ \rightarrow
[\lambda^+]^{< \omega}_\kappa$) \cite[2.5]{Sh:413}
\sn
\item[$(c)$]  $\lambda = \beth^+_\omega$ (see \cite{Sh:413},
\cite{Sh:535} more there)
\sn
\item[$(d)$]   on every large enough regular $\mu < \lambda$, there is an
algebra $M$ on $\mu$ which has no proper subalgebra with set of
elements a stationary subset of $\mu$, see \cite[3.3]{Sh:572}.
\end{enumerate}
\end{note}

\begin{note}
\label{9.5}
Sufficient condition for ``$\lambda$ not J\'onsson"
\cite[1.8,1.9]{Sh:365} for $\lambda \nrightarrow [\lambda]^{< \omega}_\sigma$
\cite[1.10,3.5,3.6,3.7]{Sh:365}.
\end{note}

\begin{note}
\label{9.6}
$\lambda$ inaccessible is not J\'onsson when: $\lambda$
not Mahlo \cite[3.8]{Sh:365}, $\lambda$ has a stationary subset $S$ not
reflecting in inaccessibles \cite[3.9]{Sh:365}, $\lambda$ not $\lambda$-Mahlo
\cite{Sh:380}, $\lambda$ not $\lambda \times \omega$-Mahlo 
\cite[1.14]{Sh:413}, there is a set $S$ of singulars satisfying, 
$\rk_\lambda(S) > \rk_\lambda(S^+)$ where $S^+ = \{\kappa < \lambda:
\kappa$ inaccessible, $S \cap \kappa$ stationary$\}$, \cite[1.15]{Sh:413}.
\end{note}

\begin{note}
\label{9.7}
If $\mu^+$ is a J\'onsson cardinal, $\mu > \cf(\mu) >
\aleph_0$, \underline{then}  $\cf(\mu)$ is ``almost" $\mu^+$-supercompact
\cite[2.8]{Sh:413} other \cite[2.10]{Sh:413}.
\end{note}

\begin{note}
\label{9.8}
If $\lambda$ is regular, and for every regular large enough
$\mu < \lambda$, for some $f:\mu \rightarrow \lambda$ we have
$\|f\|_{J^{\bd}_\mu} \ge \lambda$ (or at least this holds for ``enough"
$\mu$'s), \underline{then}  on $\lambda$ there is a J\'onsson algebra,
\cite[2.12+2.12A]{Sh:380}.  More sufficient conditions there.
\end{note}

\begin{note}
\label{9.9}
See more \cite{Sh:413}, \cite{Sh:535}.
\end{note}
\newpage

\section {Colouring $=$ negative partition relations:\\ (see \cite{Sh:282}, \cite{Sh:280},
\cite{Sh:327})} \label{A10}
\bigskip

\begin{note}
\label{10.1}
Definition of $\Pr_\ell$: $\Pr_0$, see \cite[AP,1.1]{Sh:g}, 
$\Pr^{(-)}_1$, see \cite[AP,1.2]{Sh:g}, $\Pr^{(y)}_2$, see 
\cite[AP,1.3]{Sh:g}, $\Pr^{(y)}_3$, see \cite[AP,1.4]{Sh:g}, 
$\Pr_4$, see \cite[4.3]{Sh:365}.  

For example: 
$\Pr_1(\lambda,\mu,\theta,\kappa)$ means: there is a 2-coloring of $\lambda$
by $\theta$ colours ($=$ symmetric 2-place function from $\lambda$ to
$\theta$) such that: if $\langle w_i:i < \mu \rangle$ is a sequence of
pairwise disjoint subsets of $\lambda,\bigwedge\limits_i |w_i| < \kappa$ and
$\zeta < \theta$, \underline{then}  for some $i < j$, on $w_i \times w_j$ the coloring
$c$ is constant.  In $\Pr_0(\lambda,\mu,\theta,\kappa)$ we replace $\zeta$
by $h:\kappa \times \kappa \rightarrow \theta$ and demand $\alpha 
\in w_i \text{ and } \beta \in w_j \Rightarrow c(\alpha,\beta) = 
h(\otp(w_i \cap \alpha),\otp(w_j \cap \beta))$.  If $\mu = 
\lambda$ we may omit it, if $\kappa = \aleph_0$ we may omit it.  
(See \cite[AP,1.2]{Sh:g}).
\end{note}

\begin{note}
\label{10.2}
Trivial implications \cite[AP,1.6,1.6A,1.7]{Sh:g} and
$\Pr_1 \Rightarrow \Pr_0$ by \cite[4.5(3)]{Sh:365}, $\Pr_4 \Rightarrow
\Pr_1 \Rightarrow \Pr_0$ by \cite[4.5(1)]{Sh:365}.  
For example if $\Pr_1(\lambda,\mu,\theta,\sigma),
\chi = \chi^{< \sigma} + 2^\theta \le \mu \le \lambda < 2^\chi$ then
$\Pr_0(\lambda,\mu,\chi,\sigma)$.  Other such $\Pr$ and implications
\cite[\S2,\S4]{Sh:572}.
\end{note}

\begin{note}
\label{10.3}
Colouring for successor of singular: \cite[4.1,4.7]{Sh:355},
\cite[\S2]{Sh:413} for example $\Pr_1(\lambda^+,\lambda^+,
(\cf(\lambda))^+,2)$ for $\lambda$ singular.
\end{note}

\begin{note}
\label{10.4}
Combining $\Pr_\ell$'s \cite[4.8,4.8A]{Sh:355}.
\end{note}

\begin{note}
\label{10.5}
Using $\pcf$: 
\mn
\begin{enumerate}
\item[$(a)$]   if $\lambda = \tcf(\prod {\gc}/J^{\bd}_{\gc})$
and $\big[\theta \in {\gc} \Rightarrow |{\gc} \setminus \theta| = |{\gc}| \big]$, 
\underline{then}  $\Pr_1(\lambda,\lambda,2^{|{\gc}|},\cf({\gc}))$; 
see \cite[4.1B]{Sh:355}, \cite{Sh:282}.
\sn
\item[$(b)$]   getting colouring on $\lambda \in \pcf({\ga})$
from colourings on every $\theta \in {\ga}$, see \cite[4.1D]{Sh:355}.
\end{enumerate}
\end{note}

\begin{note}
\label{10.6}
Using guessing of clubs: Definition and basic properties
of for example $(Dx)^\lambda_{\kappa,\sigma,\theta,\tau}$
\cite[4.1]{Sh:365}.
\end{note}

\begin{note}
\label{10.7}
Proof of such properties \cite[4.2]{Sh:365}, \cite[2.6]{Sh:413}
\mn
\begin{enumerate}
\item[$(a)$]   if $\lambda$ is a regular $\lambda > \sigma > \kappa$
\underline{then}  $\Pr_1(\lambda^+,\lambda^+,\sigma,\kappa)$, \cite[\S4]{Sh:365}
\sn
\item[$(b)$]   if $\lambda$ is inaccessible with a stationary subset $S$
not reflecting in inaccessibles and $\sigma \le \text{ min}_{\delta \in S}
\cf(\delta)$ and $\kappa < \lambda$ then $\Pr_1(\lambda,\lambda,\kappa,
\sigma)$, \cite[4.1+4.7]{Sh:365}
\sn
\item[$(c)$]   if $\lambda = \mu^+,\mu > 2^{\cf(\mu)},\kappa < \mu$,
\underline{then}  $\Pr_1(\lambda,\lambda,\cf(\mu),\kappa)$, \cite[2.7]{Sh:413}
\sn
\item[$(d)$]   if $\lambda = \mu^+,\mu > \cf(\mu)$ \underline{then} 
$\Pr_1(\lambda,\lambda,\cf(\mu),\cf(\mu))$, \cite[4.1(1)]{Sh:355}
\sn
\item[$(e)$]   by \cite{Sh:535} we get such 
properties for e.g. $\lambda = \beth^+_\omega$
\sn
\item[$(f)$]   if $\lambda = \aleph_2 \text{ and } \mu = \aleph_0$ or if
$\lambda = \mu^{++}$, $\mu$ regular then $\Pr_1(\lambda,\lambda,\lambda,\mu)$
(\cite[\S1]{Sh:572}). 
\end{enumerate}
\end{note}

\begin{note}
\label{10.8}
$(E2)$ implies $\Pr_4$ \cite[4.4]{Sh:365}.
\end{note}

\begin{note}
\label{10.9}
$(D2) \Rightarrow \Pr_1$ \cite[4.7]{Sh:365}.
\end{note}

\begin{note}
\label{10.10}
Concerning the results in \cite{Sh:95} on partition
relations restriction of the kind appearing there are necessary (we use
FILL) see, some day \cite{Sh:F50}.
\end{note}

\begin{note}
\label{10.11}
Galvin conjecture:
\mn
\begin{enumerate}
\item[$(a)$]   $\aleph_n \nrightarrow [\aleph_1]^{n+1}_{\aleph_0}$
(\cite[5.8(1)]{Sh:288}, more there), but
\sn
\item[$(b)$]  for the naturally defined $h:\omega \rightarrow \omega$
if $\CON(\ZFC + \lambda \rightarrow (\aleph_1)^\omega_2$) \underline{then} 
it is consistent with $\ZFC$ that: $2^{\aleph_0} = \lambda \rightarrow 
[\aleph_1]^n_{h(n)}$, (we can even get $X \in [\lambda]^{\aleph_1}$
which exhibits the conclusion simultaneously for all $n$, 
$\lambda \rightarrow [\aleph_1]^n_{h_1(n)}$, if $h_1(n) \ge n$, 
$h_1(n)/h(n) \rightarrow \infty$), \cite[3.1]{Sh:288}
\sn
\item[$(c)$]   if $\kappa$ is measurable indestructible by adding (even
many) Cohen subsets to $\kappa$, then a generalization of Halpern Lauchli 
theorem holds to ${}^{\kappa >}2$ (but using some 
$\langle <_\alpha : \alpha < \kappa \rangle$, $<_\alpha$ a well order of 
${}^\alpha 2$) (\cite[4.1,4.2 + \S2]{Sh:288}).
See more in \cite{Sh:481}, \cite{Sh:546}, \cite{Sh:585}, and \cite{Sh:1176}. 
\end{enumerate}
\end{note}

\begin{note}
\label{10.12}
More on colouring (improving results on J\'onssonness from
\cite{Sh:413} to colouring) see \cite{Sh:535}, e.g. for $\lambda = \beth^+_\omega$ 
we have Pr$_1(\lambda,\lambda,\lambda,\aleph_0)$.
\end{note}

\begin{note}
\label{10.13}
More on $\Pr_i$'s in \cite[\S3]{Sh:829}. E.g., if $\mu
> \aleph_0$ is strong limit, $\chi \ge \mu$, $\lambda = 2^\chi$ is
singular, then $\chi \in \mu \cap \Reg \setminus \{\aleph_0\}
\Rightarrow \mathrm{Ps}_1(\cf(\lambda),\lambda,\kappa)$.

See more in \cite{Sh:1163} and from there, on \cite{Ri14}, \cite{Sh:1027}.
\end{note}

\newpage

\section {Trees and linear orders} \label{A11}
\bigskip

\begin{note}
\label{11.1}
Let $\ga_\delta = \{\lambda_i:i < \delta\}$ and $\ga_i
= \{\lambda_j:j < i\}$ for every $i < \delta$.
If $\lambda = \max \pcf(\ga_\delta)$ and
$\lambda_i > \max \pcf(\ga_i)$, \underline{then}  we can find in
$\prod(\ga_\delta)$ a $<_{J_{< \lambda}[(\ga_\delta)]}$-increasing cofinal 
sequence $\langle f_\alpha:\alpha < \lambda \rangle$ such that 
$\big\{\{f_\alpha \rest \ga_j:j < i\}:i < \delta,\alpha < 
\lambda \big\}$ forms a tree with $\delta$ levels, level $i$ of cardinality
$\max \pcf(\ga_j) < \lambda_i$ and $\ge \lambda$ many $\delta$-branches
\cite[3.5]{Sh:355}. 

Note:
\mn
\begin{enumerate}
\item[$(a)$]   The lexicographic order on 
${\cF} = \{f_\alpha : \alpha < \lambda\}$ has density 
$\sum\limits_{i < \delta} \lambda_i$.
\sn
\item[$(b)$]   If $\prod \lambda_i/I$ is as in 
\cite[1.4(1)(see 1.3)]{Sh:355} then $F$ is $(\Sigma \lambda_i)^+$-free 
(see \ref{6.3}). Hence any set of cardinality $\le \Sigma \lambda_i$
is the union of $\le$ gen$(I)$ many sets $F'$ each satisfying ``for some 
$s \in I$ we have $\langle f_\alpha \rest (\delta \setminus s) : 
f_\alpha \in F'\rangle$ is increasing": i.e. 
$$\big[\alpha < \beta,\ f_\alpha, f_\beta \in F',\ i \in \delta \setminus s \big] 
\Rightarrow f_\alpha(i) < f_\beta(i).$$ 
Here $\mathrm{gen}(I) = \min\{|{\cP}|:{\cP} \subseteq I 
\text{ generates the ideal } I\}$. 
\cite[1.4(3)]{Sh:355}
\sn
\item[$(c)$]   if $\lambda > 2^{|\delta|}$, \underline{then}  we can have such
trees with exactly $\lambda$ branches \cite{Sh:276}; somewhat more: 
\cite[6.6B]{Sh:430}.
\end{enumerate}
\mn
See more in part (C).
\end{note}

\begin{note}
\label{11.2}
There are quite many $\langle \lambda_i:i < \delta \rangle$, 
$\lambda$ as in \ref{11.1}: for example, if 
$$\aleph_0 < \kappa = \cf(\mu) < \mu_0  < \mu < \lambda = \cf(\lambda) < \pp_\kappa(\mu)$$ 
\underline{then}  we can find such $\langle \lambda_i:i < \kappa \rangle$ 
with limit $\mu$ with $\mu_0 < \lambda_i < \mu$, if 
$\bigwedge\limits_{\alpha < \mu}|\alpha|^\kappa < \mu$ or at least 
$(\forall \mu' < \mu)[\pp_\kappa(\mu') < \mu]$, see
\cite[1.6(2),(4)]{Sh:371}.  Also $\pp(\aleph_{\alpha + \delta}) < 
\alpha_{\alpha + |\delta|^+}$ helps to get such examples, see 
\cite[\S5]{Sh:462}, \cite{Sh:534}.
\end{note}

\begin{note}
\label{11.3}
For $\lambda > \kappa = \cf(\kappa)$ the following cardinals are equal:

\[
\sup\{\mu:\text{some tree with } \lambda \text{ nodes has} \ge \mu
\text{ many } \kappa \text{-branches}\}
\]

\mn
and

\begin{equation*}
\begin{array}{clcr}
\sup\{\pcf({\ga}):\,&|{\ga}| < \min({\ga}),\ \cf(\otp({\ga})) = 
\kappa,\ {\ga} \subseteq \Reg \cap \lambda^+ \setminus \kappa \text{ and} \\
  &\,\theta \in {\ga} \Rightarrow \max \pcf({\ga} \cap \theta) < \theta\}  
\end{array}
\end{equation*}

\mn
see \cite[2.2]{Sh:589}.
\end{note}

\begin{note}
\label{11.4}
Definition of Ens, entangled linear order and basic facts. 
(Ens stands for \emph{entangled sequence}.)
See for example \cite[AP,2.1,2.2 more 2.3]{Sh:g}.

A linear order ${\clI}$ is $\lambda$-entangled if given any $n < \omega$
and pairwise distinct $x^e_\zeta \in {\clI}$ $(e < n,\zeta < \lambda)$
and $w \subseteq \{0,1,\dotsc,n-1\}$ there are $\zeta < \xi$ such that
for $e<n$ we have: $x^e_\zeta < x^e_\xi \Leftrightarrow e \in w$.  We say
${\clI}$ is entangled if it is $|{\clI}|$-entangled; $\Ens(\lambda,\mu)$
means there are $\mu$ linear orders ${\clI}_\zeta$ $(\zeta < \mu)$ each 
of cardinality $\lambda$ and if $n < \omega$, $\zeta_e < \mu$ distinct 
$(e < n)$ and $w \subseteq n$ and if $x^e_\zeta \in {\clI}_\zeta$ are 
distinct then for some $\alpha < \beta < \mu$ we have 
${\clI}_{\zeta_e} \models x^e_\alpha < x^e_\beta \Leftrightarrow e \in w$.

For more on $\sigma$-entangled linear orders see \cite{Sh:462}; first
\end{note}

\begin{claim}\label{11.49}
$\Ens(\cf(2^{\aleph_0}))$, by Bonnet-Shelah \cite{Sh:210}. The proof gives\\
$\Ens(\cf(2^\lambda),\aleph_0)$ when there is a linear order of cardinality 
$2^\lambda$ and density $\lambda$.
\end{claim}

\begin{note}
\label{11.5}
$\Ens(\lambda^+,\cf(\lambda))$ for $\lambda$ singular
\cite[4.9 more 4.11,4.14]{Sh:355} more \cite[5.3]{Sh:371}.
\end{note}

\begin{note}
\label{11.6}
For $\mu$ regular uncountable and a linear order ${\clI}$
of power $\mu$, $\clI$ is entangled iff the interval Boolean algebra of
${\clI}$ is $\lambda$-narrow (see Bonnet-Shelah \cite{Sh:210} 
(using different names), later \cite[2.3]{Sh:345b} or \cite[\S1]{Sh:462}).
\end{note}

\begin{note}
\label{11.7}
A sufficient condition for existence of entangled linear
order of cardinality $\lambda$ is: $\lambda = \max \pcf({\ga})$, $\kappa = |{\ga}|$, 
$[\theta \in {\ga} \Rightarrow \theta > \max\pcf(\theta \cap {\ga})]$, 
$2^\kappa \ge \sup({\ga})$, $\ga$ divisible to $\kappa$ sets not in 
$J_{< \lambda}[{\ga}]$, \cite[4.12]{Sh:355}. If we omit ``$2^\kappa \ge \sup(\ga)$" 
we can still prove $\Ens(\lambda,\kappa)$; \cite[4.10A]{Sh:355} more in 
\cite[4.10F,4.10G]{Sh:355}, \cite[5.4,5.5,5.5A]{Sh:371}.
\end{note}

\begin{note}
\label{11.8}
If $\cf(\lambda) < \lambda \le 2^{\aleph_0}$, \underline{then} 
there is an entangled linear order in $\lambda^+$,
\cite[4.13]{Sh:355}.
\end{note}

\begin{note}
\label{11.9}
If $\lambda \in \pcf({\ga})$ and 
$\big[\theta \in {\ga} \Rightarrow \theta > \max \pcf(\theta \cap {\ga})\big]$ 
and for each $\theta \in {\ga}$ there is an entangled linear order or just
$\Ens(\theta,\max \pcf(\theta \cap {\ga}))$, \underline{then} there is one 
on $\lambda$, \cite[4.10C]{Sh:355}.
\end{note}

\begin{note}
\label{11.10}
\mn
\begin{enumerate}
    \item[$(a)$]    If $\kappa^{+4} \le \cf(\lambda) < \lambda < 2^\kappa$, 
    \underline{then} there is an entangled linear order in $\lambda^+$,
    \cite[4.1 more 4.2,4.3]{Sh:410}.
\sn
    \item[$(b)$]   There is a class of cardinals $\lambda$ for which there is
    an entangled linear order of cardinality $\lambda^+$,
    \cite[\S5]{Sh:371}. 
    It is not clear if we can demand e.g. $\lambda = \lambda^{\aleph_0}$, but
    if this fails, then for $\kappa$ large enough, $\kappa^{\aleph_0} = \kappa
    \Rightarrow 2^\kappa < \aleph_{\kappa^{+4}}$ (see (a), more in \cite{Sh:462}).
\sn
    \item[$(c)$]  There is a class of cardinals $\lambda$ for which there is
    a Boolean algebra $B$ of cardinality $\lambda^+$ with neither chain nor
    antichain of cardinality $\lambda^+$; i.e. if $Y \subseteq B$, $|Y| = |B|$ then
    $(\exists x,y \in Y)[x < y]$ and $(\exists x,y \in Y)[x \nleq y \wedge y \nleq x]$. 
    In fact, for any sequence $\langle x_\alpha : \alpha < \lambda^+\rangle$ 
    of distinct members of $B$:
\sn
    \begin{enumerate}
        \item[$(i)$]  $(\exists \alpha < \beta)(x_\alpha < x_\beta)$,
\sn
        \item[$(ii)$]    $(\exists \alpha < \beta)(x_\alpha > x_\beta)$ and
\sn
        \item[$(iii)$]   $(\exists \alpha < \beta)[x_\alpha \nleq x_\beta \wedge
        x_\beta \nleq x_\alpha]$; see \cite[4.3]{Sh:462}.
    \end{enumerate}
\sn
    \item[$(d)$]   Moreover, in part (c), for any given $\lambda_0$, letting
    $\mu$ be the minimal $\mu = \aleph_\mu > \lambda_0$ then we can find $B$ as
    there with density $\mu$ (everywhere); similarly in (b).
\sn
    \item[$(e)$]   Moreover in (c) (and (b)) if the density character is
    $\mu$, $\ell \in \{0,1,2\}$, $\theta = \cf(\theta) < \mu$ and $x_\alpha \in B$ 
    (for $\alpha < \lambda$) are distinct then for some $w \subseteq \lambda$, 
    $|w| = \theta$ we have for any $\alpha,\beta \in w$, $\alpha < \beta$:

    \[
    \ell = 0 \Rightarrow x_\alpha < x_\beta 
    \]

    \[
    \ell = 1 \Rightarrow x_\alpha > x_\beta
    \]

    \[
    \ell = 2 \Rightarrow x_\alpha \nleq x_\beta \wedge x_\beta \nleq x_\alpha
    \]

    Similarly in part (b).
\sn
    \item[$(f)$]  If $2^\lambda$ is singular,  then there is an entangled linear
    order of cardinality $(2^\lambda)^+$ (the assumption implies 
    $(\forall \mu)[\lambda < \cf(\mu) < \mu \le 2^\lambda < \pp(\mu)$] 
    (i.e. $\mu = 2^\lambda$), this suffices as we can use \ref{6.6}(b), 
    \ref{6.12}, \ref{11.7}; see \cite[5.5, pg.65]{Sh:462}).
\end{enumerate}
\end{note}

\begin{note}
\label{11.11}
Universal linear orders: see Section 13, Model Theory.
\end{note}

\begin{note}
\label{11.12}
For every $\lambda$ there is $\mu,\lambda \le \mu <
2^\lambda$ such that $(A)$ or $(B)$:
\mn
\begin{enumerate}
\item[$(A)$]   $\mu = \lambda$ and for every regular $\chi \le 2^\lambda$
there is a tree $T$ of cardinality $\lambda$ with $\ge \chi$ branches (so a
linear order of cardinality $\ge \chi$ and density $\le \lambda$)
\sn
\item[$(B)$]   $\mu > \lambda$, and:
\sn
\begin{enumerate}
\item[$(\alpha)$]   $\pp(\mu) = 2^\lambda$, $\cf(\mu) \le \lambda$, 
$(\forall \theta)[\cf(\theta) \le \lambda < \theta < \mu \Rightarrow
\pp_\lambda(\theta) < \mu]$.  Hence, by \cite[\S1]{Sh:371} for every 
regular $\chi \le 2^\lambda$ there is a tree from \cite[3.5]{Sh:355}: 
$\cf(\mu)$ levels, every level of cardinality $< \mu$ and $\chi$ 
($\cf(\mu))$-branches
\sn
\item[$(\beta)$]   for every $\chi \in (\lambda,\mu)$, there is a tree 
$T$ of cardinality $\lambda$ with $\ge \chi$-branches of the same height
\sn
\item[$(\gamma)$]   $\cf(\mu) = \cf(\lambda_0)$ for $\lambda_0 =
\min\{\theta:2^\theta = 2^\lambda\}$ and even
$\pp_{\Gamma(\cf(\mu))}(\mu) = 2^\lambda$ see \cite[5.11]{Sh:355}, 
\cite[4.3]{Sh:410} and \cite[3.3]{Sh:430}; see more in \cite[2.10]{Sh:600}.
\end{enumerate}
\end{enumerate}
\end{note}

\begin{note}
\label{11.13}
If $\theta_{n+1} = \min\{\theta:2^\theta > 
2^{\theta_n}\}$ for $n < \omega,
\sum\limits_{n < \omega} \theta_n < 2^{\theta_0}$, \underline{then} 
for some $n > 0$ and regular $\mu \in [\theta_n,\theta_{n+1})$ for every
regular $\chi \le 2^{\theta_n}$, there is a tree with $\mu$ nodes and $\ge
\chi \, \mu$-branches \cite[3.4]{Sh:430}.
\end{note}

\begin{note}
\label{11.14}
Kurepa trees: there are two contexts that arise
\mn
\begin{enumerate}
\item[$(a)$]   we can get Kurepa trees of singular cardinality: if
$\bar\lambda = \langle \lambda_i:i < \delta\rangle$ and
$\delta < \lambda_i = \cf(\lambda_i)$, $\lambda_i > \max \pcf\{\lambda_j:j < i\}$ 
\underline{then} there is a tree with 
$\delta$ levels, the $i$-th level of cardinality $< \lambda_i$, 
and at least $\max \pcf\{\lambda_i:i < \delta\} \, \delta$-branches, see
\cite[3.5]{Sh:355}, hence can derive consequences
from conventional cardinal arithmetic assumptions
\sn
\item[$(b)$]   if for example $\pp(\aleph_{\omega_1}) > \aleph_{\omega_2}$
and for a club of $\delta < \omega_1,\pp(\aleph_\delta) < \aleph_{\omega_1}$,
\underline{then}  there is an $(\aleph_1)$-Kurepa tree (see \cite[2.8]{Sh:371} for
more).  We get a large family of sets with small intersection in more
general circumstances \cite[1.7]{Sh:371}.
\sn
\item[$(c)$]   If $\bigwedge\limits_{\alpha < \mu} |\alpha|^\kappa < \mu$,
cf$(\mu) = \kappa > \aleph_0$ and $\mu \le \lambda < \mu^\kappa$,
\underline{then}  there is a tree with $\mu$ nodes, $\kappa$ levels and 
exactly $\lambda$ branches, $\lambda$ of them of height $\kappa$.  
We can derive results on linear orders (really they are the same 
problems).  If we speak on the number of $\kappa$-branches (or for 
linear order number of Dedekind cuts of cofinality from at least one 
side $\kappa$), instead of ``$\bigwedge\limits_{\alpha < \mu} 
|\alpha|^\kappa < \mu$" it suffices that
\sn
\begin{enumerate}
\item[$(*)$]  $(a) \quad 2^\kappa < \mu_0 < \mu$
\sn
\item[${{}}$]   $(b)\quad$ if $\mu_0 < \chi < \mu$ and 
$\cf(\chi) \le \kappa$ then $\pp(\kappa) < \mu$.
\end{enumerate}
\end{enumerate}
\mn
See \cite{Sh:262} or \cite[6.6(1)]{Sh:430}.  (Similarly, other results can be
translated between trees and linear orders).
\end{note}
\newpage

\section {Boolean Algebras and General Topology} \label{A12}
\bigskip

\begin{note}
\label{11.15}
Concerning Boolean algebras and topology.  $\lambda$-c.c. is not 
productive and $\lambda-L$-spaces exist and $\lambda-S$-spaces exist
and more follows from $\Pr^-_1(\lambda,2)$
(or appropriate colouring) see \cite[1.6A]{Sh:282a} so
\cite[4.2]{Sh:355} is a conclusion of this.  This is translated to
results on cellularity of topological spaces 
(cellularity $\le \lambda \Leftrightarrow \lambda^+$-c.c.).

We have
\mn
\begin{enumerate}
\item[$(a)$]     if $\lambda \ge \aleph_1$, for some $\lambda^+$-c.c.
Boolean algebras $B_1,B_2$ we have: $B_1 \times B_2$ is not $\lambda^+$-c.c.
(why?  now Pr$_1(\lambda^+,\lambda^+,2,\aleph_0)$ suffice
\cite[1.6A]{Sh:g} and it holds by \cite{Sh:327} or 
\cite[4.8(1),p.177]{Sh:365} if $\lambda$ regular $> \aleph_1$, 
\cite[4.1,p.67]{Sh:355} if $\lambda$ is singular and lastly by 
\cite[\S1]{Sh:572} if $\lambda = \aleph_1$)
\sn
\item[$(b)$]   if $\lambda$ is inaccessible and has a stationary subset not
reflecting in any accessible, \underline{then}  for some $\lambda$-c.c. Boolean
algebras $B_1,B_2$ we have: $B_1 \times B_2$ is not $\lambda$-c.c. (see
\cite[4.8]{Sh:365})
\sn
\item[$(c)$]   if $\lambda$ is Mahlo, $\otimes^{\aleph_0}_\lambda$ (see
\ref{1.4}) \underline{then}  for some $\lambda$-c.c. Boolean algebras $B_n$, 
for any proper filter $I$ on $\omega$ extending $J^{\bd}_\omega$ we have 
$\prod B_n/I$ fails the $\lambda$-c.c. \cite[4.11]{Sh:365}.
\end{enumerate}
\end{note}

\begin{note}
\label{11.16}
Concerning Topology: characterizing by pp when there are
$f_\alpha \in {}^\kappa \sigma$ for $\alpha < \theta$ such that $\alpha <
\beta \Rightarrow \bigvee\limits_{i < \kappa} f_\alpha(i) < f_\beta(i)$, see
\cite[3.7]{Sh:410}, is needed for Gerlits, Hajnal and Szentmiklossy 
\cite{GHS}.
The condition is (when $\theta$ is regular for simplicity)

\[
2^\kappa \ge \theta \text{ or } (\exists \mu)\big[\cf(\mu) \le
\kappa < \mu \text{ and  }\mu \le \sigma \text{ and } \pp^+(\mu) > \theta\big]
\]

\mn
(for $\theta$ singular, just ask if it holds for every regular $\theta_1 < \theta$).

(Why not just $\theta \le \sigma^\kappa$?  Because if e.g. 
$\kappa = \beth_\omega$, $\theta = \beth_{\omega +1}$, $\sigma = \aleph_0$ 
we do not know whether $\pp^+(\kappa) = \theta^+$).
\end{note}

\begin{note}
\label{11.17}
Concerning Topology: let $X$ be a topological space, $B$ a basis of
the topology (not assuming the space to be Hausdorff or even $T_0$).  If
$\lambda = \aleph_0$ or $\lambda$ is strong limit of cofinality $\aleph_0$,
and the number of open sets is $> |B| + \lambda$, \underline{then}  it is $\ge \lambda
^{\aleph_0}$; see for $\lambda = \aleph_0$ \cite{Sh:454}, for $\lambda >
\aleph_0$, \cite{Sh:454A} relaying on \cite{Sh:460}.
\end{note}

\begin{note}
\label{11.18}
Concerning Topology: densities of box products: for example if $\mu$
is strong limit singular, $\mu = \sum\limits_{i < \cf(\mu)} \lambda_i$,
$\cf(\lambda_i) = \aleph_0$, $2^{\lambda_i} = \lambda^+_i$, $\lambda_i$ 
are strong limit cardinals, $\max \pcf\{\lambda_i : i < \cf(\mu)\} < 2^\mu$,
$\cf(\mu) < \theta < \mu$, \underline{then}  the density of the 
$\cf(\mu)^+$-box product ${}^\theta\!\mu$ is $2^\mu$ \cite[\S5]{Sh:430}.

Gitik Shelah \cite{Sh:597} prove consistency results.
\end{note}

\begin{note}
\label{11.19}
\mn
\begin{enumerate}
\item[$(a)$]   the results in \ref{11.17} come from starting to analyze
the following: 
given a $\mu_i$-complete filter $D_i$ on $\lambda_i$, for $i < \kappa$,
what is

\begin{equation*}
\begin{array}{clcr}
\min\big\{|A| : \,&A \subseteq \prod\limits_{i < \kappa} \lambda_i
\text{ such that for every} \\
  &\,\langle A_i:i < \kappa \rangle \in \prod\limits_{i < \kappa} D_i
\text{ we have } A \cap \prod\limits_{i < \kappa} A_i \ne \varnothing\big\}
\end{array}
\end{equation*}

continued in \cite{Sh:575} and then \cite{Sh:620}.
\sn
\item[$(b)$]   This is applied also to the problem of $\lambda$-Gross
spaces. These are vector spaces $V$ over a field $F$ with an inner product 
such that for any subspace $U \subseteq V$ of dimension $\lambda$,
\[
\dim\Big\{x \in V : \bigwedge\limits_{y \in U} (x,y) = 0\Big\} < \dim V.
\]
See Shelah Spinas \cite{Sh:468}.
\end{enumerate}
\end{note}

\begin{note}
\label{11.20}
A well known problem in general topology is whether every
Hausdorff space can be divided to two sets each not containing a homeomorphic
copy of Cantor's discontinuum.  In \cite{Sh:460} we have a sufficient
condition for this (e.g. $|{\ga}| \le \aleph_0 \Rightarrow |\pcf
({\ga})| \le \aleph_0$ and $2^{\aleph_0} \ge \aleph_\omega$, by
\cite[3.6(2)]{Sh:460}, the $(*)_1$ version relying on
\cite[Th.2.6]{Sh:460}). 
But we can prove: if $c \ell$ is a closure operation on ${\clP}(X)$ 
(i.e. $a \subseteq c \ell(a) = c \ell(c \ell(a))$, 
$a \subseteq b \Rightarrow c \ell(a) \subseteq c \ell(b))$ and $|a| \ge 
\aleph_0 \Rightarrow |c \ell(a)| > \beth_\omega$, \underline{then}  we can
partition $X$ to two sets, each not containing any infinite $a = c \ell(a)$.
(Can prove more).

Related weaker problem is to find large $A \subseteq {}^\omega \lambda$
containing no large closed subsets, $\lambda$ strong limit of cofinality
$\aleph_0$, if $\pp(\lambda) = 2^\lambda$ is easy (hence for higher 
cofinalities this holds and e.g. for many $\beth_\delta$, 
$\delta < \omega_1$).  
See \cite[6.9]{Sh:355}, more in \cite[3.3,3.4]{Sh:430}.  See more in
\cite{Sh:460}, \cite{Sh:668}.
\end{note}

\begin{note}
\label{11.21}
If $\pp(\lambda) > \lambda^+$, $\cf(\lambda) = \aleph_0$
(or just a consequence from \cite[\S1]{Sh:355}, see \ref{6.2} here),
\underline{then}  there is first countable $\lambda$-collectionwise Hausdorff (and
even $\lambda$-metrizable), not $\lambda^+$-collectionwise Hausdorff space
(see \cite{Sh:E9}; when we assume just $\cov(\lambda,\lambda,\aleph_1,2) >
\lambda^+$ use \cite[\S6]{Sh:355}).
\end{note}

\begin{note}
\label{11.22}
If $\lambda < \lambda^{< \lambda}$, then there is a
regular $\kappa < \lambda$ and tree $T$ with $\kappa$ levels, for each
$\alpha < \kappa$, $T$ has $< \lambda$ members of level $\le \alpha$, and $T$
has $> \lambda$ many $\kappa$-branches.  If $\lambda < \lambda^{< \lambda}$ and
$\neg(\exists \mu)\big[\mu \text{ strong limit and } \mu \le \lambda < 2^\mu\big]$,
\underline{then} this is above $2^\kappa > \lambda$; see \cite[6.3]{Sh:430}.
\end{note}

\begin{note}
\label{11.23}
Depth of homomorphic images of ultraproducts of Boolean
algebras, \cite[\S3]{Sh:506} and resolved for 
$\lambda_i > 2^{|\Dom(D)|}$ in \cite[\S3]{Sh:589}.
\end{note}

\begin{note}
\label{11.24}
If $\lambda$ is strong limit singular, $\kappa =
\cf(\lambda)$ and e.g. $2^\lambda = \lambda^+$, \underline{then}  for some
Boolean algebras $B_1,B_2$ we have: $B_1$ is $\lambda^+$-c.c., $B_2$ is
$(2^\kappa)^+$-c.c. but $B_1 \times B_2$ is not $\lambda^+$-c.c.
(see for more \cite{Sh:575}).  More constructions in \cite{Sh:620}.
\end{note}

\begin{note}
\label{11.25}
If $\lambda = \lambda^{\beth_\omega}$, $B$ a
$\beth_\omega$-c.c. Boolean algebra of cardinality $\le 2^\lambda$
\underline{then}  $B$ is $\lambda$-linked (that is $B \setminus \{0\}$ is 
the union of $\le \lambda$ sets of pairwise non-disjoint elements), 
see \cite[\S8]{Sh:575}.
\end{note}

\begin{note}
\label{11.26}
On the measure algebra, \cite{Sh:620}.
\end{note}

\begin{note}
\label{11.27}
On independent sets in Boolean Algebra, \cite{Sh:620}.
\end{note}

\begin{note}
\label{11.28}
On ultraproducts of Boolean Algebra: $s(B)$, spread, i.e.
constructing examples of inv$(\prod\limits_{i < \kappa} B_i/D) >
\prod\limits_{i < \kappa} \mathrm{inv}(B_i)/D$, see:
\mn
\begin{enumerate}
\item[$(a)$]   for inv being $s$, (spread), Roslanowski Shelah
\cite{Sh:534}, \cite{Sh:620}
\sn
\item[$(b)$]   similarly hd (hereditarily density)
\sn
\item[$(c)$]   similarly hL (hereditarily Lindelof)
\sn
\item[$(d)$]   for inv being Depth, \cite{Sh:641}, \cite{Sh:853},
  \cite{Sh:878}, \cite{Sh:956}
\sn
\item[$(e)$]   for inv being Length, \cite{Sh:641}.
\end{enumerate}
\end{note}
\newpage

\section {Strong covering, forcing, and partition calculus} \label{A13}
\bigskip

\begin{note}
\label{12.1}
Preservation under forcing: essentially $\pcf$ and $\pp$ are
preserved except for forcing notion involving large cardinals.  Specifically
if (the pair of universes) $(V,W)$ satisfies $\kappa$-covering [i.e.
$V \subseteq W$ and if $a \subseteq \Ord,W \models |a| < \kappa$ then
for some $b \in V,a \subseteq b \subseteq \Ord$ and $W \models |b| <
\kappa]$ and ${\ga} \subseteq \Ord \setminus \kappa$ is a set
from $W$ of cardinality $< \kappa$ of regulars of $W$ \underline{then} 

\[
\pcf^V\big\{\cf^V(\theta):\theta \in \pcf^V({\ga}) \big\} =
\big\{\cf^V(\lambda):\lambda \in \pcf^W(\{\cf^W(\theta):\theta \in {\ga}\})\big\}
\]

\mn
(this applies for example to $(K,V) $ 
 if there is no inner model with
measurable by Dodd and Jensen \cite{DJ1}).
\end{note}

\begin{note}
\label{12.2}
The strong covering lemma: see \cite[Ch.XIII,\S1,\S2]{Sh:f}
or better \cite[Ch.VII,\S1,\S2]{Sh:g}; see more in 
\cite[2.6,p.407]{Sh:410} and \cite{Sh:580}, each can be read independently.

Suppose $W \subseteq V$ is a transitive class of $V$ including all the
ordinals and is a model of $\ZFC$, let $\lambda > \kappa$ be cardinals of $V$.

We say $(W,V)$ satisfies the strong $(\lambda,\kappa)$-covering property if
for every model $M \in V$ with universe $\lambda$ and predicates and
function symbols there is $N \prec M$ of cardinality $< \kappa$, 
$N \cap \kappa \in \kappa$, $N \in V$ but the universe of $N$ belongs to $W$; 
we also use stronger versions (like the set of such $N$'s is positive 
or even equal to $[\lambda]^{\le \kappa}$ modulo some ideal, \underline{or} 
weaker versions like union of few sets from $W$). 

Those papers do this without using fine structure assumptions, just that
$(W,V)$ satisfies $(\lambda,\kappa)$-covering and related properties. 

\end{note}

\begin{note}
\label{12.3}
Application of ranks (see \ref{3.2}) to partition
calculus: Shelah Stanley \cite{Sh:419}.  If there is a nice filter of
$\kappa$ (see \ref{3.2}) and $\lambda,\cf(\lambda) > \kappa = 
\cf(\kappa),(\forall \mu < \lambda)\,\mu^\kappa < \lambda$ \underline{then} 
$\lambda \rightarrow (\lambda,\omega + 1)^2$.
\end{note}

\begin{note}
\label{12.11}
Also, with ranks.  If $\lambda > \cf(\lambda) 
> \aleph_0$, \underline{then}  $\lambda \rightarrow (\lambda,\omega
+1)^2$ in ZFC (see \cite{Sh:881}).
\end{note}

\begin{note}
\label{12.14}
\underline{Polarized Partition Relations} 

If $\lambda$ is strong limit singular and
$2^\lambda > \lambda^+$, \underline{then}  $\binom{\lambda^+}\lambda \rightarrow
\binom \lambda \lambda^{1,1}_2$, see \cite{Sh:586}.
\end{note}

\begin{note}
\label{12.17}
If $\lambda > \cf(\lambda)$ is a limit of measurables
and some $\pcf$ assumptions are forced, then even $\binom
{\lambda^+}{\lambda} \rightarrow \binom{\lambda^+}{\lambda}^{1,1}_2$,
see \cite{Sh:949}.
\end{note}

\begin{note}
\label{12.20}
See \cite{Sh:497}. 
\end{note}

\begin{note}
\label{12.29}
On inner models, see \cite{Sh:805}. 
\end{note}
\newpage

\section{Axiom of Choice, weak versions} \label{A14}

\begin{note}
\label{12a.2}
Set theory with weak choice: \cite{Sh:497}, \cite{Sh:938}, \cite{Sh:955}.
\end{note}

\begin{note}
\label{12a.5}
The intermediate axiom of choice: $(\mathrm{Ax})_4$ (see \cite{Sh:835}).
\begin{enumerate}[1)]
    \item We suggest considering $(\mathrm{Ax})_4$ in addition to ZF + DC, which tells us each $[\alpha]^{\aleph_0}$ is well ordered. This is orthogonal to $\bfV = \bfL[\bbR]$.
    
    \item In particular, it gives: given a set $I$, for every ordinal $\alpha$, ${}^I\!\alpha$ is covered by a sequence of few well ordered sets (depending only on $I$), even uniformly.
    
    \item There is a class of successor cardinals which are regular (and even so called ``explicitly successor'').
\end{enumerate}
\end{note}

\begin{note}\label{12a.8}
(See \cite{Sh:1005}) We have a quite strong pcf theorem which the minimal cofinality is big enough compared to the index set.
\end{note}

\begin{note}
\label{12a.11}
(See \cite{Sh:1005}) Black Boxes and constructing abelian groups.
\end{note}

\begin{note}
\label{12a.14}
On splitting stationary sets, see Larson and Shelah \cite{Sh:925}.
\end{note}

\begin{note}
\label{12a.17}
A ZFC conclusion --- we quote from \cite[\S3]{Sh:835}:

\noindent We prove that if $\mu > \kappa = \cf(\mu) > \aleph_0$, then from a well-ordering of $\clP(\clP(\kappa)) \cup {}^{\kappa>}\!\mu$ we can define a well-ordering of ${}^\kappa\!\mu$. If, e.g., $\mu$ is a strong limit singular cardinal of uncountable cofinality,  then by using a well order of $\clH(\mu)$ we can define a well-ordering of $\clP(\mu)$, hence of $\clH(\mu^+)$. Lastly, we give sufficient conditions (in ZF + DC) on singular $\mu$ for $\mu^+$ to be regular. Actually, if $\mu = \mu^{\aleph_0} + 2^{2^\kappa}$, $\kappa = \kappa^{\aleph_0}$, and $X \subseteq \mu$ codes $\clP(\clP(\kappa))$ and ${}^\omega\!\mu$, then by using $X$ as a parameter we can define a well-ordering of ${}^\kappa\!\mu$.
\end{note}

\newpage

\section {Transversals and $(\lambda,I,J)$-sequences} \label{A15}
\bigskip

See \cite{Sh:161} (and \cite{Sh:52}), a transversal is a one to one choice
function.

\begin{note}
\label{13.1}
If $I$ is an ideal on $\kappa,\lambda > \cf(\lambda)$
and $\pp_I(\lambda) > \mu$, \underline{then}  we can find a family of functions
$f_\alpha(\alpha < \mu)$ from $\kappa$ to $\lambda$, which is 
$\lambda^+$-free for $I$ i.e. any $\lambda$ of them are strictly 
increasing on each $x \in \Dom(I)$ if for each $\alpha$ we ignore a 
set $s_\alpha \in I$ such that $i \in \kappa 
\setminus (s_\alpha \cup s_\beta) \Rightarrow
f_\alpha(i) < f_\beta(i)$ (so $\{\Rang(f_\alpha):\alpha \in u\}$ has a
transversal when $u \subseteq \mu,|u| \le \lambda$) \cite[1.5A]{Sh:355} (the
case $\mu$ singular changes nothing for this purpose).  So $\NPT(\lambda^+,
\kappa)$ (see Definition below).  On weakening
``$\pp_I(\lambda) > \mu$" to ``$\pp^+_I(\lambda) > \mu$" for $\mu$ successor
of regular see \cite[\S6]{Sh:371} ($\mu$ singular-easy).  On weakening
$\pp_I(\lambda) > \mu$ to $\cov(\lambda,\lambda,\kappa^+,2) > \mu$, see
\cite[\S6]{Sh:355} for some variants; in particular
$\NPT_{J^{\bd}_\omega}(\lambda^+,\aleph_0)$ when 
$\cov(\lambda,\lambda,\kappa^+,2) > \lambda$ by
\cite[6.3,p.99]{Sh:355}.
\end{note}

\begin{note}
\label{13.2}
Definitions of variants of $\NPT$, \cite[6.1]{Sh:355},
\cite[6.3]{Sh:371} for example $\NPT(\lambda,\kappa)$ means that there is a
family $\{A_i:i < \lambda\}$ of sets each of cardinality $\le \kappa$, and
$< \lambda$ of them have a transversal, but not all.  Similarly for
$\NPT_J(\lambda,\kappa)$ we have $f_\alpha:\Dom(J) \rightarrow$ ordinals
as in \ref{13.1}.
\end{note}

\begin{note}
\label{13.3}
Trivial and easy facts \cite[6.2,6.7]{Sh:355}, why 
concentrating on\\ ``$\NPT(\lambda^+,\aleph_1)$, $\cf(\lambda) = \aleph_0$"
\cite[6.4]{Sh:355}.
\end{note}

\begin{note}
\label{13.4}
If $\lambda > \cf(\lambda) = \aleph_0$ and
$\cov(\lambda,\lambda,\aleph_1,2) > \lambda^+$, \underline{then} 
$\NPT_{J^{\bd}_\omega}(\lambda^+,\aleph_1)$, \cite[6.3]{Sh:355} 
more in \cite[6.5,6.8]{Sh:355}, \cite[6.1]{Sh:371}, \cite[6.2]{Sh:371}
application to \cite{Sh:117}, \cite[6.4,6.5]{Sh:371}.
\end{note}

\begin{note}
\label{13.5}
When $\lambda$ is a strong limit of cofinality $\aleph_0$,
there is $T \subseteq {}^\omega \lambda,|T| = 2^\lambda$ with no large dense
subset, \cite[6.9]{Sh:355} (there is a subclaim with more information).
\end{note}

\begin{note}
\label{13.6}
If $I$ is an ideal on $\kappa$, $\mu > \kappa \ge \cf(\mu)$, 
$\pp^+_I(\lambda) > \lambda = \cf(\lambda) > \mu$, $\lambda_i = \cf(\lambda_i) > \kappa$ 
(for $i < \kappa$), $\tlim_I \lambda_i = \mu$, $\langle f_i:i < \lambda \rangle$ 
is $<_I$-increasing cofinal in $\prod\limits_{i < \kappa} \lambda_i/I$, 
\underline{then} for some $A \subseteq \lambda$, $|A| = \lambda$ for every 
$B \subseteq A$ of cardinality $\lambda$
and $\delta < \mu^+$ there is $B' \subseteq A$ of order type $\delta$ and
$\langle s_\alpha:\alpha \in B' \rangle$ such that: $s_\alpha \in I$, $\alpha <
\beta$ and $\zeta \in \kappa \setminus (s_\alpha \cup s_\beta)$ and $\alpha \in B'$ and 
$\beta \in B' \Rightarrow f_\alpha(\zeta) < f_\beta(\zeta)$
so $\big\langle \Rang(f_\alpha \rest (\kappa \setminus s_\alpha):
\alpha \in B' \big\rangle$ is a sequence of pairwise disjoint sets. (For somewhat
more, see \cite[6.2,6.2A(3)]{Sh:430}).
\end{note}

\begin{note}
\label{13.7}
($\kappa$-MAD families)

Let $\kappa = \cf(\kappa) > \aleph_0$.  For any $\mu \ge 2^\kappa$,
letting $$\chi = \chi^\kappa_\mu = \sup\{\pp_{J^{\bd}_\kappa}(\mu'):
2^\lambda \le \mu' \subseteq \mu,\ \cf(\mu') = \kappa\}$$ we have:
\mn
\begin{enumerate}
\item[$(a)$]   every $\kappa$-almost disjoint subfamily of $[\mu]^\kappa$
(i.e. intersection of two has cardinality $< \kappa$) has cardinality
$\le \chi$; also $\chi^\kappa_\mu = T_{J^{\bd}_\kappa}(\mu)$
\sn
\item[$(b)$]   trivially there is maximal $\kappa$-almost disjoint family
$\subseteq [\mu]^\kappa$ and all such families have the same cardinality
which is in $\chi$
\sn
\item[$(c)$]   if $\chi_0 = \chi,\chi_{n+1} = \chi^\kappa_{(\chi_n)},
\chi_\omega = \sum\limits_{n < \omega} \chi_n$ \underline{then} 
\sn
\sn
\begin{enumerate}
\item[$(\alpha)$]  $\chi_n = \sup\{\pp_{J_n}(\mu'):2^\lambda \le \mu' 
\le \mu,\cf(\mu') = \kappa\}$ where $$J_n = \{A \subseteq \kappa^n:
(\exists^{< \kappa} \alpha_0)(\exists^{< \kappa} \alpha_1) \ldots
(\exists^{< \kappa} \alpha_{n-1})\langle \alpha_0,\dotsc,\alpha_{n-1} \rangle
\in A\}$$
\sn
\item[$(\beta)$]   $\chi^\kappa_{(\chi_\omega)} = \chi_\omega$ 
(hence is doubtful if it is consistent to have $\chi_n \ne \chi_{n+1}$).
\end{enumerate}
\end{enumerate}
\end{note}

\begin{note}
\label{13.8}
$\bar \eta = \langle \eta_\alpha:\alpha < \lambda \rangle$
is a $(\lambda,I,J)$-sequence for $\bar I = \langle I_i:i < \delta \rangle$
iff each $\eta_\alpha \in \prod\limits_{i < \delta} \Dom(I_i),J$ is
an ideal on $\delta,I$ is an ideal on $\lambda$, each $I_i$ is an ideal on
$\Dom(I_i)$, and

\[
X \in I^+ \Rightarrow \big\{i < \delta:\{\eta_\alpha(i):\alpha \in X\} \in I_i \big\}
\in J.
\]

\mn
The definition was introduced in \cite{Sh:575} and considered again in
\cite{Sh:620}.  In \cite{Sh:620} first the case of the Erd\"os-Rado ideal 
defined there was considered.  For the case $\bar I = 
\langle J^{\bd}_{\lambda_i}:i < \delta \rangle$ and $\lambda = 
\tcf(\prod\limits_{i \le \delta} \lambda_i/J^{\bd}_\delta)$ and
$J = J^{\bd}_\lambda,I = J^{\bd}_\delta$, the existence of a 
$(\lambda,I,J)$-sequence comes from pcf theory.  Also the case $I_i = 
\prod\limits_{\ell < n_i} J^{\bd}_{\lambda_{i,\ell}}$ for 
$\langle \lambda_{i,\ell}:\ell < n_i \rangle$ increasing a sufficient $\pcf$ 
condition for the existence of a $(\lambda,I,J)$-sequence was given in 
\cite{Sh:620} which holds sometimes (for any given $\langle n_i:i < \delta
\rangle$.  Also in \cite{Sh:620} the case 
$I_i = J^{\bd}_{\langle \lambda_{i,\ell}:\ell < n \rangle}$, for 
$\langle \lambda_{i,\ell}:\ell < n \rangle$ a decreasing sequence
of regulars was considered, giving a sufficient condition which requires $\pcf$
to be reasonably complicated.  A most case $I_i = \prod\limits_{\ell
< n} J^{\mathrm{nst},\theta}_{\lambda_{i,\ell}},\lambda_{i,\ell}$ regular
decreasing, $J^{\mathrm{nst},\theta}_\lambda$ is the ideal of non-stationary
sets $+ \{\delta < \lambda:\cf(\delta) \ne \theta\}$, when e.g.
$\delta < \kappa < \lambda_{i,\ell}$ and we prove existence for some
$\langle \lambda_{i,\ell}:\ell < n,i < \delta \rangle$.  Many applications
for Boolean algebras can be found in \cite{Sh:620}.
\end{note}

\begin{note}
\label{13.9}
The family $\{\kappa:\NPT(\kappa,\aleph_1)\}$
is not too small, see \cite{Sh:108}, Magidor Shelah \cite{Sh:204}, 
\cite{Sh:523}.
\end{note}
\newpage

\section {Model Theory, algebra, and Black Boxes} \label{A16}
\bigskip

\begin{note}
\label{14.1}
$L_{\infty,\lambda}$-equivalent non-isomorphic models
in $\lambda$: if $\lambda > \cf(\lambda) > \aleph_1$ there are such
models of cardinality $\lambda$ (if $\cf(\lambda) = \aleph_1$, it suffices to
have: $\big\langle \lambda_i : i < \cf(\lambda) \big\rangle$ is an
increasing sequence of regulars with limit $\lambda$ and that
$$\big\{\delta < \cf(\lambda):\text{there is an unbounded } a \subseteq \delta \text{ with }
\lambda > \max \pcf\{\lambda_i:i \in a\}\big\}$$ is stationary; not known
if this fails in some universe of set theory, see \cite[\S7]{Sh:355}.
\end{note}

\begin{note}
\label{14.2}
\underline{Universal Models}: for example, the class of
\underline{linear orders}.  

If $\lambda$ is regular and $\exists \mu\ [\mu^+ < \lambda < 2^\mu]$, 
\underline{then} there is in $\lambda$ no universal linear order, not
even a universal model (for elementary embeddings) for $T$ in $\lambda$ where
$T$ is a first order theory with the strict order property.  For almost all
singular $\lambda$ we have those results, more specifically if $\lambda$
is not a fixed point of the second order the result holds; and if it fails
for $\lambda$ the consequences for pp are not known to be consistent, see
\cite{Sh:409} which rely on guessing clubs.
\end{note}

\begin{note}
\label{14.3}
A much weaker demand on the first order $T$ suffices
in \ref{14.2}: $\mathrm{NSOP}_4$, see \cite[\S2]{Sh:500}. On the remaining cardinals
see some information in \cite[\S3]{Sh:457}; on complimentary consistency 
(only for $\lambda = \aleph_1$) see \cite[\S4]{Sh:100}.
\end{note}

\begin{note}
\label{14.4}
\underline{Universal models for $(\omega +1)$-trees with
$(\omega +1)$-levels and or stable} 

$\qquad$ \underline{unsuperstable $T$}:

Similar results: if $\lambda$ regular $(\exists \mu)[\mu^+ < \lambda <
\mu^{\aleph_0}]$ then there is no universal member; also for most singular
\cite{Sh:447}. 

Similarly if $\kappa = \cf(\kappa) < \kappa(T)$,
$(\exists \mu)[\mu^+ < \lambda < \mu^\kappa]$.
\end{note}

\begin{note}
\label{14.5}
Universal abelian groups have similar results for pure
embedding (under reasonable restrictions (mainly the groups are reduced,
because there are divisible universal abelian groups the interesting 
cardinals are $\lambda^{\aleph_0} > \lambda > 2^{\aleph_0}$).  For torsion
free reduced abelian groups, ${\gK}^{\mathrm{rtf}}$, or reduced separable
$p$-groups, ${\gK}^{\rs(p)}$ if $2^{\aleph_0} + \mu^+ < \lambda =
\cf(\lambda) < \mu^{\aleph_0}$, \underline{then}  there is no universal.  For
``most" $\lambda$, $\lambda$ regular can be omitted.

(This and more \cite{Sh:455}).
\end{note}

\begin{note}
\label{14.6}
We can use the usual embedding but restrict the class 
of abelian groups.  The natural classes: ${\gK}^{\mathrm{rtf}}$ 
(torsion free, reduced i.e. has no divisible subgroups) and 
${\gK}^{\rs}(p)$ (reduced separable $p$-groups).  But in addition we 
restrict ourselves to the abelian groups which are $(< \lambda)$-stable 
(see \cite{Sh:456}; club guessing is used).
\end{note}

\begin{note}
\label{14.7}
For classes ${\gK}^{\mathrm{rtf}}, {\gK}^{\rs(p)}$ from
\ref{13.7} of abelian groups under embeddings see
\cite{Sh:552}: mainly if $\lambda^{\aleph_0} > \lambda > 2^{\aleph_0}$ there
are negative results except when some pcf phenomena not known to 
be consistent (also club guessing is used).  Below the continuum there 
are independence results.  More on the existence of universals see 
\cite{Sh:457} on metric spaces see \cite{Sh:552} and on normed 
spaces \cite{Sh:614}.
\end{note}

\begin{note}
\label{14.8}
For cardinals $\ge \beth_\omega$, for the classes
${\gK}^{\mathrm{rtf}}, {\gK}^{\rs(p)}$ the results in \ref{14.7} are
improved to have demands on the cardinals like \ref{14.4}, see
\cite{Sh:622}.
\end{note}

\begin{note}
\label{14.9}
There exists a reflexive abelian group, whose
cardinality is the first measurable cardinal, \cite{Sh:904}. 
We do not succeed in proving the existence of arbitrarily large reflexive groups, but show that it would follow from relatively weak assumptions. Furthermore, we show that any condition which would preclude their existence must be quite stringent.
\end{note}

\begin{note}\label{14.32}
For a survey on the existence of universal abelian groups under various embedding relations, see \cite[\S10]{Sh:1151} with references and table (see \cite[pg. 302]{Sh:1151}.)
\end{note}

\begin{note}\label{14.35}
\begin{enumerate}[1)]
    \item Existence of a $\lambda$-free but not free abelian group of cardinality $\lambda$ for arbitrarily large $\lambda <$ [first fixed point]. (See Magidor and the author, \cite{Sh:204}.)
    
    \item Consistency of ``for the first $\lambda = \aleph_\lambda$, every $\lambda$-free algebra is free'' for any variety; also \cite{Sh:204}.
    
    \item See more in \ref{14.44}, \cite{Sh:1028}.
\end{enumerate}
\end{note}

\begin{note}
\label{14.38}
\underline{Diamonds and Omitting Types}:  In the omitting type
theorem for $\bbL(Q)$ in the $\lambda^+$ interpretation, not only $\lambda =
\lambda^{<\lambda}$ (needed even for the completeness) was used in
\cite{Sh:82} but $(D \ell)_\lambda$ [for $\lambda$ successor this 
is $\diamondsuit_\lambda$, generally it means: there is 
$\langle {\cP}_\alpha:\alpha < \lambda \rangle$, ${\cP}_\alpha$ a 
family of $< \lambda$ subsets of $\lambda$ such
that for every $A \subseteq \lambda$ for stationarily many $\delta < \lambda$,
$A \cap \delta \in {\cP}_\delta$].  Now by \cite{Sh:460}: if $\lambda > \beth_\omega$ 
then $\lambda = \lambda^{< \lambda} \Leftrightarrow (D \ell)_\lambda$.  
In fact: if $\lambda = \lambda^{< \lambda}$ and 
$(\forall \mu < \lambda)(\mu^{<\kappa>_{\tr}} < \lambda) 
\Rightarrow (D \ell)_{S^\lambda_\kappa}$
(where $(D \ell)_{S^\lambda_\kappa}$ is defined as above but for 
$\alpha \in S^\lambda_\kappa =: \{\delta < \lambda:\cf(\delta) = \kappa\}$, and
$\mu^{<\kappa>_{\tr}} = \sup\{\lambda:\text{there is a tree with } \mu 
\text{-nodes and } \lambda \text{ many } \kappa \text{-branches}\}$). 
\end{note}

\begin{note}
\label{14.41}
There are uses for proving Black Boxes (see 
\cite[Ch.III,\S6]{Sh:e}), those are construction principles 
provable in ZFC, and have quite many applications, see there for references.
\end{note}

\begin{note}
\label{14.44}
\begin{enumerate}[1)]
    \item In \cite{Sh:898} we prove\footnote{The proof of \cite[2.6, Case 2]{Sh:898} appears incomplete, but the claim is proved in \cite[1.11]{Sh:775}.} that for strong limit singular $\mu$ we can find quite large and quite free sets $\subseteq {}^{\cf(\mu)} \! \mu$ and quite strong Black Boxes.
    
    \item We define the following objects as in \cite[\S1]{Sh:898}: 
    \begin{enumerate}
        \item Let $\bfC = \{ \text{strong limit singular } \mu : \pp(\mu) =^+ 2 ^\mu \}$, with $=^+$ as on pg.\pageref{A}
        
        \item $\bfC_\kappa = \{ \mu \in \bfC : \cf(\mu) = \kappa \}$ 
        
        \item The set $\clF \subseteq {}^\kappa\!\mu$ is called $(\theta, \sigma, J)$-free, where $J$ is an ideal on $\kappa$, \underline{when} $$ f_1 \neq f_2 \in \clF \Rightarrow \{i < \kappa : f_1(i) = f_2(i) \} \in J$$
        and every $\clF' \subseteq \clF$ of cardinality $< \theta$ is $[J, \sigma]$-free, which means that:
        \begin{itemize}
            \item there is a sequence $\LL u_f : f \in \clF' \RR$ of members of $J$ such that for every pair $(\gamma, i) \in \mu \times \kappa$, the set $\{ f \in \clF' : f(i) = \gamma \wedge i \notin u_f \}$ has cardinality $< 1 + \sigma$.
        \end{itemize}
        
        \item We may replace ``$\clF \subseteq {}^\kappa\!\mu$" by a sequence 
        $\bar C = \LL C_\delta : \delta \in S \RR$ with $C_\delta$ a set of order type 
        $\kappa$, or even a by a set $\{ C_\delta : \delta \in S \}$. This means that the 
        definition applies to $\{ f_\delta : \delta \in S \}$, where $f_\delta$ is an 
        increasing function $\kappa \to C_\delta$ for each $\delta$; similarly for the other parts.
    \end{enumerate}
\end{enumerate}
\end{note}

\begin{definition}\label{14.47}
[See \cite[0.5=LOp.14]{Sh:898}]

Assume we are given a quadruple $(\lambda,\mu, \theta, \kappa)$ of cardinals.\footnote{but we may replace $\lambda$ by an ideal $I$ on $S \subseteq \lambda = \sup(S)$, so writing $\lambda$ would mean that $S=\lambda$; also, we may replace $\kappa$ by an ideal $J$ on $\kappa$, so writing $\kappa$ would mean that $J = J_\kappa^\bd$.} Let $\BB^-(\lambda,\mu, \theta, \kappa)$ mean that some pair $(\bar C, \bar \bfc)$ satisfies clauses (A) and (B) below; we call the pair 
$(\bar C, \bar \bfc)$ a \emph{witness} for $\BB^-(\lambda,\mu, \theta, \kappa)$. Let $\BB(\lambda,\mu, \theta, \kappa)$ mean that some witness $(\bar C, \bar \bfc)$ 
satisfies clause (A) below and for some sequence $\LL S_i : i < \lambda \RR$ of pairwise disjoint subsets of $\lambda$ (or of $S$), each $(\bar C \rest S_i, \bar \bfc \rest S_i)$ satisfies clause (B) below,\footnote{thus replacing $S$ and $\bar \bfc$ by $S_i$ and $\bar \bfc \rest S_i$} where:
\begin{enumerate}
    \item 
    \begin{enumerate}
        \item $\bar C = \LL C_\alpha : \alpha \in S \RR$ and 
        $S = S(\bar C) \subseteq \lambda = \sup(S)$
    
        \item $C_\alpha \subseteq \alpha$ has order type $\kappa$
    
        \item $\bar C$ is $\mu$-free (see \ref{0.4}) [\underline{but} when we replace $\kappa$ by $J$ then we say ``$\bar C$ is $(\mu,J)$-free".]
    \end{enumerate}
    
    \item    
    \begin{enumerate}
        \item $\bar \bfc =  \LL \bfc_\alpha : \alpha \in S \RR$
        
        \item $\bfc_\alpha$ is a function from $C_\alpha$ to $\theta$
    
        \item if $\bfc : \bigcup\limits_{\alpha \in S} C_\alpha$, then $\bfc_\alpha = \bfc \rest C_\alpha$ for each $\alpha \in S$ 
        
        [\underline{but} when we replace $\lambda$ by an ideal $I$ on $S$, then we demand that the set $\{ \alpha \in S : \bfc_\alpha = \bfc \rest C_\alpha \}$ is not in $I$].
    \end{enumerate}
\end{enumerate}
\end{definition}

\begin{remark}\label{14.47a}
The reader may recall that if $S$ is a stationary subset of
$$\{ \delta < \lambda : \cf(\delta) = \kappa \}$$ for a regular cardinal $\lambda$, 
$S$ is non-reflecting, and $\bar C = \LL C_\alpha : \alpha \in S \RR$ satisfies 
$C_\delta \subseteq \delta - \sup(C_\delta)$ and $\otp(C_\delta) = \kappa$, 
\underline{then} $\diamondsuit_S$ implies $\BB(\lambda, \lambda, \lambda, \kappa)$. Therefore, if $\bfV = \bfL$ then for every regular $\kappa < \lambda$ with $\lambda$ a non-weakly compact cardinal we have $\BB(\lambda, \lambda, \lambda, \kappa)$.
\end{remark}

\begin{theorem}\label{14.50}
(From \cite[1.18-LOp.15]{Sh:750}.) We have $\BB(\lambda, \bar C, (\lambda,\theta), {<}\mu)$ \underline{when}:
\begin{enumerate}
    \item $\mu \in \bfC_\kappa$, $\lambda = \cf(2^\mu)$ and $\theta < \mu$, 
    $\sigma = \cf(\sigma) < \mu$;
    
    \item $S \subseteq S_\sigma^\lambda$ is stationary;
    
    \item $\bar C = \LL C_\alpha : \alpha \in S \RR$, $C_\delta \subseteq \delta$, 
    $|C_\delta| \leq \mu$;
    
    \item $\chi < 2^\mu \Rightarrow \chi^{\LL \sigma \RR_\tr} < 2^\mu$;
    
    \item $\bar C$ is \emph{shallow}: that is, 
    $\big| \{ C_\delta \cap \alpha : \alpha \in C_\delta \} \big| < \lambda$ for 
    $\alpha < \lambda$.
\end{enumerate}
\end{theorem}

\noindent\textbf{The BB Trichotomy Theorem \ref{14.47}.}\label{14.53} If $\mu \in \bfC_\kappa$ and $\kappa > \sigma = \cf(\sigma)$, \underline{then} at least one of the following holds:
\begin{enumerate}
    \item There is a $\mu^+$-free $\clF \subseteq {}^\kappa \! \mu$ of cardinality $2^\mu$
    
    \item 
    \begin{enumerate}
        \item $\lambda := 2^\mu = \lambda^{< \lambda}$ (so $\lambda$ is regular) and 
        $\chi < \lambda \Rightarrow \chi^\sigma < \lambda$
        
        \item if $S \subseteq S_\sigma^\lambda$ is stationary and
        $\bar C = \LL C_\alpha : \alpha \in S \RR$ is a weak ladder system (i.e., $C_\delta \subseteq \delta$)\footnote{Bear in mind that a choice of $C_\delta = \delta$ would satisfy this} \underline{then}
    
        \item letting $J_S^{\mathrm{nst}} = \{ A \subseteq \lambda : A \cap S \text{ is not stationary in } \lambda \}$ we have\footnote{What about freeness? We may get it by the choice of $\bar C$; also, if $\bar C$ is a ladder system (particularly if \emph{strictly}), we will get a weak form (e.g. stability).}
        \begin{enumerate}
            \item $\BB(J_S^{\mathrm{nst}}, \bar C, \theta,{\leq}\mu)$ for every 
            $\theta < \mu$ provided that\\ $\delta \in S \Rightarrow |C_\delta| < \mu$
        
            \item $\BB(J_S^{\mathrm{nst}}, \bar C, (2^\mu,\theta),{<}\lambda)$ for any 
            $\theta < \mu$
        \end{enumerate}
    \end{enumerate}
    
    \item 
    \begin{enumerate}
        \item $\lambda_2 = 2^\mu$ is regular, $\chi < \lambda_2 \Rightarrow \chi^\sigma < \lambda_2$, and $\lambda_1 = \min\{ \partial : 2^\partial > 2^\mu \}$ is both regular and strictly less than $2^\mu$
        
        \item like (B)(b) for $\lambda = \lambda_2$ but $|C_\delta| < \lambda_1$ for 
        $\delta \in S$ (so $C_\delta = \delta$ would not work)
        
        \item $\BB(J_S^{\mathrm{nst}}, \mu^+, \theta, \kappa)$ for every 
            $\theta < \mu$ and any stationary subset $S$ of $\lambda_1$
    
        \item[(c)$'$] like (B)(b), but for $\lambda = \lambda_1$, $S$ a club or simply not in the weak diamond ideal \cite{Sh:65}.
    \end{enumerate}
\end{enumerate}

\begin{note}\label{14.56}
In \cite{Sh:1028}
\begin{enumerate}[1)]
    \item We get scales with some two-cardinal freeness properties. This is used to get somewhat free $n$-dimensional scales.
    
    \item Hence for every $n < \omega$ there is an $\aleph_{\omega_1n+1}$-free (but not $\aleph_{\omega_1n+2}$-free) abelian group $G$ such that 
    $\Hom(G,\bbZ) = \{0\}$.
    
    \item Also, a complementary consistency result.
\end{enumerate}
\end{note}

\begin{note}
\label{14.61}
\underline{On tiny models}:  On tiny models see 
Laskovski, Pillay and Rothmaler \cite{LaPiRo}, 
$M$ is tiny if $\mu = \|M\| < |T|,T$ categorical 
in $|T|^+$, where $|T|$ is the number of formulas up to equivalence.  Assume
further that for $T$ not every regular type is trivial, \underline{then}  existence of
such $T$ for given $\mu$ is equivalent to the existence of 
$A_i \in [\mu]^\mu$ for $i < \mu^+$ such that 
$\bigwedge\limits_{i < j} |A_i \cap A_j| < \aleph_0$, hence 
necessarily $\mu < \beth_\omega$.  (Proved in the appendix of 
\cite{Sh:460}).
\end{note}

\begin{note}
\label{14.64}
\underline{On cofinalities of the symmetric group}:  
Let $\mathrm{Sp}$ be the family of regular $\lambda$ such that the permutation 
group of $\omega$ is the union of a strictly increasing chain of 
subgroups.  Now $\mathrm{Sp}$ has closure properties under $\pcf$, say if 
$n < \omega \Rightarrow \lambda_n \in \mathrm{Sp}$ then 
$\pcf\{\lambda_n : n < \omega\} \subseteq \mathrm{Sp}$ (Shelah and Thomas \cite{Sh:524}).
\end{note}

\begin{note}
\label{14.67}
\underline{Hanf number}

On application to Hanf numbers see Grossberg Shelah \cite{Sh:238}.
\end{note}

\begin{note}
\label{14.70}
On the number of non-isomorphic models: see \cite[\S2]{Sh:600}.  See
more in \cite{Sh:1029}, e.g. on groups.  A survey on existence of
universal, see \cite{Dj05}, more recently \cite{Sh:F1808}.
\end{note}

\begin{note}\label{14.73}
In addition, if $\mu = \aleph_\mu > \cf(\mu)$ then there is no universal 
linear order of cardinality $\lambda = \mu^+$. See \cite{Sh:F2150}.
\end{note}

\begin{note}\label{14.76}
We can sum: (see \cite[6.1]{Sh:F2150}) 

there is $T$ such that if $\lambda \notin \mathrm{Univ}(T)$, 
\underline{then} for some cardinal $\mu$
\begin{enumerate}
    \item[$(*)$] $\lambda = \mu^+$ and (a) or (b), where:
    \begin{enumerate}
        \item $\mu$ is singular strong limit such that $\mu = \aleph_\delta$, 
        $\delta < \mu$ and $2^\mu > \lambda$,
        
        \item $\mu$ is regular, $2^{<\mu} \leq \lambda < 2^\mu$ and 
        $\gb_\mu = \lambda_1 < \gd_\mu$.
    \end{enumerate}
\end{enumerate}
(see \cite[6.1]{Sh:F2150})
\end{note}

\newpage

\section {Discussion} \label{A17}
\bigskip

\begin{thesis}
\label{15.1} 
Artificially/naturality thesis.

Probably you will agree that, for a polyhedron, $v$ (number of vertices) $e$
(number of edges) and $f$ (number of faces) are natural measures, whereas
$e+v+f$ is not, but from a deeper point of view $v-e+f$ runs deeper than
all.  In this vein we claim: for $\lambda$ regular, $2^\lambda$ is the right
measure of ${\clP}(\lambda)$, and $\lambda^\kappa$ is a good measure of
$[\lambda]^{\le \kappa}$.  However, the various cofinalities are
better measures.  $\lambda^\kappa$ is an artificial combination of more basic
things of two kinds: the function $\lambda \mapsto 2^\lambda$ ($\lambda$
regular which is easily manipulated) and the various cofinalities we discuss
(which are not).  For example $\pp(\aleph_\omega) < \aleph_{\omega_4}$ is the
right theorem, not $\aleph^{\aleph_0}_\omega < \aleph_{\omega_4} +
(2^{\aleph_0})^+$ (not to say: $2^{\aleph_\omega} < \aleph_{\omega_4}$ when
$\aleph_\omega$ is strong limit).  Also the equivalence of the different
definitions which give apparently weak and strong measures, show naturality:
\mn
\begin{enumerate}
\item[$(a)$]  $\cf([\aleph_\omega]^{\aleph_0}) = \pp(\aleph_\omega)$
\sn
\item[$(b)$]   $\min\{|S|:S \subseteq [\lambda]^{\le \kappa}
\text{ stationary}\} = \cf([\lambda]^{\le \kappa},\subseteq)$
for $\kappa < \lambda$
\sn
\item[$(c)$]  if $\lambda \ge \mu > \theta > \sigma = \cf(\sigma)
> \aleph_0$ \underline{then} 
\[
\cov(\lambda,\mu,\theta,\sigma) = \sup\{\pp_{\Gamma(\cf(\chi),
\theta,\sigma)}(\chi):\mu \le \chi \le \lambda,\sigma \le \cf(\chi)
< \theta\}.
\]
\end{enumerate}
\mn
Note, $\chi \le \pp_\theta(\lambda)$ says $[\lambda]^{\le \kappa}$
is at least as large as $\chi$ in a strong sense, whereas $\chi \ge
\min\{|S|:S \subseteq [\lambda]^{\le \kappa}$ stationary$\}$
says that $[\lambda]^{\le \kappa}$ can be exhausted very well by
$\chi$ ``points" (for the right filters: measure 1).

We tend to think the $\pp$'s are enough, but there is a gap 
is our understanding concerning cofinality $\aleph_0$, mainly: 
is it true that
\mn
\begin{enumerate}
\item[$(*)$]   $\lambda > \cf(\lambda) = \aleph_0 \Rightarrow 
\cf([\lambda]^{\le \kappa},\subseteq) = \pp_{\aleph_0}(\lambda)$.
\end{enumerate}
\mn
We have many approximations saying that this holds in many cases
(see \ref{7.5}).

More generally, we should replace power by products, and cardinality by
cofinality, and therefore deal with $\pcf({\ga})$.
\end{thesis}

\begin{note}
\label{15.2}
\underline{The Cardinal Arithmetic below the continuum thesis}:

We should better investigate our various cofinalities without assuming
anything on powers (for example, the difference between the old result
$\pp(\aleph_\omega) < \aleph_{(2^{\aleph_0})^+}$ and the latter result
$\pp(\aleph_\omega) < \aleph_{\omega_4}$ is substantial); as
\mn
\begin{enumerate}
\item[$(a)$]   you should try to get the most general result (when it has
substance of course)
\sn
\item[$(b)$]   if we add many Cohen reals, all non-trivial products are
$\ge 2^{\aleph_0}$, but our various cofinalities do not change, so we should
not ignore this phenomenon
\sn
\item[$(c)$]   even if we want to bound $2^\lambda$ for $\lambda$ strong
limit singular, we need to investigate what occurs in the interval
$[\lambda,2^\lambda]$ which is a problem of the form indicated above; this
is central concerning the problem (see \cite{Sh:430}): if $\lambda$ is the
$\omega_1$-fixed point then $2^\lambda$ is $<$ the $\omega_4$-th fixed
point
\sn
\item[$(d)$]   looking at cardinal arithmetic without assumptions on the
function $\lambda \mapsto 2^\lambda$, makes induction on cardinality
more useful.
\end{enumerate}
\end{note}

\begin{thesis}
\label{15.3}
\mn
\begin{enumerate}
\item[$(A)$]  $\pp(\lambda)$ is the right power set operation.
\end{enumerate}
\mn
$\lambda \mapsto 2^\lambda$ ($\lambda$ regular) is very elastic, you can
easily manipulate it, but $\pp(\lambda)$ ($\lambda$ singular) and
$\cov(\lambda,\mu,\theta,\sigma)$ are not; it is hard to manipulate them, and
we can prove theorems about them in ZFC.
\mn
\begin{enumerate}
\item[$(B)$] pcf, cov are basic operations, with non-trivial ZFC results. 
\end{enumerate}
\mn
1) Consider $[\lambda]^{\le \kappa}$, the
family of subsets of $\lambda$ of cardinality $\le \kappa$, when $\lambda >
\kappa$ (see \ref{15.1}). 

\noindent
2)  $\lambda^\kappa$ is the crude measure of $[\lambda]^{\le \kappa}$.

It is very interesting to measure it, and cardinality is generally a very
crude measure; $\pp_\kappa(\lambda)$ is a fine measure; and we have
intermediate ones: $\cf([\lambda]^{\le \kappa},\subseteq)$,
min$\{|S|:S \subseteq [\lambda]^{\le \kappa} \text{ stationary}\}$ and
more. 
The best is when we can compute cruder 
numbers from finer ones; particularly when they are equal, so we could use
different definitions for the same cardinal depending on what we want to
prove.  So we want to show that the $\pp_{\Gamma(\cf(\lambda))}(\lambda)$
for $\lambda$ singular is enough.
\end{thesis}

\begin{note}
\label{15.4}
$\pp_{\Gamma(\theta,\sigma)}(\lambda)$ is the finest we
have for what we want; they are like the skeleton of set theory; you can 
easily change your dress and even can manage to change how much flesh you
have; but changing your bones is harder.  You may take hypermeasurable
$\lambda$, blow up $2^\lambda$ and make it singular; this does not affect
for example $\pp_{\Gamma(\aleph_1)}(\lambda^*)$ when $\lambda^* > \lambda,
\cf(\lambda^*) = \aleph_1$ (even if $\lambda^* <$ new $2^\lambda$), nor
$\cov(\lambda^*,\lambda^*,\aleph_1,\sigma)(\sigma = 2,\aleph_1)$; 
they measure really how many subsets of $\lambda^*$ of cardinality 
$\aleph_1$ there are - not through some $\lambda' < \lambda^*$ 
having many subsets of cardinality $\le \aleph_1$.
\end{note}

\begin{note}
\label{15.5}
Subconscious remnants of GCH have continued to influence
the research: concentration on strong limit cardinals; but from our point of
view, even if $2^{\aleph_0}$ is large and $\mu < 2^{\aleph_0} 
\Rightarrow 2^\mu = 2^{\aleph_0}$, the cardinal arithmetic below 
$2^{\aleph_0}$ does not become simpler.

Also GCH was used as an additional assumption (or semi-axiom), but rarely
was the negation of CH used like this: 
simply because one didn't know to prove interesting theorems from
$\neg \CH$.  But now we know that violations of $\GCH$
have interesting consequences (see below).
\end{note}

\begin{note}
\label{15.6}
Up to now we have many consequences of $\GCH$ (or instances of
it) and few of the negations of such statements.  We now begin to have
consequences of the negation, for example see here \ref{11.9}; so we can
hope to have proofs by division to cases.  For example, let $\lambda$ be a
strong limit singular; if $\pp(\lambda) > \lambda^+$ then $\NPT(\lambda^+,
\cf(\lambda))$ and if $\pp(\lambda) \le \lambda^+$ then $2^\lambda =
\lambda^+$ (and $\diamondsuit^*_{\{\delta < \lambda^+:\cf(\delta) \ne
\cf(\lambda)\}}$) and so various constructions are possible (see
here \ref{11.9}(b) and \cite{Sh:462} on more, also \cite{Sh:E9},
\cite{Sh:534}).
\end{note}

\begin{note}
\label{15.7}
The right problems.

An outside viewer may say that the main problem,

\[
(\aleph_\omega = \beth_\omega \Rightarrow 2^{\aleph_\omega} < 
\aleph_{\omega_1})
\]

\mn
was not solved.  As an argument we may accuse others: maybe 
$\aleph_{\omega_4}$ is the right bound.  But more to the point 
is our feeling that this is not the right problem. The right problems are:
\mn
\begin{enumerate}
\item[$(\alpha)$]   Does $\pcf({\ga})$ always have cardinality
$\le |{\ga}|$? 
\sn
\item[$(\beta)$]   Is $\cov(\lambda,\lambda,\aleph_1,2) =^+
\pp(\lambda)$ when $\cf(\lambda) = \aleph_0$?
\end{enumerate}
\mn
Now $(\alpha)$ is just a member of a family of problems quite 
linearly ordered by implication discussed in \cite[\S6]{Sh:420}, 
\cite{Sh:460}, which seem unattackable both by the forcing methods 
and ZFC methods.  The borderline between chaos and order seems
\mn
\begin{enumerate}
\item[$(\alpha)^-$]  Can $\pcf({\ga})$ have an accumulation point which
is an inaccessible cardinal? (Hopefully not.)
\end{enumerate}
\mn
Similarly $(\beta)$ is the remnant of the conjecture that all
cov$(\lambda,\mu,\theta,\sigma)$ can be expressed by the values of
$\pp_{\Gamma(\theta,\sigma)}(\lambda')$ and even 
$\pp_{\Gamma(\cf(\lambda'))}(\lambda')$; this has been proved in many
cases (see \ref{7.5}).  On an advance see \cite{Sh:460}.

Also though $(\alpha), (\beta)$ have not been solved, much of what we want
to derive from them has been proved.

Another problem on which no light was shed is:
\mn
\begin{enumerate}
\item[$(\gamma)$]   if $\lambda$ is the first fixed point, find a bound on
$\pp(\lambda)$ (or better $\cov(\lambda,\lambda,\aleph_1,2)$).
\end{enumerate}
\mn
We can hope for the $\omega_4$-th fixed point, to serve as a bound but will
be glad to have the first inaccessible as a bound.  Even getting a bound
assuming GCH below $\lambda$ would open our eyes.  This becomes a problem
after \cite{Sh:111}, \cite[Ch.XII,\S5,\S6]{Sh:b}.
\mn
\begin{enumerate}
\item[$(\delta)$]   Generalize \cite[\S1]{Sh:355} to deal with what occurs
above $\tlim_I \, \lambda_i$ 
(for example \ref{4.1}, $(\lambda,\sigma)$-entangled linear order).
\end{enumerate}
\mn
More accurately, assume $\prod\limits_{i < \delta(*)} \lambda_i/J$ has true
cofinality $\lambda,\mu = \tlim_I(\lambda_i) = \sup(\lambda_i),
\lambda_i$ regular $> \delta(*)$, and $\sup_{i < \delta(*)} 
\lambda_i < \theta = \cf(\theta) < \lambda$.  We can 
find regular $\lambda'_i < \lambda_i$ such that 
$\tcf(\prod \, \lambda'_i/J) = \theta$) as exemplified by $\bar f$, 
which is $\mu^+$-free (hence tlim$(\lambda'_i) = \lambda_i$) in addition: if
$\delta < \theta,\cf(\delta) < \theta$ and $\cf(\delta) > 2^{|\delta(*)|}$
(or just $\bar f \restriction \delta$ has a $<_J$-lub) then \wilog \,
$f_\delta/J$ is the $<_J$-lub of $\bar f \rest \delta$, we want to
know something on $\langle \cf(f_\delta(\alpha)):\alpha < \delta
\rangle$.  For more information see \cite[4.1,4.1A]{Sh:400}. 

Note that we also do not know, for example
\mn
\begin{enumerate}
\item[$(\eps)$]   if $\cf(\lambda) \le \kappa < \lambda$, is
$\cf(\pp_\kappa(\lambda)) > \lambda$? (we know that it is $> \kappa$)
\sn
\item[$(\zeta)$]   we believe pcf considerations will eventually have
impact on cardinal invariants of the continuum, but this has not materialized
so far.
\end{enumerate}
\end{note}

\begin{note}
\label{15.8}
The perspective here led to phrasing some hypotheses, akin to $\GCH$ or $\mathrm{SCH}$.

The ``strong hypothesis" says $\pp(\lambda) = \lambda^+$ for (every) singular
$\lambda$. Note it is like $\GCH$ but is not affected by, say, c.c.c. forcing;
it follows from $\neg 0^\#$ and from $\GCH$; its negation is known to be
consistent and I feel it is a natural axiom.  Other hypotheses may still follow from 
$\ZFC$: for example, the medium hypothesis says $|\pcf({\ga})| \le |{\ga}|$, and the weak says 
$\{\mu:\pp(\mu) \ge \lambda,\ \mu < \lambda,\ \cf(\mu) = \aleph_0[> \aleph_0]\}$ 
is countable [finite]. There are intermediate ones;
such hypotheses and consequences are dealt with in \cite[\S6]{Sh:420}, see
more in \cite{Sh:460}, \cite{Sh:513}.  Particularly concerning the connection
of the medium and weak ones, (see \ref{12.3}, \ref{7.18}),
$\mathrm{ZF} + \mathrm{DC} + [\alpha]^{\aleph_0}$ well ordered suffice, see \cite{Sh:835}.
\end{note}

\begin{note}\label{15.11}
In a major advance, Gitik has proven in \cite{Git20} that the following 
version of the weak hypothesis fails: $\{ \mu : \cf(\mu) = \kappa < \mu,\ \pp_{\kappa\text{-complete}}(\mu) \geq \lambda \}$ may be large.

Still open is whether, in the RGCH (see \cite{Sh:460}), we can replace $\beth_\omega$ by $\aleph_\omega$.
\end{note}

\bigskip

\centerline {$* \qquad * \qquad *$}
\newpage

\section {Part B - Corrections to the book \cite{Sh:g}} \label{B}
\bigskip

\noindent
\underline{page 50,line 22}:  see more in Part C.
\bigskip

\noindent
\underline{page 51,line 12}:  replace $\lambda^+$ by $\mu$.
\bigskip

\noindent
\underline{page 51,line 13}:  see more in Part C.
\bigskip

\noindent
\underline{page 66,Theorem 3.6}:  second line of theorem: 

replace $\lambda^{\beta +1}$ by $\lambda^{\beta +1}_0$
\smallskip

\noindent
\underline{add} after the second line of Remark 3.6A: 

\noindent
2)  This is essentially the proof from \cite[Ch.XIII,\S6]{Sh:b} and more
appears in Ch.IX 
\medskip

\noindent
first line of the proof: 

\underline{replace} $\lambda > \aleph_0$ by ``$\lambda_0 > |\alpha|^+$ 
(why? as we can replace $\lambda_0$ by $\lambda^+_0$ and deduce 
the result on the original $\lambda_0$ from the result on $\lambda^+_0$)"
\medskip

\noindent
\underline{replace} fifth line of the proof: 

$N_h = \bigcap\big\{\text{Skolem Hull}_M\big(\lambda_0 \cup 
\bigcup\limits_{\beta < \alpha} C_\beta\big) : C_\beta 
\text{ a club of } f(\lambda^{+\beta+1}) \text{ for } \beta < \alpha \big\}$
\medskip

\noindent
\underline{add} in the end of the proof: 

Clearly this family is a family of subsets of $\lambda$ each of cardinality
at most $\lambda_0$ of the right cardinality.  So we have to prove just that
it is cofinal.  So let $X$ be a subset of $\lambda$ of cardinality at most
$\lambda_0$, and we shall find a member of the family which includes it.
Let $\chi$ be large enough.  By \ref{3.4} we can find an elementary
submodel $N_i$ of $(\clH(\chi),\in,<^*_\chi)$, for $i \le \delta =:
|\alpha|^+$ each of cardinality such that $\{F,\lambda_0,\alpha^*,\lambda,X,
f,g\} \in N_i$ and $i < j \rightarrow N_i \in N_j$ increasing continuous with
$i$ and condition (b) form 3.4 holds for $f \in F$.

It is enough to prove that
\mn
\begin{enumerate}
\item[$(*)$]   $N_f$ includes $N_\delta \cap \lambda$
\end{enumerate}
\mn
for this it is enough to prove
\mn
\begin{enumerate}
\item[$(**)$]   if $C_\beta$ is a club of $\lambda^{\beta +1}_0$ for each
$\beta < \alpha^*$ and $M'$ is the Skolem Hull in $M$ of $\lambda_0 \cup
\bigcup \{C_\beta:\beta < \alpha^*\}$ then $M'$ include $N_\delta \cap 
\lambda$.
\end{enumerate}
\mn
For this we prove by induction on $\gamma \le \alpha$ that
\mn
\begin{enumerate}
\item[$(**)_\gamma$]  $M'$ includes $\lambda \cap \lambda^{+ \gamma}_0$.
\end{enumerate}
\medskip

\noindent
\underline{Case 1}:  $\gamma = 0$.

In this case as $M$ includes $\lambda_0$ this is trivial.
\medskip

\noindent
\underline{Case 2}:  $\gamma$ a limit cardinal ordinal.

In this case the induction hypothesis implies the conclusion trivially.
\medskip

\noindent
\underline{Case 3}:  $\gamma = \beta +1$.

Use the induction hypothesis and the choice of the functions $f$ and $g$. 
(See more Ch.IX, \ref{3.3})
\bigskip

\noindent
\underline{page 136,lines 21,22,23}:\label{pg136}

replace by: 

No problem to define.  We define $B^\alpha_i$ (for $i < \lambda,\alpha \in S$) 
by induction on $\alpha$:

\[
B^\alpha_i = \begin{cases} \big\{\beta:\cf(\beta) \ne \lambda \text{ and }
\beta \in A^\alpha_i \vee \beta = \sup(\beta \cap A^\alpha_i)\big\} 
&\text{ \If \, } \cf(\alpha) \ne \aleph_1 \\
  \bigcap \big\{\bigcup\limits_{\beta \in C} B^\beta_i:
C \text{ a club of } \alpha \text{ such that } 
\bigwedge\limits_{\beta \in C} \cf(\beta) = \aleph_0 \big\} &\text{ \If \, } 
\cf(\alpha) = \aleph_1
\end{cases}
\]

\mn
(or see \cite[4.1]{Sh:351}).
\bigskip

\noindent
\underline{page 210, line 15}:

add: or $\lambda$ is not Mahlo and we can use Ch.III.
\bigskip

\noindent
\underline{page 222,line 24}: replace by:

Definition 1.4.  1) We say $D$ is strongly nice if it is strongly nice to
every 
\bigskip

\noindent
\underline{page 224,line 8}: replace by:

\[
\sup\Big\{\prod\limits_{i < \omega_1} f(i)/D : D 
\text{ is a normal filter extending } D^*\Big\}
\]
\bigskip

\noindent
\underline{page 228,line 1}: replace $D^*$ by $D^* \in V^*$.
\bigskip

\noindent
\underline{pages 334-337}:  see a rewriting in \cite{Sh:E11}
\bigskip

\noindent
\underline{page 334,line -4} replace by: 

(2)  The first phrase follows from part 1 and check the second
\bigskip

\noindent
\underline{page 335,line 4}: replace ``$f \restriction {\gb}
\mu[{\ga}] \le f^\mu_\alpha$" by ``$f \restriction {\gb}_\mu[{\ga}] 
\le f^\mu_\alpha$" 
\bigskip

\noindent
\underline{page 335,line 18}:  space after $\varnothing$; replace 
$\bigcap\limits^{n}_{\ell =1} {\gb}_{\sigma_\ell}
[{\ga}]$ by $\bigcup\limits^{n}_{\ell =1} {\gb}_{\sigma_\ell}[{\ga}]$
\bigskip

\noindent
\underline{page 336,line 3}: replace ${\gb}$ by ${\gc}$
\bigskip

\noindent
\underline{page 336,line -7}: replace $\square_{3.3}$ by $\square_{3.2}$
\bigskip

\noindent
\underline{page 381,lemma3.5 and page 383,line 21}: 
No!  But see \cite[5.12]{Sh:400} and \cite[\S6]{Sh:513}
\bigskip

\noindent
\underline{page 410,line -1}: 
replace by: $\{\delta < \sigma:\cov(\lambda_\delta,\lambda_\delta,
\theta^+,2) < \mu_\delta\}$ contains a club of $\sigma$, where
\mn
\begin{enumerate}
\item[$(*)(i)$]   let $\mu_\delta$ be $\pp^{\cer}_\theta(\lambda_\delta)$ the
first regular $\mu > \lambda_\delta$ such that: 

if ${\ga} \subseteq \Reg \cap \lambda_\delta \setminus |{\ga}|^+$,
then $\sup\{\max \pcf({\gb}):{\gb} \subseteq {\ga},
|{\gb}| \le \theta \text{ and } (\forall \chi < \lambda_\delta)
\max \pcf({\mathfrak b} \cap \chi) < \lambda_\delta\}$ 
(so normally this means $\cov(\lambda_\delta,\lambda_\delta,\theta^+,2) =^+
\pp_\theta(\lambda_\delta)$).
\end{enumerate}
\bigskip

\noindent
\underline{page 411,line 1}: replace by: 
\mn
\begin{enumerate}
\item[$(ii)$]   $\cov(\lambda,\lambda,\theta^+,2) < \pp^{\cer}_\theta
(\lambda)$ which normally means $\cov(\lambda,\lambda,\theta^+,2) =^+
\pp_\theta(\lambda)$, e.g. if $\cov(\lambda_i,\theta^+,\theta^+,2) 
< \lambda$ for a club of $i < \sigma$
\sn
\item[$(iii)$]   if e.g. $\sigma^{\aleph_0} < \lambda$, then we can add 
$\{\delta < \sigma$: if $\cf(\delta) = \aleph_0 \text{ then } 
\pp^{\cer}_{J^{\bd}_\omega}(\lambda_\delta) > 
\cov(\lambda_\delta,\lambda_\delta,\theta,2)\}$ contains a club 
(for the changes needed for the proof see below, Part C).
\end{enumerate}
\bigskip

\noindent
\underline{page 417,line 11}: add:

Here examples are constructed for $\lambda$ singular and in \cite{Sh:572}
for $\lambda = \aleph_1$ which was the last case.
\bigskip

\noindent
\underline{page 418,line 20}: sequence of \underline{not} sequence of ...
\newpage

\section{Part C - Expansions for \cite{Sh:g}} \label{C} 
\bigskip

\noindent
\S17 \quad Short Expansions
\label{C1}
\bigskip

\noindent
\underline{page 50,line 22}: add: [this is the proof of II,1.4(3)].
\bigskip

\noindent
\underline{Case 1}:  otp$(A)$ is zero.

Trivial.
\bigskip

\noindent
\underline{Case 2}:  otp$(A)$ is a successor ordinal.

Let $\alpha$ be the last member of $A$ and let $A'$ by $A \setminus
\{\alpha\}$.  Clearly the order type of $A'$ is (strictly smaller than that
of $A$) hence by the induction hypothesis we can find $s'_\beta \in i$ for
$\beta \in A'$ as required.  Define $s_\beta$ for $\beta \in A$ as follows:

if $\beta = \alpha$, then $s_\beta = \varnothing$ and if $\beta \in A'$ then
$s_\beta =: \{i < \kappa:i \in s'_\beta$ or $f_\alpha(i) \le
f_\beta(i)\}$.  Now $s_\beta$ is a subset of $\kappa$ and if $\beta = 
\alpha$ is the union of two sets: $s'_\beta$ and $\{i <
\kappa:f_\alpha(i) \le f_\beta(i)\}$, now the first belongs 
to $I$ by its choice and the second as
we know $f_\beta <_I f_\alpha$ (because $\beta < \alpha$).  So $S_\beta$,
their union is in $I$, too. 

This holds also in the case $\beta = \alpha$.  So $s_\beta \in I$ for $\beta
\in A$, and it is easy to check the requirements.
\bigskip

\noindent
\underline{Case 3}:  otp$(A)$ is a limit ordinal.

Let $\delta$ be $\sup(A)$, so is a limit ordinal.  So by 
\ref{1.3}(ii)$(\delta)$ there is a closed unbounded subset $C$ of 
$\delta$ and sets $\tau_\alpha \in I$ for $\alpha \in C$ such that 
$i \in \kappa \setminus \tau_\alpha \setminus s_\beta$ and $\alpha 
< \beta$ implies $f(i) < f_\beta(i)$. 

Without loss of generality $0 \in C$ (let $t_0 =: \{i < \kappa:f_0(i) \ge
f_{\min(A)}(i)\}$).

Now for every $\alpha \in C$ let $A_\alpha =: A \cap [\alpha,\min(A
\setminus (\alpha +1))$.  Clearly $\otp(A_\alpha) < \otp(A)$, let
$A'_\alpha =: A_\alpha \cup \{\alpha\}$.  So $\otp(A'_\alpha) = 1 +
\otp(A_\alpha) < \otp(A)$ (as the latter is a limit ordinal).
So we can apply the induction hypothesis, getting $s'_\beta$ for $\beta \in
A'_\alpha$ as guaranteed there.

Now we define $s_\beta$ for $\beta \in A$ as follows: let $\alpha_\beta =:
\sup(C \cap \beta)$ and 
$\gamma_\beta =: \min(A \setminus (\alpha + 1))$.  
So $\beta \in A_{\alpha_\beta}$, hence $s'_\beta$ is well defined, and let 
$$s_\beta =: s'_\beta \cap \{i < \kappa: \text{ it is not true that }
f_{\alpha_\beta}(i) \le f_\beta(i)\}.$$

Now check.

\bigskip
\centerline {$* \qquad * \qquad *$}
\bigskip

\noindent
\underline{page 51, line 13}: add to the end of line 
(this is line 7 of the proof of II,1.5A).

Of course, we do not have knowledge on the relation between $f_\alpha(i)$
and $f_\beta(j)$, so we just e.g. use $f'_\alpha$ defined by
$f'_\alpha(i) =: \kappa f_\alpha(i) + i$ (so $f'_\alpha$ is a function from
$\kappa$ to $\lambda$, as $\kappa < \lambda$).  Now 
$\langle f'_\alpha : \alpha < \mu \rangle$ is as required (note that 
$\big\langle \{f_\alpha(i) : i < \mu\} : i < \kappa \big\rangle$ is a sequence 
of pairwise disjoint subsets of $\lambda$).


\newpage

\section {More on II,3.5} \label{C2}

This refinement is used in \cite{Sh:810}.
\begin{claim}
\label{k.1}
Assume
\mn
\begin{enumerate}
    \item[$(a)$]   ${\ga} = \{\lambda_i:i < \delta\}$ is an
    increasing sequence of regular cardinals $> \delta$
\sn
    \item[$(b)$]   $\lambda = \tcf_\pi({\ga},<_{J^{\bd}_\delta})$
\sn
    \item[$(c)$]  $\lambda_0 > 2^{|i|}$ for $i < \delta$ or just
    $i < \delta \Rightarrow \lambda_0 > |\pcf({\ga})\lambda_i|$
\sn
    \item[$(d)$]   $\cf(\delta) > \aleph_0$
\sn
    \item[$(e)$]   $S =: \Big\{i < \delta:\text{for some } i_0 <
    i_1,\ \pcf\{\lambda_j:i_0 < j < i\} \setminus \sum\limits_{j<i}
    \lambda_j$ is a singleton cardinal $< \underset{j < \delta} {}\to \sup
    \lambda_j\Big\}$ is stationary.
\end{enumerate}
\mn
\underline{then}  we can find $\langle f_\alpha:\alpha < \lambda \rangle$ such that
\mn
\begin{enumerate}
    \item[$(a)$]   $f_\alpha \in \prod\limits_{i < \red{\delta}} \lambda_i$ is
    $<_{J^{\bd}_\delta}$-increasing and cofinal
\sn
    \item[$(b)$]   if $f \in \prod\limits_{i < \delta} \lambda_i$ and
    $(\forall i < \delta)(\exists \alpha < \lambda)[f \rest i =
    f_\alpha \rest i]$ then $f \in \{f_\alpha:\alpha < \lambda\}$.
\end{enumerate}
\end{claim}

\begin{remark}
This is just the proof of \cite[Ch.II,3.5]{Sh:g}, just we use more of it.
\end{remark}

\begin{PROOF}{\ref{k.1}}
Let $\mu = \sum\limits_{i < \delta} \lambda_i$ and
$\mu_j = \sum\limits_{i < j} \lambda_i$ for $j < \delta$.  

Recall ${\ga} = \{\lambda_i:i < \delta\}$, so $\min({\ga}) > 
|\mu \cap \pcf({\ga})|$.  Let $\bar{\gb} = \langle {\gb}_\theta:
\theta \in \pcf({\ga}) \rangle$ be a generating sequence for 
$\pcf({\ga})$.  Choose $\langle \bar f^\theta:\theta \in \pcf({\ga}) 
\rangle$ as in claim \ref{k.3} below.  Now we let 

$$\begin{aligned}
{\cF} = 
\Big\{f \in \prod\limits_{i < \delta} \lambda_i : &\text{ for every } 
\theta \in \pcf({\ga}) \text{, for some } n < \omega\\
&\text{and } \theta_0 < \ldots < \theta_{n-1} \text{ from } \pcf({\gb}_\theta)\\ &\text{and } \alpha_0 < \theta_0,\dotsc,\alpha_{n-1} < \theta_{n-1}\\ 
&\text{we have } f \rest {\gb}_\theta = \max\{f^{\theta_\ell}_{\alpha_\ell}:\ell < n\}\Big\}.
\end{aligned}$$

Let $f_\alpha := f^\lambda_\alpha$ for $\alpha < \lambda$.

First clearly
\mn
\begin{enumerate}
\item[$(*)_1$]   $\alpha < \lambda \Rightarrow f_\alpha \in {\cF}$.
\end{enumerate}
\mn
Secondly, the main point is
\mn
\begin{enumerate}
\item[$(*)_2$]   if $f',f'' \in {\cF}$ then $f'
  <_{J^{\bd}_{\ga^\delta}} f''$ or $f' =_{J^{\bd}_\delta} f''$ or $f''
<_{J^{\bd}_{\ga}} f'$.
\end{enumerate}
\mn
Why does $(*)_2$ hold? Given $f',f'' \in {\cF}$, let 
$${\gc}_1 = \{\theta \in {\ga} : f'(\theta) < f''(\theta)\},$$ 
$${\gc}_2 = \{\theta \in {\ga} : f'(\theta) = f''(\theta)\},$$ 
$${\gc}_3 = \{\theta \in {\ga} : f'(\theta) > f''(\theta)\},$$ 
so 
$\langle {\gc}_1,{\gc}_2,{\gc}_3 \rangle$ is a partition of ${\ga}$.
Let $E = \{i < \delta:\text{ for } \ell = 1,2,3$ if $\sup({\gc}_\ell) 
= \sup({\ga})$ then $\sup({\gc}_\ell \cap \lambda_i) =
\sup({\ga} \cap \lambda_i)$ and if $\sup({\gc}_\ell) < 
\sup({\ga})$ then $\sup({\gc}_\ell) < \lambda_j$ for some $j < i\}$.  

Clearly $E$ is a club of $\delta$; by clause (c) of the assumption, $S
\subseteq \delta$ is stationary hence $S \cap E \ne \varnothing$, so let
$i \in S \cap E$ and let $\theta_i$ be the single member of
$\pcf({\ga} \cap \lambda_i) \setminus \mu_i = \pcf(\{\lambda_j:
j < i\}) \setminus \mu_i$ (recall the definition of $S$).  
So ${\gb}_{\theta_i}$ contains an end-segment of ${\ga} \cap \lambda_i$ --- 
say ${\gb}'$.  By the choice of ${\cF}$ and the assumption
$f',f'' \in {\cF}$ and the choice of $\theta_i$, we know that for 
some end segment ${\gb}''$ of ${\gb}',f' \rest {\gb}'' 
\in \{f^{\theta_i}_\alpha \rest {\gb}'':\alpha < \theta_i\}$ 
and \wilog \, also $f'' \rest {\gb}'' \in 
\{f^{\theta_i}_\alpha \rest {\gb}'':\alpha < \theta_i\}$.  
So for some $\beta'$, $\beta'' < \theta_i$ we have 
$f' \rest {\gb}'' = f^{\theta_i}_{\beta'} \rest {\gb}''$ 
and $f'' \rest {\gb}'' = f^{\theta_i}_{\beta^{''}} \rest {\gb}''$.

Now $\beta' < \beta'' \vee \beta' = \beta'' \vee \beta' > \beta''$ and
accordingly we get one of the three possibilities in $(*)_2$.

Now clearly we are done.  
\end{PROOF}

\begin{claim}
\label{k.2}
1) In \ref{k.1} we can weaken assumption (e) to
\mn
\begin{enumerate}
\item[$(e)^-_{\mathfrak a}$]    letting $\langle \mu_i:i < \sigma \rangle$ be
increasing continuous with limit $\sup({\ga})$ so $\sigma = 
\cf(\sup({\ga}))$ for some normal filter $D$ on $\cf(\sup({\ga})) 
= \cf(\delta)$ we have:
\sn
\item[$(e)^-_D$]   if ${\ga}' \subseteq {\ga}(=:\{\lambda_i:
i < \delta \})$ and $\sup({\ga}') = \sup({\ga})$ then 
$$\{i < \sigma:\max(\pcf({\ga}' \cap \mu_i)) = \max(\pcf({\ga} 
\cap \mu_i))\} \in D.$$
\end{enumerate}
\mn  
2) Assume ${\ga}$ has no last element, $\cf(\sup({\ga})) > \aleph_0$, 
and $\lambda = \tcf(\pi{\ga}/J^{\bd}_{\ga})$ and 
$\mu < \sup({\ga}) \Rightarrow \max \pcf(\ga \cap \mu) < \sup(\ga)$,
(e.g. ${\ga} = \{\lambda_i:i < \delta\}$ from \ref{k.1} assuming
clauses (a)-(d) of \ref{k.1}.

\underline{then}  for some unbounded ${\ga}^* \subseteq {\ga}$, we
have (clause (a), (b), (c) of \ref{k.1} and) clause $(e)^-_{\ga^*}$ 
of part (1) holds (hence the conclusion of \ref{k.1}).
\end{claim}

\begin{PROOF}{\ref{k.2}}
1) Let ${\ga} = \{\lambda_i:i < \delta\}$ such that 
$\lambda_i$ is regular increasing with $i$.

We repeat the proof of \ref{k.1}.  So our problem is that in proving
$(*)_2$, so we have $f',f''\in {\cF}$ and having defined 
the partition ${\gc}_1,{\gc}_2,{\gc}_3$ of ${\ga}$, at least 
two parts are unbounded in ${\ga}$ say ${\gc}_{\ell_1},{\gc}_{\ell_2}$
\mn
\begin{enumerate}
\item[$\boxtimes$]    if $\theta \in \pcf({\mathfrak c}_\ell)
\setminus \{\lambda\}$ then ${\gb}_\theta \setminus {\gc}_\ell 
\in J_\theta[{\ga}]$.
\end{enumerate}
\mn
[Why $\boxtimes$?  As in the proof of \ref{k.1}, we know that
$f',f'' \in {\cF}$ hence for some ${\gc} \in \mathbf J_{<
\theta}[{\ga}]$ we have $f' \restriction ({\gb}_\theta
\setminus {\gc}),f'' \restriction ({\gb}_\theta \setminus
{\gc})$ belongs to $\{f^\theta_\alpha \restriction ({\gb}_\theta 
\setminus \theta):\alpha <\theta\}$ and we continue as there.]
\mn
Now for $\ell=1,2,3$ we have $\sup({\gc}_\ell) = \sup({\ga}) 
\Rightarrow S_\ell := \{j < \sigma:\max \pcf({\gc}_\ell \cap \mu_\ell) 
= \max \pcf({\ga} \cap \mu_\ell)\} = \cf(\sup({\ga}))\} \in D$
also $E := \{i:\sup({\ga} \cap \mu_i) = \mu_i\}$ is a club of $\sigma$.
Hence $S = E \cap S_1 \cap S_2 \cap S_3 \in D$.

So for the $D$-majority of $j < \sigma$ we have sup$({\gc}_{\ell_1}
\cap \mu_j) = \mu_j = \sup({\gc}_{\ell_2} \cap \mu)$ and 
$\max \pcf({\gc}_{\ell_j} \cap \mu_j) = \max \pcf({\ga} \cap
\mu_j) = \max \pcf({\gc}_{\ell_2} \cap \mu_j)$ and we get
contradiction by $\boxtimes$. 

\noindent
2) We try to choose $\langle {\ga}_\eta:\eta \in {}^n \sigma
\rangle$ by induction on $< \omega$ such that
\mn
\begin{enumerate}
\item[$(i)$]   ${\ga}_{\LL\ \RR} = {\ga}$
\sn
\item[$(ii)$]   $\ga_\eta \subseteq \ga_{\eta \rest n}$ 
for $\eta \in {}^{n+1} \sigma$
\sn
\item[$(iii)$]  $\sup(\ga_\eta) = \sup(\ga)$
\sn
\item[$(iv)$]   for every $\eta \in {}^n \sigma$ for some club
$E_\eta$ of $\sigma$ we have: for every $j \in E_\eta$ there is $i < j$ 
such that $\max \pcf(\ga_{\eta \caret \LL i\RR} \cap \mu_j) < \max \pcf(\ga_n \cap \mu)$.
\end{enumerate}
\mn
Now for $n=0$ there is no problem and if ${\ga}_n,\langle {\ga}_\eta:
\eta \in {}^n \sigma \rangle$ has been chosen but there is no
suitable $\langle {\ga}_\eta:\eta \in {}^{n+1} \sigma \rangle$
then for some $\eta \in {}^n \sigma$ letting 
$${\cP}_\eta = \big\{\{i < \sigma : \max \pcf({\gb} \cap i) < \max \pcf({\ga}_\eta \cap i)\}:{\gb} \subseteq {\ga}_\eta,\sup({\gb}) = \sup({\ga}_\eta) \big\},$$ the normal 
ideal $D_\eta$ (on $\sigma$) which
${\cP}_\eta$ generates satisfies $\varnothing \notin D_\eta$ so 
${\ga}_\eta,D_\eta$ are as required.  Lastly, not all the
${\ga}_\eta$'s are defined as then we let 
$$E = \{i < \sigma : i \text{ a limit ordinal such that } \eta \in {}^{\omega >} i \Rightarrow i \in E_\eta\}$$ 
Clearly $E$ is a club of $\sigma$.  Now for any $i \in E$, we choose, by induction 
on $n < \omega$, a sequence $\eta_n \in {}^n i$ such that $\eta_n \triangleleft \eta_{n+1}$ 
and $\max \pcf({\ga}_{\eta_n}) > \max \pcf({\ga}_{\eta_{n+1}})$.
We let $\eta_0 = \LL\ \RR$, and $\eta_{n+1}$ will exist by clause (iv).  So
$\langle \max \pcf({\ga}_{\eta_n}):n < \omega \rangle$ is a
strictly decreasing sequence of cardinals, a contradiction.  So we are
done.
\end{PROOF}

\begin{claim}
\label{k.3}
Assume
\mn
\begin{enumerate}
\item[$(a)$]   $|\pcf({\ga})| < \,\min({\ga})$, ${\ga}$ as usual 
a set of regular cardinals
\sn
\item[$(b)$]    $\bar{\gb} = \langle {\gb}_\theta:\theta \in
[{\ga}] \rangle$ a generating sequence for $\pcf(\ga)$
\magenta{(exists by x.x - FILL)} which is closed (i.e. $\mu \in {\gb}_\theta
\Rightarrow {\gb}_\mu \subseteq b_\theta)$ and smooth (i.e. 
$\pcf({\gb}_\mu) \cap {\ga} = {\gb}_\mu$).
We can choose by induction on $\theta \in \pcf({\ga}),\bar f^\theta 
= \langle f^\theta_\alpha:\alpha < \lambda \rangle$ such that
\sn
\begin{enumerate}
\item[$(\alpha)$]   $f^\theta_\alpha \in \prod {\gb}_\theta$
is $<_{J_{< \theta}[{\gb}_\theta]}$-increasing and cofinal
\sn
\item[$(\beta)$]   if $\theta \in \pcf({\ga})$, $\alpha < \theta$ and 
$\mu \in {\gb}_\theta$, \underline{then} for some $n < \omega$, 
$\mu_0,\dotsc,\mu_{n-1} \in \pcf({\gb}_\mu)$,
$\beta_0 < \mu_0,\dotsc,\beta_{n-1} < \mu_{n-1}$, and $\beta < \mu$ we have 
$$f \rest {\gb}_\mu = \max[\{f^{\mu_\ell}_{\beta_\ell}:\ell < n\}].$$
\end{enumerate}
\end{enumerate}
\end{claim}

\begin{proof}
This is a restatement of \cite[Ch.VII,\S1]{Sh:g}.
\end{proof}

\begin{claim}
\label{k.5}
Assume $\kappa$ is regular and $\bar\theta = \langle \theta_i:
i < \kappa\rangle$ is a sequence of regular cardinals
$> \kappa^+$.  \underline{then}  for some $u,E,\bar \lambda,\lambda$ and $D$
we have
\mn
\begin{enumerate}
\item[$(a)$]   $u \subseteq \kappa$ is unbounded
\sn
\item[$(b)$]   $\lambda = \tcf(\prod\limits_{i \in u}
\theta_i,<_{J^{\bd}_u}\rangle$
\sn
\item[$(c)$]  $E := \{\delta < \kappa:\delta$ a limit ordinal and
$\delta = \sup(u \cap\delta)$
\sn
\item[$(d)$]  $\bar \lambda = \langle \lambda_\delta:\delta \in E \rangle$
\sn
\item[$(e)$]   $\lambda_\delta = \max \pcf\{\theta_i:i \in
u\cap \delta \setminus j\}$ for every $j \in [j_\delta,\delta)$
\sn
\item[$(f)$]  $\lambda = \tcf(\prod\limits_{i \in u}
\theta_i,<_{J^{\bd}_u})$
\sn
\item[$(g)$]  ${\cD}$ is a normal filter on $\kappa$ extending ${\cD}_\kappa$
\sn
\item[$(h)$]   if $A \in {\cD}^+$, $v_\delta \subseteq u \cap \delta$, 
$j_\delta < \delta$, $\lambda_\delta > \max \pcf\{\theta_i : i \in u \cap \delta \setminus v_\delta \setminus j_\delta\}$ for
$\delta \in A$ then $\bigcup\{v_\delta:\delta \in A\}$ is a co-bounded
subset of $u$.
\end{enumerate}
\end{claim}

\begin{remark}
We can add:
\mn
\begin{enumerate}
\item[$(i)$]   if $v$ is an unbounded subset of $u$ then the set
$$\Big\{i < \kappa : \max \pcf \big( \{\theta_j:j \in i \cap v\} \big) = 
\max\pcf \big(\{\theta_j : j \in i \cap u\} \big) \Big\}$$ 
belongs to ${\cD}$.
\end{enumerate}
\end{remark}

\begin{PROOF}{\ref{k.5}}
By the $\pcf$ theorem there is $u_0 \in [\kappa]^\kappa$ such that
\mn
\begin{enumerate}
\item[$(*)$]   $\lambda = \tcf\big(\prod\limits_{i \in u_0}
\theta_i, {<_{J^{\bd}_{u_0}}}\big)$ is well defined.
\end{enumerate}
\mn
Now for every $u \in [u_0]^\kappa$ we define $E_u$ and $\langle
\lambda^u_\delta:\delta \in E_u\rangle$ as in clauses (c),(e) and
stipulate $\lambda^u_i=0$ for $i \in \kappa \setminus E_u$ and let
$\bar \lambda^u = \langle  \lambda^u_i:i <\kappa \rangle$.  So
$\gamma_u = \rk_{\cD_\kappa}(\langle \lambda^u_i:i <\kappa\rangle)$ is 
a well defined ordinal and we can choose $u_1 \in
[u]^\kappa$ such that $\gamma_{u_1}$ is minimal.  Let 
$${\cD}^*_u = \big\{A \subseteq \kappa:A \in {\cD}_\kappa \text{ or } A \in {\cD}^+_\kappa 
\setminus {\cD}_\kappa \text{ and } \gamma_{u_1} < \rk_{D+(\kappa \setminus A)} 
(\langle \lambda^u_i:i < \kappa \rangle) \big\}.$$

As for clause (h), but \cite{Sh:589}, ${\cD}_u$ is a normal filter
on $\kappa$ (extending ${\cD}_\kappa$).  For proving \magenta{[?]} assume that
$A\in {\cD}^+_u$, $\bar v = \langle v_\delta : \delta \in A\rangle$, 
$\bar j = \langle j_\delta : \delta \in A\rangle$ and $v_\delta \subseteq u_1 \cap \delta$,
$j_\delta \cap \delta$, $j_\delta < \delta$ and
$\lambda^{u_1}_\delta > \max \pcf\{\theta_i:i \subseteq
\delta \cap u_1 \setminus v_\delta \setminus j_\delta\}$.

We should prove that $v := u_1 \setminus \bigcup\{v_\delta:\delta \in
A\}$ is bounded in $\kappa$.  Toward contradiction assume $\kappa =
\sup(v)$ and we shall prove that $\gamma_v < \gamma_u$, thus deriving
the desired contradiction
\mn
\begin{enumerate}
\item[$(**)_1$]   $\gamma_{u_1} = \rk_{{\cD}_\kappa}(\bar \lambda^{u_1})$.
\end{enumerate}
\mn
But by the choice of ${\cD}_u$
\mn
\begin{enumerate}
\item[$(**)_2$]   $\rk_{{\cD}_\kappa}(\bar \lambda^{u_1}) = 
\rk_{\cD_\kappa +A}(\bar \lambda^{u_1})$.
\end{enumerate}
\mn
Now clearly by our assumption
\mn
\begin{enumerate}
\item[$(**)_3$]  $\delta \in A \Rightarrow \lambda^v_\delta <
\lambda^{u_1}_\delta$
\end{enumerate}
\mn
hence
\mn
\begin{enumerate}
\item[$(**)_4$]  $\bar \lambda^v < \bar \lambda^{u_1} 
\mod({\cD}_\kappa +A)$ hence
\sn
\item[$(**)_5$]   $\rk_{{\cD}_\kappa +A}(\bar \lambda^{u_1}) >
\rk_{\cD_\kappa +A}(\bar \lambda^v)$.
\end{enumerate}
\mn
Now by a monotonicity property of $\rk_D(\bar \lambda^v)$ in $D$
\mn
\begin{enumerate}
\item[$(**)_6$]   $\rk_{\cD_\kappa +A}(\bar \lambda^v) \ge
\rk_{\cD_\kappa}(\bar \lambda^v)$.
\end{enumerate}
\mn
But
\mn
\begin{enumerate}
\item[$(**)_7$]   $\rk_{\cD_\kappa}(\bar \lambda^v) = \gamma_v$.
\end{enumerate}
\mn
Together $(**)_1 - (**)_7$ gives $\gamma_{u_1} > \gamma_v$,
contradicting the choice of $u_1$.  The contradiction comes from
assuming that $v$ is unbounded in $\kappa$, so $\sup(v) < \kappa$, thus
finishing the proof of clause (h) and of the claim.
\end{PROOF}

\begin{remark}
\label{k.6}
We can replace $(J^{\bd}_\kappa,\cD_\kappa)$ by other such 
pairs (on $\kappa$ or on $[\mu]^{< \kappa}$).
\end{remark}

\begin{observation}
\label{k.7}
Assume $\theta = \cf(\theta)$ and $\bar \lambda = 
\langle \theta_i:i < \kappa \rangle$ is an increasing
sequence of regular cardinals $> \kappa^{++}$ and $\lambda = 
\tcf(\prod \theta_i,<_{J^{\bd}_\sigma})$.  Then we can find an $u,{\cF},
\bar f$ such that
\mn
\begin{enumerate}
\item[$(a)$]   $u \subseteq \kappa$ is unbounded
\sn
\item[$(b)$]   $\bar f = \langle f_\alpha:\alpha < \lambda\rangle$
such that
\sn
\item[$(c)$]   $f_\alpha \in \prod\limits_{i \in u} \theta_i$
\sn
\item[$(d)$]  $\langle f_\alpha:\alpha < \lambda\rangle$ is
$<_{J^{\bd}_u}$-increasing cofinal in $(\prod\limits_{i \in u}
\theta_i,<_{J^{\bd}_u})$
\sn
\item[$(e)$]   ${\cF} \subseteq \prod\limits_{i \in u} \theta_i$
includes $\{f_\alpha:\alpha < \lambda\}$ and 
$\big|\{f \rest \delta : f \in {\cF}\}\big| \le \lambda_\delta$ for $\delta \in E$
\sn
\item[$(f)$]   if $f \in \prod\limits_{i \in u} \theta_i$ for every
$j < \theta$ for some $g \in {\cF}$ we have $f_\alpha \restriction
(j \cap u) = g \restriction (j \cap u)$ then $f \in {\cF}$
\sn
\item[$(g)$]   ${\cF}$ is linearly ordered by $<_{J^{\bd}_u}$.
\end{enumerate}
\mn
Let $u,E,\langle \lambda_i:i < \kappa \rangle$, ${\cD}$ be as in the
previous claim.  As we can ...?

For $j \in E_i$ let $J_j = \{v \subseteq u \cap j:\max \pcf(\lambda_i:
i \in [j',j) \cap u) < \lambda_i$ for some $j' < j\}$.
For each $\delta \in E_u$ choose $\langle f^\delta_\alpha:\alpha <
\lambda^u_j\rangle$ such that
\mn
\begin{enumerate}
\item[$\circledast_\delta$]  $(a) \quad f^\delta_\alpha \in 
\prod\limits_{i \in u \cap \delta} \theta_i$
\sn
\item[${{}}$]  $(b) \quad \bar f^\delta = \langle
f^\delta_\alpha:\alpha < \lambda_j\rangle$ is
$<_{J_\delta}$-increasing and cofinal in $(\prod\limits_{j \in u \cap
\delta} \theta_i,<_{J_i})$
\sn
\item[${{}}$]   $(c) \quad$ if $\bar f^j \restriction \delta$ has a
$<_{J_j}$-l.u.b. then $f^j_\delta$ is an increasing $<_{J_j}$-l.u.b.
\end{enumerate}
\mn
Let
\mn
\begin{enumerate}
\item[$\circledast$]  ${\cF}^* = \{f \in \prod\limits_{i \in u}
\lambda_i: (\forall\delta \in E) (\exists\alpha <
\lambda_\delta)\ f \restriction (u \cap \delta) =
f^\delta_\alpha \mod J_\delta\}$.
\end{enumerate}
\mn
For $f \in {\cF}$ let $g^+_f$ be the function with domain
$E$, $g^+_f(\delta) = \alpha$, $\alpha$ as above (clearly it is unique).  
\magenta{By \cite[xxx]{Sh:e} - FILL}
\mn
\begin{enumerate}
\item[$\circledast$]   we can find $f_\alpha \in {\cF}$ for
$\alpha < \lambda$ such that $\bar f = \langle f_\alpha:\alpha <
\lambda\rangle$ is $<_{J^{\bd}_u}$-increasing cofinal in
$(\prod\limits_{i \in u} \lambda_i,<_{J^{\bd}_u})$ (and if
$\alpha > \kappa,\bar f \restriction \alpha$ has a
$<_{J^{\bd}_u}$-e.u.b. then $f_\alpha$ is such $<_{J^{\bd}_u}$-e.u.b..
\end{enumerate}
\mn
Assume toward contradiction
\mn
\begin{enumerate}
\item[$\boxtimes$]   $f_1,f_2 \in {\cF}$ and $u_1 = \{i \in
u:f_1(i) < f_2(i)\}$ is unbounded in $\kappa$ and also $u_2 := u
\setminus u_1$ is unbounded $\theta \kappa$.
\end{enumerate}
\mn
Now $E$ is partitioned to 
$$\begin{aligned}
A_1 = \{\delta \in E:g_{f_1}(\delta) < g_{f_2}(\delta)\}&\text{ and}\\ 
A_2 = \{\delta \in E:g_{f_1}(\delta) \ge g_{f_2}(\delta)\}&.
\end{aligned}$$

Hence for some $\ell \in \{1,2\}$ we have $A_\ell
\in D^+$.  So for each $\delta \in A_\delta$ we can find $v_\delta$
such that
\mn
\begin{enumerate}
\item[$(a)$]  $u \setminus v_\delta \in J_\delta$ and $v_\delta
\subseteq u \cap\delta$
\sn
\item[$(b)$]   $f_k \restriction v_\delta =
f^\delta_{g_{f_k}(\delta)} \restriction v_\delta$ for $k=1,2$
\sn
\item[$(c)$]   if $\ell = 1$ then $g_{f_1}(\delta) <
g_{f_2}(\delta)$ and $f^\delta_{g_{f_1}(\delta)} \restriction v_\delta
< f^\delta_{g_{f_2}(\delta)} < f^\delta_{g_{f_2}(\delta)} \restriction
v_\delta$
\sn
\item[$(d)$]   if $\ell = 2$ then $g_{f_1}(\delta) \ge
g_{f_2}(\delta)$ and $f^\delta_{g_{f_1}(\delta)} \restriction v_\delta
\ge f^\delta_{g_{f_2}(\delta)} \restriction v_\delta$.
\end{enumerate}
\mn
This is clearly possible.

Now if $i \in v := \bigcup\{v_\delta:\delta \in A_i\}$ then $[f_1(i) <
f_2(i) \Leftrightarrow \ell=1]$ but by the previous claim (clause (b))
and clause (a), $v$ is a co-bounded subset of $u,f_1 < f_2$ mod
$J^{\bd}_u$ or $f_2 \le f_1 \mod J^{\bd}_u$ so we are done.
\end{observation}

\begin{conclusion}
\label{k.8}
Assume $\mu > \kappa = \cf(\mu) > \aleph_0,\langle \mu_i:
i < \kappa \rangle$ is increasing continuous
sequence with limit $i,\cf(\mu_i) \le \kappa$ and $\pp(\mu_i) <
\mu_{i+1}$ for $i < \kappa$.  \underline{then}  we can find ${\cF}$ as in
\ref{k.7} of cardinality (and cofinality) $\pp(\mu)$.
\end{conclusion}
\newpage

\section {More on III,4.10: Densely running away from Colours} \label{C3}

\begin{question}
\label{ac.1}
[Hajnal]:  Let $\lambda = (2^{\aleph_0})^+$.  Is there 
$c:[\lambda]^2 \rightarrow \omega$ such that

\[
(\forall A \in [\lambda]^\lambda)(\forall n < \omega)(\exists B \in
[A]^\lambda)\big[n \notin \Rang(c \restriction [B]^2)\big]?
\]

\mn
Answer: yes.

Clearly it is equivalent to the property
$P_7(\lambda,\aleph_0,2)$ defined below for $\lambda = (2^{\aleph_0})^+$.
Now Claim \ref{ac.3} covers the case $\lambda = (2^{\aleph_0})^+$ and 
then we have more.  We look again at 
\cite[Ch.III,4.9-4.10C,pp.177-181]{Sh:e}.
\end{question}

\begin{definition}
\label{ac.2}
$\Pr_7(\lambda,\sigma,\theta)$ where $\lambda \ge \theta \ge 1,
\lambda \ge \sigma = \cf(\sigma)$ means that
there is $c:[\lambda]^2 \rightarrow \sigma$ such that  

\[
(\forall A \in [\lambda]^\lambda)(\forall \alpha < \sigma)
(\exists B \in [A]^\lambda)\big[\min \Rang(c \rest [B]^2) > \alpha\big]
\]

\mn
(So far, $\theta$ is redundant).
Moreover, if $w_\alpha \in [\lambda]^{< 1 + \theta}$ for $\alpha < \lambda$
are pairwise disjoint and $\zeta < \sigma$ \underline{then}  for some $X \in 
[\lambda]^\lambda$ we have
\mn
\begin{enumerate}
\item[$(*)$]   if $\alpha < \beta$ are from $X$ then 
$(\forall i \in w_\alpha)(\forall j \in w_\beta)(c\{i,j\} \ge \zeta)$.
\end{enumerate}
\end{definition}

\begin{claim}
\label{ac.3}
Assume $\lambda$ is a regular uncountable cardinal, 
$2 \le \kappa < \lambda$ and $\otimes^\kappa_\lambda$ holds or just
$\oplus^\kappa_\lambda$ (see below). 

\underline{then}  there is a symmetric 2-place function $c$ from $\lambda$ to
$\aleph_0$ such that:
\mn
\begin{enumerate}
\item[$(*)$]   if $\langle w_i:i < \lambda \rangle$ is a sequence of
pairwise disjoint non-empty subsets of $\lambda,|w_i| < \kappa$ and $n <
\omega$, \underline{then}  for $Y \in [\lambda]^\lambda$ for every $i < j$ from $Y$
we have:

\[
\max(w_i) < \min(w_j)
\]

\[
\bigwedge\limits_{\alpha \in w_i} \bigwedge\limits_{\beta \in w_j}
c(\alpha,\beta) > n.
\]
\end{enumerate}
\mn
(i.e. $\Pr_7(\lambda,\aleph_0,\kappa))$. 

Note that Definition \ref{ac.3A}(1)
is from \cite[Ch.III,4.10,p.178]{Sh:e}.
\end{claim}

\begin{definition}
\label{ac.3A}
1) For a Mahlo (inaccessible) cardinal $\lambda$ and $\kappa < \lambda$ let
\mn
\begin{enumerate}
\item[$\otimes^\kappa_\lambda$]   there is $\bar C = \langle C_\delta:
\delta \in S^\lambda_{\in} \rangle$, where 
$S^\lambda_{\in} =: \{\delta < \lambda:\delta \text{ is inaccessible}\},
C_\delta$ a club of $\delta$, such that: 
for every club $E$ of $\lambda$ for some $\delta \in \acc(E) 
\cap S^\lambda_{\in}$ of cofinality $\ge \kappa$, 
for \underline{no} $\zeta < \kappa$ and $\alpha_\eps \in 
S^\lambda_{\in}$ (for $\eps < \zeta)$ do we have
\sn
\begin{enumerate}
\item[$(*)$]   $\nacc(E) \cap \delta \setminus \bigcup\limits_{\eps <
\zeta} C_{\alpha_\eps}$ is bounded in $\delta$.
\end{enumerate}
\end{enumerate}
\mn
2) For $\lambda$ regular $> \kappa = \cf(\kappa) \ge \aleph_0$, let 
\mn
\begin{enumerate}
\item[$\oplus^\kappa_\lambda$]  there is $\bar C = \langle C_\delta:
\delta \in S \rangle,S = \{\delta < \lambda:\delta \mathrm{limit} \},C_\delta$
a club of $\delta$ such that: for every club $E$ of $\lambda$ for some
$\delta \in \acc(E)$ of cofinality $\ge \kappa$, for \underline{no} 
$\zeta < \kappa$ and $\alpha_\eps \in S$ (for 
$\eps < \zeta$) do we have
\sn
\begin{enumerate}
\item[$(*)'$]  $S^\lambda_{\ge \kappa} \cap E \setminus
\bigcup\limits_{\eps < \zeta} C_{\alpha_\eps}$ is bounded 
in $\delta$ where $S^\lambda_{\ge \kappa} = \{\delta < \lambda:
\cf(\delta) \ge \kappa\}$.
\end{enumerate}
\end{enumerate}
\end{definition}

\begin{remark}
\label{ac.3B}
1) For $\lambda$ Mahlo, the property 
$\otimes^2_\lambda$ holds \If \, there are stationary
subsets $S_i$ of $\lambda$ for $i < \lambda$ such that for no $\delta <
\lambda,\bigwedge\limits_{i < \delta} [S_i \cap \delta$ a stationary in
$\delta$] (we can consider only $\delta$ inaccessible). 

[Why?  Choose $C_\delta$ a club of $\delta$ disjoint to $S_i$ for some
$i(\delta) < \delta$, such that $\min(C_\delta) > i(\delta)$]. 

\noindent
2) This is close to \cite[\S3]{Sh:276}, see \cite[Ch.III,2.12]{Sh:g}.  As in
\cite[\S3]{Sh:276}, the proof is done such that from appropriate failures of
Chang conjectures or existence of colourings we can get stronger colourings
here.  For the result as stated also $c(\beta,\alpha) = \ell g[\rho(\beta,
\alpha)]$ is O.K., but the proof as stated is good for utilizing failure of
Chang conjecture (as in \cite[\S3]{Sh:276}). 

\noindent
3) Note that $\otimes^2_\lambda$ is closely related to $\otimes_{\bar C}$
from \cite[Ch.III,2.12]{Sh:g}.  
Also if $\kappa \le \aleph_0$, then in $\otimes^\kappa_\lambda$ 
we can replace nacc$(E)$ by $E$. 

\noindent
4) Note that $\lambda$ weakly compact fails 
even $\otimes^2_\lambda$ and forcing notion $P$ which is 
$\theta$-c.c. for some $\theta < \lambda$ preserves this.
\end{remark}

\begin{observation}
\label{ac.3C}
In Definition \ref{ac.3A} in $(*)$ and $(*)'$ if
$\kappa \le \aleph_0$ it does not matter whether we write $E$ or 
$\nacc(E)$. 
\end{observation}

\begin{observation}
\label{ac.3D}
1) $\otimes^2_\lambda$ implies $\otimes^{\aleph_0}_\lambda$. 

\noindent
2)  $\oplus^2_\lambda$ implies $\oplus^{\aleph_0}_\lambda$. 

\noindent
3)  If $\kappa_1 < \kappa_2 < \lambda$ then $\otimes^{\kappa_2}_\lambda
\Rightarrow \otimes^{\kappa_1}_\lambda$ and $\oplus^{\kappa_2}_\lambda
\Rightarrow \oplus^{\kappa_1}_\lambda$. 

\noindent
4) $\otimes^\kappa_\lambda \Rightarrow \oplus^\kappa_\lambda$ if $\lambda$
is inaccessible $> \aleph_0$.
\end{observation}

\begin{PROOF}{\ref{ac.3C}}
1) Let $\bar C$ exemplify $\otimes^2_\lambda$ and we shall show
that it exemplifies $\otimes^{\aleph_0}_\lambda$, assume not and let $E$ be
a club of $\lambda$ which exemplifies this.  We choose by induction on
$k < \omega$ a club $E_k$ of $\lambda:E_0 = E$, if $E_k$ is defined let

\begin{equation*}
\begin{array}{clcr}
A_k =: \{\delta < \lambda:&\delta \in \acc(E_k) \cap S^\lambda_{\in} 
\text{ and for no} \\
  &\alpha \in S^\lambda_{\in} \text{ is } E_k \cap \delta \setminus
C_\alpha \text{ bounded in } \delta\}.
\end{array}
\end{equation*}

\mn
As $\bar C$ exemplifies $\otimes^2_\lambda$, clearly $A_k$ is a stationary
subset of $\lambda$ and let

\[
E_{k+1} = \{\delta \in E_k:\delta = \sup(A_k \cap \delta)\}.
\]

\mn
Let $\delta(*) \in \bigcap\limits_{k < \omega} E_k$ which 
necessarily belong $\subseteq E$.  
By the choice of $E$ we can find $n < \omega = \kappa$ and $\alpha_\ell 
\in S^\lambda_{\in}$ for $\ell < n$ such that $\nacc(E) \cap \delta(*) 
\setminus \bigcup\limits_{\ell < n} C_{\alpha_\ell}$ is bounded in 
$\delta(*)$.   Now we choose by induction on $k \le n,\delta_k \in
\acc(E_{n+1-k})$ such that $\delta_k < \delta(*)$ and
$\nacc(E_{n+1-k}) \cap \delta_k \setminus \bigcup\limits_{\ell < n-k} 
C_{\alpha_\ell}$ is bounded in $\delta_k$.  For $k=0$ any large enough 
$\delta \in \delta(*) \cap E_{n+1}$ is O.K.  For $k+1$ use the
definition of $E_{n+1-k}$.  For $k=n,\delta_n$ gives a contradiction
to the choice of $E$.  

\noindent
2) Same proof replacing $S^\lambda_{\in}$ by $S^\lambda_{\ge \kappa}$.

\noindent
3) The same $\bar C$ witnesses it. 

\noindent
4) Here $\lambda$ is inaccessible.  That is, we have to show that: 

\underline{the version with $(*) \Rightarrow$ the version with $(*)'$}

Let $\bar C' = \langle C'_\delta:\delta \in S^\lambda_{\in} \rangle$
exemplifies $\otimes^\kappa_\delta$.  We define $S = \{\delta < \lambda:
\delta\ \mathrm{limit} \}$ and $\bar C = \langle C_\delta:\delta \in S \rangle$
as follows: if $\delta \in S^\lambda_{\in} (\subseteq S)$ we let $C_\delta =
C'_\delta$ and if $\delta \in S \setminus S^\lambda_{\in}$ let $C_\delta$ be
a club of $\delta$ of order type $\cf(\delta)$ with 
$\cf(\delta) < \delta \Rightarrow \min(C_\delta) > \cf(\delta)$
and if $\delta$ is a successor cardinal, say $\theta^+$ then
$\min(C_\delta) > \theta$ (possible as $\delta \notin S^\lambda_{\in} 
\Rightarrow \cf(\delta) < \delta \vee 
(\exists \theta < \delta)(\delta = \theta^+))$.
We shall show that $\langle C_\delta:\delta \in S 
\rangle$ exemplify $\oplus^\kappa_\lambda$.

Given a club $E$ of $\lambda$, let 
$$E_0 = \{\delta \in E:\delta \text{ a limit
cardinal, } \otp(\delta \cap E) = \delta \text{ and } \delta > \kappa\}$$ 
$$\text{and } E_1 = \{\delta \in E_0:\otp(\delta \cap E_0)\text{ is divisible by }
\kappa^+\},$$ so $E_1$ is a club of $\lambda$ so by the
version with $(*)$ there is $\delta \in \acc(E_1) \cap 
S^\lambda_{\in}$ hence $\cf(\delta) > \kappa$ satisfying $(*)$, i.e. 
the requirement in \ref{ac.3A}(1); we shall
show that it satisfies the requirement in \ref{ac.3A}(2) thus finishing.

So let $\zeta < \kappa$ and $\alpha_\eps \in S$ for $\eps <
\zeta$ and we should prove that $Y =: S^\lambda_{\ge \kappa} \cap E
\setminus \bigcup\limits_{\eps < \zeta} 
C_{\alpha_\eps}$ is unbounded
in $\delta$, so fix $\beta^* < \delta$ and we shall prove that $Y \cap
(\beta^*,\delta) \ne \varnothing$ thus finishing.

Let $\zeta $ be the disjoint union of $u_0,u_1,u_2$, where 
$$\begin{aligned}
u_0 =& \{\eps < \zeta:\alpha_\eps < \delta\}\\
u_1 =& \{\eps < \zeta: \alpha_\eps \ge \delta \text{ and } 
\alpha_\eps \in S \setminus S^\lambda_{\in}\}\\
u_2 =& \{\eps < \zeta:\alpha_\eps \ge \delta 
\text{ and } \alpha_\eps \in S^\lambda_{\in}\}.
\end{aligned}$$

By the choice of $\delta$ we know that $Y_2 = \nacc(E_1) \cap \delta
\setminus \bigcup\limits_{\eps \in u_2} C_{\alpha_\eps}$ is unbounded
in $\delta$.  As $\cf(\delta) \ge \kappa$ (see its choice, i.e. 
$\delta \in \acc(E) \cap S^\lambda_{\in} \wedge \min(E) > \kappa)$, 
we can find $\beta \in Y_2$ such that $\beta < \delta$,
$\beta > \beta^*$ and $\beta > \alpha_\eps$ for $\eps \in u_0$.

Now $\cf(\beta) = \kappa^+$ as $\beta 
\in Y_2 \subseteq \nacc(E_1)$ and by the choice of $E_1$.  
Also $$\eps \in u_0 \Rightarrow \sup(C_{\alpha_\eps}) < \beta\ 
\text{ and }\ \eps \in u_2 \Rightarrow \sup(C_{\alpha_\eps} \cap \beta) 
< \beta$$ (as otherwise $\beta \in C_{\alpha_\eps}$, contradicting 
$\beta \in Y_2$), so we can find $\beta_0 < \beta$ such that 
$\eps \in u_0 \cup u_2 \Rightarrow \sup(C_{\alpha_\eps} \cap \beta) < \beta_0$.  
Now for $\eps < \zeta$, if 
$C_{\alpha_\eps} \cap (\beta_0,\beta) \ne \varnothing$ then
$\eps \in u_1$, so by the choice of $C_{\alpha_\eps}$ we know
$|C_{\alpha_\eps}| = \cf(\alpha_\eps) < \min(C_{\alpha_\eps}) < \beta$, 
noting that $\beta$ is a cardinal as $E_0$ is a set of cardinals.  
By the definition of $E_0,E_1$ we know that 
$E \cap S^\lambda_{\ge \kappa} \cap \beta$ has cardinality $\beta$ 
hence $E \cap S^\lambda_{\ge \kappa} \setminus \beta_0$ has cardinality 
$\beta$, so we finish.  
\end{PROOF}

\begin{PROOF}{\ref{ac.3}}
By \ref{ac.7}(4) \wilog \, $\oplus^\kappa_\lambda$, so let 
$\bar C$ be as required in $\oplus^\kappa_\lambda$.

We define $e_\alpha$ for every ordinal $\alpha < \lambda$ as follows:
\mn
\begin{enumerate}
\item[$(a)$]   if $\alpha = 0$, $e_\alpha = \varnothing$
\sn
\item[$(b)$]  if $\alpha = \beta +1$, $e_\alpha = \{0,\beta\}$
\sn
\item[$(c)$]  if $\alpha$ is a limit ordinal, then we let 
$e_\delta = C_\delta \cup \{0\}$.
\end{enumerate}
\mn
Let $S$ be the set of limit ordinals $< \lambda$.
For $\alpha < \beta$ we define by induction on $\ell < \omega$ the ordinals
$\gamma^+_\ell(\beta,\alpha)$, $\gamma^-_\ell(\beta,\alpha)$.
\medskip

\noindent
\underline{$\ell = 0$}:  $\gamma^+_\ell(\beta,\alpha) = \beta$, 
$\gamma^-_\ell(\beta,\alpha) = 0$
\medskip

\noindent
\underline{$\ell = k+1$}:  $\gamma^+_\ell(\beta,\alpha) = \min
(e_{\gamma^+_k(\beta,\alpha)} \setminus \alpha)$ if $\alpha < \gamma^+_k
(\beta,\alpha)$ and $\gamma^-_\ell(\beta,\alpha) = 
\sup(e_{\gamma^+_k(\beta,\alpha)} \cap \alpha)$ \If \, $\alpha <
\gamma^+_k(\beta,\alpha)$ and $\alpha \notin 
\acc(e_{\gamma^+_k(\beta,\alpha)})$.

Note that $\gamma^-_\ell(\beta,\alpha) < \alpha \le \gamma^+_\ell
(\beta,\alpha)$ if they are defined and then $\ell > 0 \Rightarrow \gamma^+
_\ell(\beta,\alpha) < \gamma^+_{\ell -1}(\beta,\alpha)$ (prove by induction).
So if $\alpha < \beta < \lambda$ for some
$k = k(\beta,\alpha) < \omega$ we have: $\gamma^+_\ell(\beta,\alpha)$ is
defined iff $\ell \le k$ and: $\gamma^-_\ell(\beta,\alpha)$ is defined iff
$\ell < k \vee [\ell = k \text{ and } \gamma^+_k(\beta,\alpha) = \alpha]$ and:
$\gamma^+_k(\beta,\alpha) = \alpha$ or $\alpha \in \acc
(e_{(\gamma^+_k(\beta,\alpha))})$.  Let $\rho(\beta,\alpha) =
\langle \gamma^+_\ell(\beta,\alpha):\ell \le k(\beta,\alpha) \rangle$.  Note
(we shall use it freely):
\mn
\begin{enumerate}
\item[$\otimes_1$]   if $\gamma < \alpha < \beta,k \le k(\beta,\alpha)$
and $\gamma^-_k(\beta,\alpha)$ is defined and 
$\bigwedge\limits_{\ell \le k} \gamma^-_\ell(\beta,\alpha) < \gamma$
\underline{then} 
\sn
\begin{enumerate}
\item[$(\alpha)$]  $\ell \le k \Rightarrow \gamma^+_\ell(\beta,\alpha)
= \gamma^+_\ell(\beta,\gamma)$
\sn
\item[$(\beta)$]  $\ell \le k \Rightarrow \gamma^-_\ell(\beta,\alpha) =
\gamma^-_\ell(\beta,\gamma)$
\sn
\item[$(\gamma)$]  $k(\beta,\gamma) \ge k(\beta,\alpha) 
\text{ and } \rho(\beta,\alpha) \trianglelefteq \rho(\beta,\gamma)$.
\end{enumerate}
\end{enumerate}
\mn
Now we define $c\{\alpha,\beta\} = c(\beta,\alpha) = c(\alpha,\beta)$
for $\alpha < \beta < \lambda$ as follows:

\[
c(\beta,\alpha) = k(\beta,\alpha) + 1.
\]

\mn
So assume $\bar w = \langle w_i:i < \lambda \rangle$ is a sequence of
pairwise disjoint subsets of $\lambda$, $|w_i| < \kappa$ and $n(*) < \omega$.
Without loss of generality for some $\kappa^* < 1 + \kappa$, 
$\bigwedge\limits_{i < \lambda} |w_i| = \kappa^*$ and $i < \min(w_i)$ 
and $[i < j \Rightarrow \sup(w_i) < \min(w_j)]$.  
Let $w_i = \{\alpha^i_\eps : \eps < \kappa^*\}$.  Let $\chi \ge (2^\lambda)^+$, and we choose by induction on $n < \omega$ and for each $n$ by induction on $i < \lambda$, 
$N^n_i \prec (\clH(\chi),\in,<^*_\chi)$ such that $\|N^n_i\| < \lambda$, 
$\{\langle N^n_\eps: \eps \le j \rangle:j < i\} \subseteq N^n_i$, $\bar w \in N^0_i$, 
$N^n_i$ increasing continuous in $i$ and $\langle N^m_i : i < \lambda \rangle \in 
N^n_0$ for $m < n$.

Let us define for $\ell < \omega$

\[
E^\ell = \{\delta < \lambda:N^\ell_\delta \cap \lambda = \delta\}
\]

\begin{equation*}
\begin{array}{clcr}
S^\ell = \{\delta \in S^\lambda_{\ge \kappa} \cap \acc 
(E^\ell):&\text{for \underline{no} } \zeta < \kappa \text{ and } 
\alpha_\eps < \lambda \\
  &\text{for }\eps< \zeta \text{ do we have} \\
  &\delta > \sup[S^\lambda_{\ge \kappa} \cap E^\ell \cap \delta \setminus
\bigcup\limits_{\eps < \zeta} C_{\alpha_\eps}]\}.
\end{array}
\end{equation*}

\mn
Note that $\alpha < \lambda \Rightarrow (E^\ell,S^\ell) \in N^{\ell +1}
_\alpha$ hence $\delta \in E^{\ell +1} \Rightarrow \delta = \sup(\delta
\cap S^\ell)$.

We know that $S^\ell$ is a stationary subset of $\lambda$ 
as $E^\ell$ is a club of $\lambda$ because $\oplus^\lambda_\kappa$ is 
exemplified by $\bar C$. 

Choose $\delta_{n(*)} \in E^{2(n(*)+1)} \cap S^{2(n(*)+1)}$ and then choose
$\alpha(*) < \lambda$ such that $\alpha(*) > \delta_{n(*)}$.  We
now choose by downward induction on $m < n(*)$ ordinals $\delta_m,
\zeta^*_m$ such that:
\mn
\begin{enumerate}
\item[$(*)(i)$]  $\delta_m < \delta_{m+1}$
\sn
\item[$(ii)$]  $\delta_m \in E^{2m} \cap S^{2m}$
\sn
\item[$(iii)$]  $\delta_m > \sup\{\gamma^-_\ell(\beta,\delta_{m+1}) : \beta \in w_{\alpha(*)},\ 
\ell \le k(\beta,\delta_{m+1}),\ \gamma^-_\ell(\beta,\delta_{m+1}) \text{ well defined} \}$

\item[$(iv)$]  $\delta_m \notin \bigcup \{C_\gamma : \gamma = \gamma^+_{k(\beta, \delta_{m+1})}(\beta, \delta_{m+1})$ for some $\beta \in w_{\alpha(*)}\}$
\sn
\item[$(v)$]   $\zeta^*_m < \delta_m,\zeta^*_m < \zeta^*_{m+1}$ if $m+1
< n(*)$
\sn
\item[$(vi)$]   if $\alpha \in [\zeta^*_m,\delta_m)$ then $(\forall
\beta' \in w_\alpha)(\forall \beta'' \in w_{\alpha(*)})
\big[\rho(\beta'',\delta_m) \triangleleft \rho(\beta'',\beta')\big]$
\end{enumerate}
\mn
[Why can we do it?  Assume $\delta_{m+1} \in S^{2(m+1)}$ has already been
defined and we shall find $\delta_m,\zeta_m$ as required.
Let 
$$Y_m = \{\gamma^-_\ell(\beta,\delta_{m+1}) : \beta \in w_{\alpha(*)},\ 
\ell \le k(\beta,\delta_{m+1}),\ \gamma^-_\ell(\beta,\delta_{m+1}) \text{ well defined} \}$$ 
so $Y_m$ is a subset of $\delta_{m+1}$ of cardinality $< \kappa$, but 
$\delta_{m+1} \in S^{2(m+1)}$ (if $m = n(*)-1$ by the choice of $\delta_n$, 
if $m < n-1$ by the induction hypothesis).  But\\ $S^{2(m+1)} \subseteq S^\lambda_{\ge \kappa}$, 
hence $(\forall \delta \in S^{2(m+1)})[\cf(\delta) \ge \kappa]$, hence
sup$(Y_m) < \delta_{m+1}$.  Also as $\delta_{m+1} \in E^{2(m+1)} \cap
S^{2(m+1)}$ by the definition of $S^{2(m+1)}$, there is 

$$
    \xi^*_m \in  S^\lambda_{\ge \kappa} \cap E^{2(m+1)} \cap
    \delta \setminus \bigcup \big\{ {\le_\gamma} : (\exists\beta \in w_{\alpha(*)}) [ \gamma = \gamma^+_{k(\beta, \delta_{m+1})} (\beta, \delta_{m+1})] \big\} \setminus \sup(Y_m)
$$

As each $e_\gamma$ is closed and there are $< \kappa$ of them,$\zeta^*_m< \xi^*_m$. where
$$\zeta^*_m = \sup\Big\{ \{\sup Y_m\} \cup \big\{\sup (e_\gamma \cap \xi^*_m):
(\exists\beta \in w_{\alpha(*)}) [\gamma = \gamma^+_{k(\beta, \delta_{m+1})} (\beta, \delta_{m+1})] \big\} \Big\}$$ 
So we can find $\delta_m \in (\zeta^*_m,\xi^*_m) \cap 
S^\lambda_{\ge \kappa} \cap E^{2m} \cap S^{2m}$ as required and 
choose $\zeta_m < \delta_m$ large enough.]
\mn
\begin{enumerate}
\item[$(**)$]  For every $\alpha \in [\zeta^*_0,\delta_0)$ we have

\[
(\forall \beta' \in w_\alpha)(\forall \beta'' \in w_\alpha)
[c\{\beta',\beta''\} \ge n].
\]
\end{enumerate}
\mn
[Why?  By clause (vi) above.]

Let 

\begin{equation*}
\begin{array}{clcr}
W = \Big\{\delta < \lambda:&\delta > \zeta^*_0 \text{ and for some }
\alpha'' \ge \delta \text{ we have} \\
  &\text{for every } \alpha' \in (\zeta^*_0,\delta) \text{ we have} \\
  &(\forall \beta' \in w_{\alpha'})(\forall \beta'' \in w_{\alpha''})
\big[c\{\beta',\beta''\} \ge n\big]\Big\}.
\end{array}
\end{equation*}

\mn
As $\delta_0 \in E_0$ (see $(*)(ii)$) so by $E_0$'s definition, $\delta_0
= N^0_{\delta_0} \cap \lambda$ hence $\zeta_0 \in N^0_{\delta_0}$.
Now $\bar w \in N^0_{\delta_0}$ (read definition) hence $W \in 
N^0_{\delta_0}$ and by $(*) + (**)$ and $W$'s definition
$\delta_0 \in W$, hence $W$ is a stationary subset of $\lambda$.   For
$\delta \in W$, let $\alpha''(\delta)$ be as in the definition of $W$.  So
$E = \big\{\delta^* : (\forall \delta \in W \cap \delta^*)[\alpha''(\delta) <
\delta^*]\big\}$, it is a club of $\lambda$ hence $W' = W \cap E$ is a stationary
subset of $\lambda$ and $\{\alpha''(\delta):\delta \in W'\}$ is as required.
\end{PROOF}

\begin{conclusion}
\label{ac.5}
If $\lambda = \cf(\lambda) > \aleph_0$ is not Mahlo (or is Mahlo as in 
\ref{ac.3A}(1) or \ref{ac.3A}(2)), $\kappa$ \underline{then}  
$\Pr_7(\lambda,\aleph_0,\aleph_0)$.
\end{conclusion}

\begin{PROOF}{\ref{ac.5}}
By \ref{ac.3} it suffices to prove $\oplus^{\aleph_0}_\lambda$.  This holds by \ref{ac.5A}, \ref{ac.7} and \ref{ac.8} below.
\end{PROOF}

\begin{claim}
\label{ac.5A}
1) If $\lambda = \mu^+$ \underline{then}  $\oplus_\lambda^{\cf(\mu)}$. 

\noindent
2) If $\lambda$ is (weakly) inaccessible, not Mahlo or Mahlo as in
\ref{ac.3A}(1), e.g. as in \ref{ac.3B}(1), and $\aleph_0 \le \kappa < 
\lambda$ \underline{then}  $\oplus_\lambda^\kappa$.
\end{claim}

\begin{PROOF}{\ref{ac.5A}}
1) Choose $C_\delta$ a club of $\delta$ of order type $\cf(\delta)$. 

Repeat the proof of \ref{ac.3C}(2), using 
$$E_0 = \{\delta < \lambda : \delta > \mu \text{ and } 
\otp(E \cap \delta) = \delta \text{ is divisible by } \mu^2\}.$$  
The only point slightly different is $|C_{\alpha_\eps} \cap (\beta_0,
\beta)| < |\beta|$ (now $\beta$ is not a cardinal).  For $\mu$ singular,
$|C_{\alpha_\eps}| < \mu = |\beta| = |\beta \cap 
S^\lambda_{\ge \kappa} \cap E \setminus \beta_0|$, and 
for $\mu$ regular we choose $\delta$ of
cofinality $\mu$ and everything is easy. 

\noindent
2) Now $\oplus^\lambda_\kappa$ holds trivially 
(choose a club $E^*_0$ of $\lambda$ with no inaccessible member 
and choose $C_\delta$ a club of $\delta$ of order type $\cf(\delta)$ 
such that cf$(\delta) < \delta \Rightarrow \min(C_\delta) 
> \cf(\delta)$ and $\delta \notin E^* \Rightarrow \min(C_\delta) > 
\sup(E^* \cap \delta)$, now for any club $E$
choose $\delta \in \acc(E \cap E^*))$  So we can apply \ref{ac.3C}(2). 
\end{PROOF}

\begin{definition}
\label{ac.6}
$\Pr_8(\lambda,\mu,\sigma,\theta)$ means:

there is $c:[\lambda]^2 \rightarrow [\sigma]^{< \aleph_0} \setminus
\{\varnothing\}$ such that \Iff \, $w_\alpha \in [\lambda]^{< \theta}$ for
$\alpha < \lambda$ are pairwise disjoint and $\zeta < \sigma$ \underline{then}  for
some $Y \in [\lambda]^\mu$ we have $\alpha', \alpha'' \in Y$ and  
$$\alpha' < \alpha'' \Rightarrow (\forall \beta' \in w_{\alpha'},\ \forall 
\beta'' \in w_{\alpha''})\big[\zeta \in c\{\beta',\beta''\}\big].$$
\end{definition}

\begin{observation}
\label{ac.7}
Note that $\Pr_8(\lambda,\lambda,\sigma,\theta) \Rightarrow 
\Pr_7(\lambda,\sigma,\theta)$ because we can use 
$c'\{\alpha,\beta\} = \max[c\{\alpha,\beta\}]$.
\end{observation}

\begin{claim}
\label{ac.8}
1) If $\lambda$ is regular and $\aleph_0 \le \sigma \le \lambda$ \then
\, $\Pr_8(\lambda^+,\lambda^+,\sigma,\lambda)$. 

\noindent
2) If $\mu$ is singular, $\lambda = \mu^+$ and $\aleph_0 \le \sigma \le
\cf(\mu)$ \underline{then}  $\Pr_8(\lambda,\lambda,\sigma,\cf(\mu))$. 

\noindent
3) If $\lambda$ is inaccessible $> \aleph_0$, $S \subseteq \lambda$
stationary not reflecting in inaccessibles and $\sigma < \lambda$, 
$\theta = \min \{\cf(\delta):\delta \in S\}$ \underline{then}  $\Pr_8(\lambda,\lambda,\sigma, \theta)$.
\end{claim}

\begin{PROOF}{\ref{ac.8}}
The proofs in \cite[Ch.III,\S4]{Sh:g} gives this - in fact this
is easier. E.g. 

\noindent
1) Follows by Claim \ref{ac.3} (and \cite[Ch.III,4.2(2),p.162]{Sh:g})
but let us give some details.

Let $\bar e$ be as there (i.e. $\bar e = \langle e_\alpha : \alpha < \lambda^+ \rangle$, 
$e_0 = \varnothing$, $e_{\alpha+1} = \{\alpha\}$, $e_\delta$ a club
of $\delta$ of order type $\cf(\delta)$).  Let $h:\lambda^+ \rightarrow
\sigma$ be such that $(\forall \zeta < \sigma)(\exists^{\stat} \delta
< \lambda^+)(\cf(\delta) = \lambda \text{ and } h(\delta) = \zeta)$, 
$h_\alpha = h \rest e_\alpha$, $\bar h = \langle h_\alpha : \alpha < \lambda^+ \rangle$.

Let $\gamma(\beta,\alpha),\gamma_e(\beta,\alpha),\rho_{\bar h}$ be as there
(Stage A,p.164) and also the colouring $d$: for $\alpha < \beta < \lambda^+$

\[
d(\beta,\alpha) = \max\{h(\gamma_{\ell +1}(\beta,\alpha):
\gamma_{\ell +1}(\beta,\alpha) \text{ well defined}\}.
\]

\mn
By Stage B there the result should be clear.  
\end{PROOF}

\noindent
Hajnal has shown the following: 
\begin{theorem}
\label{ac.38}
Assume $\lambda = (2^{< \kappa})^+$, $\kappa = \cf(\kappa) > \omega$, $I$ is 
a normal ideal concentrating on $S_{\kappa,\lambda} = \{\alpha < \lambda:\cf(\alpha) = \kappa\}$, 
$\Lambda \subset [\lambda]^2$ is such that $\Lambda \cap [B]^2 \ne \varnothing$ for all $B \in I^+$ 
and $\Lambda = \bigcup\limits_{\eta < \xi} \Lambda_\eta$ for some $\xi < \kappa$.

Then there exist $I$ and $T$ such that $I \subset J$, $T \subset \xi$, $J$ is a
normal ideal and for all $\eta \in T$ and $B \in J^+$ we have 

\[
[B]^2 \cap \Lambda_\eta \ne \varnothing \text{ and }
G \cap [B]^2 \subset \bigcup\{\Lambda_\eta:\eta \in T\}.  
\]
\end{theorem}

\noindent
This comes from the following:
\begin{lemma}
\label{ac.41}
Assume $\lambda = (2^{< \kappa})^+,\kappa = \cf(\kappa) > \omega$.

$I$ is a normal ideal concentrating on $S_{\kappa,\lambda},P$ is a partial
order not containing decreasing sets of type $\kappa$. 

Assume further that

\[
p:{\clP}(\lambda) \rightarrow P \text{ and}
\]

\[
p(A) \le_p p(B) \text{ for } A \subset B.
\]

\mn
Then there is an $A \in I^+$ and a normal ideal $J \supset I$ satisfying
$B \in J$ iff $B \in I$ or $p(B) \prec_P p(A)$ for $B \subset A$.
\end{lemma}
\bigskip

\noindent
\underline{page 412}:  add at the end 

The following improves \cite[Ch.IX,5.12,p.410]{Sh:g}.
\begin{claim}
\label{ac.44}

\noindent
1) Assume
\mn
\begin{enumerate}
\item[$(a)$]   $\sigma = \cf(\sigma) > \aleph_0$
\sn
\item[$(b)$]  $\langle \lambda_i:i < \sigma \rangle$ increasing continuous,
$\lambda = \sup\{\lambda_i:i < \sigma\}$
\sn
\item[$(c)$]  $\sigma \le \theta < \lambda$ and $\sigma^{\aleph_0} < \lambda$
\sn
\item[$(d)$]  $\cov(\lambda_i,\lambda_i,\theta^+,2) < \lambda$ for
$i < \sigma$.
\end{enumerate}
\mn
\underline{then} 
\mn
\begin{enumerate}
\item[$(\alpha)$]   $\pp_\theta(\lambda) =^+ \cov(\lambda,\lambda,
\theta^+,2)$ and $\pp^{\cer}_{J^{\bd}_\sigma}(\lambda) = \cov(\lambda,
\lambda,\theta^+,2)^+$ (on $\pp^{\cer}$ see below).
\sn
\item[$(\beta)$]  $S^* = \{ \delta < \sigma:\cov(\lambda_\delta,
\lambda_\delta,\theta^+,2)^+ = \pp^{\cer}_{J^{\bd}_{\cf(\delta)}}
(\lambda_\delta)\}$ contains a club of $\sigma$. 
\end{enumerate}
\mn
2) Instead of ``$\sigma^{\aleph_0} < \lambda"$ it suffices
\mn
\begin{enumerate}
\item[$\otimes$]  for some club $C$ of $\sigma$, \Iff \, $i < \delta 
\in C,\delta$ of cofinality $\aleph_0$ and set ${\ga} 
\subseteq \lambda_\delta$ of cardinality $\le \lambda_i$ and 
${\ga}$ is a set of regular cardinals, \underline{then}  
$$\lambda >  \Big|\big\{\tcf(\textstyle\prod{\gb}/J^{\bd}_{\gb}) : \gb \subseteq {\ga},\ \sup({\gb}) 
= \lambda_\delta,\ \otp({\gb}) = \omega,\ \textstyle\prod \gb/J^{\bd}_\gb \text{ has true cofinality}\big\}\Big|.$$  
(So \wilog \, $\lambda_{\delta + 1}$ is above this cardinality.)
\end{enumerate}
\end{claim}

\begin{definition}
\label{ac.47}

Let $J$ be an ideal on some ordinal $\Dom(J)$.  We let 
\begin{align*}
    \pp^{\cer}_J(\lambda) = \min\big\{ \mu :\ & \mu \text{ regular } > \lambda,   \text{ and }\\ 
&  \sup\{\tcf \textstyle\prod\limits_t \lambda_t/J : \bar \lambda = 
\langle \lambda_t : t \in \Dom(J) \rangle\\ 
&\text{is strictly increasing with limit } \lambda\} < \mu\big\}.
\end{align*}
\end{definition}

\begin{PROOF}{\ref{ac.47}}
Proof of \ref{ac.47}
1) Similar to the proof of \cite[Ch.IX,5.12]{Sh:g}.
We assume toward contradiction that the desired conclusion fails.

Without loss of generality
\mn
\begin{enumerate}
\item[$(*)_0(a)$]   each $\lambda_i$ is singular of cofinality $< \sigma$
\sn
\item[$(b)$]   $\theta^{+3} < \lambda_0$ and $\sigma^{\aleph_0} < \lambda_0$
\sn
\item[$(c)$]  $\cov(\lambda_i,\lambda_i,\theta^+,2) < \lambda_{i+1}$
\sn
\item[$(d)$]   $\mu \in (\lambda_0,\lambda_{i+1}) \Rightarrow 
\pp_\theta(\mu) < \lambda_{i+1}$.
\end{enumerate}
\mn
[Why?  Clearly we can replace $\langle \lambda_i:i < \sigma \rangle$ by
$\bar \lambda \restriction C = \langle \lambda_i:i \in C \rangle$ for any
club $C$ of $\sigma$, hence it is enough to show that each of the demands
holds for $\bar \lambda \restriction C$ for any small enough club $C$ of
$\sigma$.  Now (a) holds whenever $C \subseteq \{i < \lambda:
i \mathrm{limit} \}$, clause (b) holds for $C \subseteq [i_0,\sigma)$ 
when $\theta^{+3} < \lambda_{i_0}$ and clause (c) holds as 
$\cov(\lambda_i,\lambda_i,\theta^+,2) < \lambda$ and use
\cite[Ch.II,5.3,10]{Sh:g} + Fodor's lemma and monotonicity of $\cov$.
  
Lastly, clause (d) holds as if $\{\mu < \lambda:\pp_\theta(\mu) 
\ge \lambda \cf(\mu) \le \theta\}$ is unbounded in 
$\lambda$, we get a contradiction by \cite[Ch.II,2.3(4)]{Sh:g}.]

Let $\lambda_\sigma =: \lambda$.  By \cite[Ch.VIII,1.6(3)]{Sh:g} we have 
(but shall not use)
\mn
\begin{enumerate}
\item[$(*)_1$]   if $\delta \le \sigma$ and $\cf(\delta) > \aleph_0$ 
\underline{then}  $\pp^+_\theta(\lambda_\delta) = \pp^{\cer}_{J^{\bd}_{\cf\lambda_\delta}}\!\!\!(\lambda_\delta)$ 
(and $\cf(\lambda_\delta) = \cf(\delta)$).
\end{enumerate}
\mn
Now by clause (d)
\mn
\begin{enumerate}
\item[$(*)_2$]   ${\ga} \subseteq \Reg \cap \lambda_i \setminus
\lambda_0$, $|{\ga}| \le \theta$ and $\sup({\ga}) \le \lambda_i$ 
\underline{implies} $\max \pcf({\ga}) < \pp^+_\theta(\lambda_i)$.
\end{enumerate}

Let

\[
S =: \big\{i \le \sigma:\cov(\lambda_i,\lambda_i,\theta^+,2) \ge
\pp^{\cer}_{J^{\bd}_{\cf\lambda_i}}\!\!\!(\lambda_i) \big\}.
\]

\mn
So it is enough to prove that $S$ is not stationary. 

Let for $i \le \sigma,\mu_i =: \pp^{\cer}_{J^{\bd}_{\cf(\lambda_i)}}
(\lambda_i)$, so $\lambda_{i+1} > \mu_i > \lambda_i,\mu_i$ 
is regular.   Note that $\mu_\sigma = pp_\theta(\lambda_\sigma) =
\pp^{\cer}_{J^{\bd}_\sigma}(\lambda_i)$ by \cite[Ch.VIII,1.6(3)]{Sh:g}. 

Clearly
\mn
\begin{enumerate}
\item[$(*)_3$]   $\lambda_i < \mu_i = \text{ cf}(\mu_i) \le \cov
(\lambda_i,\lambda_i,\theta^+,2)^+$.
\end{enumerate}
\mn
We can find $\bar A = \langle A_\zeta:\zeta < \lambda \rangle$ such that:
\mn
\begin{enumerate}
\item[$(*)_4(a)$]    $\zeta < \lambda_0 \Rightarrow A_\zeta = \varnothing$
\sn
\item[$(b)$]   $\lambda_i \le \zeta < \lambda_{i+1} \Rightarrow A_\zeta
\subseteq \lambda_i \text{ and } |A_\zeta| < \lambda_i$
\sn
\item[$(c)$]   for every $A \subseteq \lambda_i$ of cardinality 
$\le \theta$, for some $\zeta$, $\lambda_i < \zeta < \cov(\lambda_i,
\lambda_i,\theta^+,2)$ (which is $< \lambda_{i+1}$) we have $A \subseteq
A_\zeta$.
\end{enumerate}
\mn
Choose $\chi$ regular large enough, now choose by
induction on $i \le \sigma$ an elementary submodel $M^*_i$ of 
$(\clH(\chi)$, $\in,<^*_\chi)$, $\|M^*_i\| < \mu_i$, $M^*_i \cap \mu_i$ is an ordinal such that
\mn
\begin{enumerate}
\item[$(*)_5$]   if $i \le \sigma$, \underline{then}  
\[
\bigcup\limits_{j < i} M^*_j \cup \{\zeta:\zeta \le \lambda_i\} \cup
\{\langle \lambda_i : i < \sigma \rangle,\ \bar A,\ \langle M^*_j : j < i \rangle\}
\subseteq M^*_i.
\]
\end{enumerate}
\mn
Let ${\cP}_i = M^*_i \cap [\lambda_i]^{< \lambda_i}$.  It is enough to show
that

\begin{equation*}
\begin{array}{clcr}
S_1 = \{i \le \sigma:&\text{for some } Y \subseteq \lambda_i,\ 
|Y| \le \theta \text{ and} \\
  &Y \text{ is not a subset of any member of } {\cP}_i\}
\end{array}
\end{equation*}

\mn
is not stationary and $\sigma \notin S_1$ (in fact $S,S_1$ are equal). 

\noindent
[Why?  As clearly $S \subseteq S_1$.]

We assume $S_1$ is a stationary subset of $\sigma$ or $\sigma \in S_1$ and
eventually will finish by getting a contradiction.

For each $i \in S_1$ choose $Y_i \subseteq \lambda_i$ of cardinality 
$\le \theta$ which is not a subset of any member of ${\cP}_i$.  
Let $Y = \bigcup\limits_{i \in S_1} Y_i$, so $Y \subseteq \lambda$, 
$|Y| \le \theta$; and for each $i < \sigma$ we can find an ordinal 
$\zeta(i)$ such that $\lambda_i \le \zeta(i) < \cov(\lambda_i,\lambda_i,\theta,2)$ 
(which is $< \lambda_{i+1}$) and $Y \cap \lambda_i \subseteq A_{\zeta(i)}$.  
Now $|A_{\zeta(i)}| < \lambda_i$, hence by Fodor's Lemma 
there is $i(*) < \sigma$ such that

\[
S_2 =: \{i < \sigma:|A_{\zeta(i)}| < \lambda_{i(*)}\}.
\]

\mn
is a stationary subset of $\sigma$.  Let $Z =: \{\zeta(i):i \in S_2\}$.  Now
if $\sigma \in S_1$, then by \cite[Ch.IX,II,5.4]{Sh:g} and
\cite[Ch.II,\S1]{Sh:g} we have 
$$\pp^{\cer}_{J^{\bd}_\sigma}(\lambda) = \cov(\lambda,\lambda,\sigma^+,\sigma) =^+ 
\pp_{\Gamma(\sigma^+,\sigma)}(\lambda)$$ so there are
$j^* < \sigma$ and $B_j \in {\cP}_\sigma = M_\sigma \cap [\lambda]
^{< \lambda}$ for $j < j^*$ such that $Z \subseteq 
\bigcup\limits_{j < j^*} B_j$.  So for some $j<j^*$ we have 
$|Z \cap B_j| = \sigma$.  Now the set

\[
A^* = \bigcup\{A_\gamma:\gamma \in B_j,|A_j| \le \lambda_{i(*)}\}
\]

\mn
belongs to $M_\sigma$, has cardinality $\le \lambda_{i(*)} \times |B_j| <
\lambda$ and

\begin{equation*}
\begin{array}{clcr}
Y = &\bigcup\{Y \cap \lambda_i:i \in S_2 \text{ and } \zeta(i) \in B_j\}
\subseteq \\
  &\bigcup \{A_{\zeta(i)}:i \in S_2 \text{ and } \zeta(i) \in B_j\}
\subseteq A^* \in {\cP}_\sigma
\end{array}
\end{equation*}

\mn
contradiction.  So we have finished the case $\sigma \in S_1$ and from now
on we shall deal with the case $\sigma \notin S_1$ hence $S_1$ is a stationary
subset of $\sigma$, hence without loss of generality
$S_2 \subseteq S_1$.  Note that if $\delta < \sigma \text{ and } \cf(\delta) 
> \aleph_0$, we can apply this proof to $\lambda_\delta,\langle \lambda_i:
i < \delta \rangle$ (for $\sigma' = \cf(\delta))$ hence
\mn
\begin{enumerate}
\item[$(*)_6$]   $i \in S_2 \Rightarrow \text{ cf}(i) = \aleph_0$. 
\end{enumerate}
\mn
Clearly
\begin{enumerate}
\item[$(*)_7$]   for no $i \in S_2$ and $Z' \subseteq Z \cap \lambda_i$ is
$Z'$ unbounded in $\lambda_i$ and is contained in a member of $M^*_i$ of
cardinality $< \lambda_i$.
\end{enumerate}
\mn
Now we want to work as in the proof of \cite[CH.IX,3.5]{Sh:g}, but 
for $\sigma$ places at once with ``nice" behavior on a club of
$\sigma$, in the end the model is the Skolem Hull of the union 
of $\aleph_0$ sets, so one ``catches" an unbounded subsets of $Z$.  
Let $\bar \lambda = \langle \lambda_i:i \le \sigma \rangle$.

We shall choose by induction on $k < \omega$,

\[
N^a_k,N^b_k,g_k,\langle A^k_\ell:\ell < \omega \rangle,
\left< \langle A^k_{\ell,i}:i \le \sigma \rangle:\ell < \omega \right>
\]

\mn
such that:
\mn
\begin{enumerate}
\item[$(a)$]   for $x \in \{a,b\},N^x_k$ is an elementary submodel of
$(\clH(\chi),\in,<^*_\chi,\sigma,\bar \lambda)$ of cardinality $\le
\sigma$ and $N^x_k$ is the Skolem Hull of $N^x_k \cap \lambda$ and
$N^a_k \prec N^b_k$
\sn
\item[$(b)$]   $N^a_0[N^b_0]$ is the Skolem Hull of $\{i:i \le \sigma\}$
[of $Z \cup \{i:i \le \sigma\}$] in $(\clH(\chi),\in,<^*_\chi,
\sigma,\bar \lambda)$
\sn
\item[$(c)$]   $g_k \in \prod(\Reg \cap N^a_k \cap \lambda \setminus
\lambda^+_0)$
\sn
\item[$(d)$]   for $x \in \{a,b\}:N^x_{k+1} \text{ is the Skolem Hull of}$

\[
N^x_k \cup \{g_k(\kappa):\kappa \in \Dom(g_k)\} \cup (N^b_k \cap \lambda_0)
\]
\sn
\item[$(e)$]   $N^a_k \cap \lambda = \bigcup\limits_{\ell < \omega} A^k_\ell$
\sn
\item[$(f)$]   $A^k_\ell = \bigcup\limits_{i < \sigma} A^k_{\ell,i}$ and
$\langle A^k_{\ell,i}:i < \sigma \rangle$ is continuous increasing (in $i$)
and $A^k_{\ell,i} \subseteq \lambda_i$ and $|A^k_{\ell,i}| < \sigma$
\sn
\item[$(g)$]   if $\kappa \in \Reg \cap \lambda \cap N^a_k \setminus
\lambda^+_0$ then $\sup(N^b_k \cap \kappa) < g_k(\kappa) < \kappa$
\sn
\item[$(h)$]   if ${\ga} \subseteq A^k_\ell$ has order type $\omega$
and sup$({\mathfrak a}) = \lambda_i$ and ${\mathfrak a}$ is a subset of some
${\gb} \in M^*_i$ of cardinality $\le \lambda_0$, \underline{then} 
for some infinite ${\gb} \subseteq {\ga},g_k \restriction {\gb}$ 
is included in some function $h^k_{\ga} \in M^*_i$ such that 
$|\Dom(h^k_{\mathfrak a})| \le \lambda_0$.
\end{enumerate}
\mn
For $X \in \clH(\chi)$ and a function $F$ we let

\[
A(X,F) =: \{F(x_1,\dotsc,x_n):x_1,\dotsc,x_n \in X\}.
\]

\mn
Let us carry the induction for $k=0$; we define $N^a_0,N^b_0$ by 
clause (b) and define $\{A^0_\ell:\ell < \omega\}$ as

\[
\big\{A(\sigma +1,F):F \text{ a definable function in }
(\clH(\chi),\in,<^*_\chi,\sigma,\bar \lambda) \big\}.
\]

\mn
For $k+1$, let $g'_k \in \prod(\Reg \cap \lambda \cap N^a_k \setminus
\lambda^+_0)$ be defined by $g'_k(\kappa) = \sup(N^b_k \cap \kappa)$ 
(note: the domain of $g'_k$ is determined by $N^a_\kappa$, the values --- by
$N^b_k$).

We now shall find $g_k$ satisfying:
\mn
\begin{enumerate}
\item[$(\alpha)$]   $\Dom(g_k) = \Dom(g'_k),g_k \in \prod(\Dom(g'_k))$
\sn
\item[$(\beta)$]  $g'_k < g_k$
\sn
\item[$(\gamma)$]  if $i < \sigma,\ell < \omega$ and ${\mathfrak a} \subseteq
\Reg \cap A^k_\ell \setminus \lambda^+_0$ is unbounded in
$\lambda_i$ and is a subset of some $\gb \in M^*_i$ of cardinality
$\le \lambda_0$ and is of order type $\omega$, \underline{then}  for some 
infinite $\gb \subseteq {\ga}$ we have $g_k \restriction \gb$ 
is included in some $h_\gb \in M^*_i$ such that $|\Dom(h_\gb)| 
\le \lambda_0$
\sn
\item[$(\delta)$]   if ${\ga} \subseteq \lambda_i \cap \Reg \cap
A^k_{\ell,i} \setminus \lambda^+_0$ and ${\ga} \in M^*_{i+1}$ \underline{then} 
$g_k \rest {\ga} \subseteq h$ for some function from $M^*_{i+1}$.
\end{enumerate}
\mn
Note: a function choosing $\langle \bar f^{{\ga},\mu}:\mu \in
\pcf({\ga}) \rangle$ satisfying $(*)_{\ga}$ below for each
${\ga} \subseteq \Reg \cap \lambda \setminus \theta^+,
|{\ga}| \le \theta$ is definable in $(\clH(\chi),\in,<^*_\chi)$, so
each $M^*_i$ is closed under it where
\mn
\begin{enumerate}
\item[$(*)_{\ga}$]   $\bar f^{{\ga},\mu} = \langle 
f^{\ga,\mu}_\alpha:\alpha < \mu \rangle$ satisfies
\sn
\begin{enumerate}
\item[$(\alpha)$]   $f^{\ga,\mu}_\alpha \in \prod {\ga}$,
\sn
\item[$(\beta)$]   $\alpha < \beta \Rightarrow f^{\ga,\mu}_\alpha
<_{J_{< \mu}[\ga]} f^{{\ga},\mu}_\beta$
\sn
\item[$(\gamma)$]  if $\theta < \cf(\alpha) < \min(\ga)$ \underline{then} 
$f^{{\ga},\mu}_\alpha(\kappa) = \min\{\bigcup\limits_{\beta \in C} 
f^{{\ga},\mu}_\beta(\kappa):C \text{ a club of } \alpha\}$
\sn
\item[$(\delta)$]   if $f \in \prod {\ga}$ then for some $\alpha <
\mu$ we have $f < f^{{\ga},\mu}_\alpha \mod J_\mu[{\ga}]$.
\end{enumerate}
\end{enumerate}
\mn
Let $\langle \ga_{i,\zeta} : \zeta < \zeta_i \le \sigma^{\aleph_0}\rangle$ 
list the $\ga$ such that $\tcf(\prod\ga/J^{\bd}_\ga)$ is well defined 
and for some $n < \omega$, $\ga \subseteq A^k_n$, 
$\ga \subseteq \Reg \cap \lambda_i \setminus \lambda^+_0$, 
$\otp(\ga) = \omega$, $\lambda_i = \sup(\ga)$ and there is 
$\gb \subseteq \Reg \cap \lambda_i \setminus \lambda^+_0$, 
$\gb \in M^*_i$, $|\gb| \le \lambda_0$ such
that $\ga \subseteq \gb$: note that the number of such
$\ga$-s is $\le \sigma^{\aleph_0}$.
Let $\{\gb_{i,\zeta} : \zeta < \zeta_i \le \sigma^{\aleph_0}\}$ be such
that $\gb_{i,\zeta} \subseteq \Reg \cap \lambda_i \setminus
\lambda^+_0$, $\gb_{i,\zeta} \in M^*_i$, $|\gb_{i,\zeta}| \le \lambda_0$
and $\ga_{i,\zeta} \subseteq \gb_{i,\zeta}$.

So apply \cite[CH.VIII,\S1]{Sh:g}; i.e. let $\theta_1 = \theta + \sigma^{\aleph_0}$ choose 
$\langle M^k_\zeta : \zeta < \theta^{++}_1 \rangle$ increasing continuous, 
$M^k_\zeta \prec (\clH(\chi),\in,<^*_\chi)$, 
$\langle M^k_\zeta : \zeta \le \xi \rangle \in M^k_{\xi +1}$, 
$\|M^k_\zeta\| \le \lambda_0$ and $g'_k$, 
$\langle A^k_{\ell,i} : i < \sigma,\ell < \omega \rangle$, $Z$, 
$\langle \gb_{i,\zeta}:i \in S_2,\zeta < \zeta_i \rangle$ belong to $M^k_0$; 
and the function 
$g_k(\kappa) =: \sup(\kappa \cap \bigcup\limits_{\zeta < \theta^{++}_1 } M^k_\zeta)$ 
satisfies clauses $(\alpha),(\beta),(\gamma),(\delta)$ above.  
Now $N^a_{k+1},N^b_{k+1}$ are defined by clause (d).  Note that by the definition 
of $\mu'_i$ we have: for every $i < \sigma$, $\zeta < \zeta_i$, for some infinite 
$\ga = \ga^*_{i,\zeta} \subseteq \ga_{i,\zeta}$ we have $\mu_{i,\zeta}$, 
$k = \max \pcf(\ga) < \mu_i$.  Moreover $\prod\ga/J^{\bd}_\ga$ 
has true cofinality.  So our main demand on $g_k$ is:
$g_{k+1} \rest {\ga}^*_{i,\zeta} = 
f^{{\gb}_{i,\zeta},\mu_{i,\zeta,k}}_\delta \mod 
J^{\bd}_{{\ga}^*_{i,\zeta}}$ for a suitable $\delta$, so
$\delta = \sup(\mu_{i,\zeta,k} \cap M^k_{\theta^{++}})$ is O.K.  (For clause
$(\gamma)$ use (b) + (c) above.)

Now let $\{A^{k+1}_\ell:\ell < \omega\}$ be a list of:

\begin{equation*}
\begin{array}{clcr}
\Big\{\lambda \cap A\big(\bigcup\limits_{m \le n} A^k_m \cup \Rang
[g_k \rest \bigcup\limits_{m < n} A^k_m],F \big):&n < \omega \text{ and $F$ a} \\
  &\text{definable function in}\\ 
  &\big(\clH(\chi),\in,<^*_\chi,\theta, \bar \lambda\big) \Big\}
\end{array}
\end{equation*}

\mn
and if

\[
A^{k+1}_\ell = \lambda \cap A \big( \bigcup\limits_{m < n} A^k_m \cup \Rang
[g_k \rest \bigcup\limits_{m < n} A^k_m],F^{k+1}_\ell \big) \text{ and }
i < \sigma
\]

\mn
$g_k \rest (\bigcup\limits_{m < n} A^k_{m,i} \cap \lambda_i \setminus \lambda^+_0)$ 
is included in some function: $\Dom(h^k_{\ell,i}) = {\gb}^k_{\ell,i}$, 
$h^k_{\ell,i}(\kappa) = \sup(\kappa \cap M^k_{\theta^{++}})$.

Having finished the inductive definition note that:
\mn
\begin{enumerate}
\item[$(*)_8$]   $\bigcup\limits_k N^a_k \prec \bigcup\limits_k 
N^b_k \prec (\clH(\chi),\in,<^*_\chi,\theta,\bar \lambda)$.
\end{enumerate}
\mn 
[Why?  As $N^a_k \prec N^b_k \prec (\clH(\chi),\in,<^*_\chi,\theta,
\bar \lambda)$ by clause (a) and clause (d).]
\mn
\begin{enumerate}
\item[$(*)_9$]    $\bigcup\limits_k N^a_k \cap \lambda_0 = 
\bigcup\limits_k N^b_k \cap \lambda_0$.
\end{enumerate}
\mn
[Why?   $N^b_k \cap \lambda_0 \subseteq N^a_{k+1} \cap \lambda_0$ (see
clause (d)).]
\mn
\begin{enumerate}
\item[$(*)_{10}$]   if $\mu \in \Reg \cap \lambda^+ \setminus
\lambda^+_0$ and $\mu \in \bigcup\limits_k N^a_k$ then 
$\bigcup\limits_{k < \omega} N^a_k$ contains an unbounded subset 
of $\mu \cap \bigcup\limits_{k < \omega} N^b_k$.
\end{enumerate}
\mn 
[Why?  By clauses (d) + (g).]

So clearly (as usual)

\[
\bigcup\limits_k N^a_k \cap \lambda = \bigcup\limits_{k} N^b_k \cap \lambda.
\]

\mn
but $Z \subseteq N^b_0 \subseteq \bigcup\limits_{k < \omega} N^b_k$ and
$Z \subseteq \lambda$ hence $Z \subseteq \bigcup\limits_{k < \omega} N^a_k \cap 
\lambda$.  So for each $i \in S_2$, we can
find $\langle (\bar a^{i,k},w^{i,k},u^{i,k},\bar F^{i,k}):
k \le k(i) \rangle$ such that:
\mn
\begin{enumerate}
\item[$(a)$]   $\bar a^{i,k(i)} = \langle \zeta(i) \rangle$
\sn
\item[$(b)$]    $\bar a^{i,k} = \langle a^{i,k}_n:n < n^{i,k} \rangle$
\sn
\item[$(c)$]   each $a^{i,k}_n$ belongs to $N^a_k \cup (\lambda_0 \cap
N^b_{k+1})$
\sn
\item[$(d)$]   $w^{i,k} = \{n < n^{i,k}:a^{i,k}_n \in \lambda_0 \cap
N^b_{k+1}\}$
\sn
\item[$(e)$]   $u^{i,k} = \{n < n^{i,k}:a^{i,k}_n \in N^a_k \cap
\Reg \cap \lambda \setminus \lambda^+_0\}$
\sn
\item[$(f)$]   $\bar F^{i,k} = \langle F^{i,k}_n:n \in n^{i,k} \setminus
w^{i,k} \rangle$, and $F^{i,k}_n$ is a definable function in \\
$(\clH(\chi),\in,<^*_\chi,\theta,\bar \lambda)$
\sn
\item[$(g)$]   if $k > 0$, then $a^{i,k}_n = 
F^{i,k}_n(\ldots,a^{i,k-1}_m,\dotsc,g_{k-1}(a^{i,k-1}_{m'}),
\ldots)_{m < n^{i,k-1},m' \in u^{i,k-1}}$.
\end{enumerate}
\mn
Let $a^{i,k}_n \in A^k_{\ell(i,k,n)}$.  Note $(*)$ 
We can find stationary $S_3 \subseteq S_2$ such that:
\mn
\begin{enumerate}
\item[$(*)$]   if $i \in S_3$ then $k(i) = k(*)$ and for $k \le k(*)$
we have $n^{i,k} = n^k$, $w^{i,k} = w^k$, $u^{i,k} = u^k$, 
$\bar F^{i,k} = \bar F^k$, $\ell(i,k,n) = \ell(k,n)$.
\end{enumerate}
\mn
We can also find a stationary $S_4 \subseteq S_3$ such that:
\mn
\begin{enumerate}
\item[$(*)$]   if $i_1 < i_2$ are in $S_4$ then $a^{i_1,k}_n \in
A^k_{\ell(k,n),i_2}$
\sn
\item[$(**)$]   if $k < k(*),n \in u^k$ then $\langle a^{i,k}_n:i \in
S_4 \rangle$ is constant or strictly increasing and if it is strictly
increasing and its limit is $\ne \lambda$ (hence is $< \lambda$) then it
is $< \lambda_{\min(S_4)}$.
\end{enumerate}
\mn
Let $E = \{\delta < \sigma:\delta = \sup(\delta \cap S_4) \text{ and if }
n \in u^k, \text{ and } \langle a^{i,k}_n:i \in S_4 \rangle 
\text{ is strictly}$ 

$\qquad \qquad \quad \qquad \text{ increasing with limit } \lambda 
\text{ then } \langle a^{i,k}_n:i \in S_4 \cap \delta \rangle 
\text{ is strictly}$ 

$\qquad \quad \qquad \qquad \text{ increasing with limit } 
\lambda_\delta\}$. 

Now choose $\delta(*) \in E \cap S_1$, and choose $b$, a subset of
$\delta(*) \cap S_4$ of order type $\omega$ with limit
$\delta(*)$.  We can choose $b^{k,n} \in [b]^{\aleph_0}$ for $k \le k(*)$,
$n \le n^k$ such that: $b^{0,0} = b$, $b^{k,n+1} \subseteq b^{k,n}$, 
$b^{k+1,0} = b^{k,n^k}$, and if $n \in u^k$, 
$\langle a^{i,k}_n : i \in S_4 \rangle$ strictly increasing with 
limit $\lambda$ then $\prod\{a^{i,k}_n : i \in b^{k,n+1}\} / J^{\bd}_{b^{k,n+1}}$ 
has true cofinality which necessarily is $< \mu_{\delta(*)}$.  

So (recall $n^{k(*)} = 1$) $b^* = b^{k(*),1}$ is a subset of 
$S_4 \cap \delta(*)$ of order type $\omega$ with limit $\delta(*)$ and 
$b^* \subseteq b^{k,n}$ for $k \le k(*)$, $n \le n^k$ and $b^* \subseteq b^{k,n+1}$ 
hence $n \in u^k$ and 
$$\langle a^{i,k}_n : i \in S_4 \rangle \text{ strictly increasing }
\Rightarrow \mu_{\delta(*)} > \max \pcf\{a^{i,k}_n:i \in b^*\}.$$

Now we prove by induction on $k \le k(*)$ that for each $n < n^k$ for some
${\gB}_{k,n} \in M^*_{\delta(*)}$ with $\|{\gB}_{k,n}\| \le \lambda_0$ 
we have $\{a^{i,k}_n : i \in b^*\} \subseteq \gB_{k,n}$.  For $k=0$ 
clearly $A^0_{\ell(k,n)} \in M^*_{\delta(*)}$ has cardinality $\le \sigma$.
For $k > 0$, for each $n < n^k$ we use the ``$b^{k,n+1} \subseteq b$ and 
the choice of $g_{k-1}$ and clause $(*)_0(c)$.  So we get a
contradiction to $(*)_7$ so we are done.

\noindent
2) A variant of the proof of part (1).  First, it is enough to prove, for
each $i(*) < \sigma$ restrict ourselves to $S^* \cup \{\delta < \sigma$: the
cardinal appearing in $\otimes$ is $\ge \lambda_{i(*)}\}$, then \wilog \,
$i(*)=0$ and see that $\zeta_i \le \lambda_0$ is O.K.
\end{PROOF}

\begin{remark}
\label{ac.50}

\noindent
1) Note that if we just omit ``$\sigma^{\aleph_0} < \lambda"$ we still get 
that for a club of $\delta < \sigma$, $\cf(\delta) > \aleph_0$ or 
$\cf(\delta) = \aleph_0$ and $\pp^{\cer}_{J^{\bd}_\omega}(\lambda_\delta)$; 
if $< \cov(\lambda_i,\lambda_i,\theta^+,2)$ is still 
$\ge \lambda^{+ \lambda_\delta}_\delta$.
\end{remark}

\begin{conclusion}
\label{ac.53}

If $\mu$ is strong limit singular of uncountable cofinality \underline{then} 
 for a club of $\mu' < \mu$ we have $(2^{\mu'})^+ =^+ 
\pp^{\cer}_{J^{\bd}_{\cf(\mu')}}\!\!\!(\mu')$.
\end{conclusion}

\begin{conclusion}
\label{ac.56}

If $\beth_\delta$ is a singular cardinal of uncountable cofinality,
\underline{then}  for a club of $\alpha < \delta$, if $\cf(\alpha) = \aleph_0$
then
\mn
\begin{enumerate}
\item[$(*)_1$]   $2^{\beth_\alpha} =^+ \pp(\lambda)$
\sn
\item[$(*)_2$]   there is $S \subseteq {}^\omega(i_\alpha)$ of cardinality
$2^{\beth_\alpha}$ containing no perfect subset (and more --- see 
\cite[\S6]{Sh:355}).
\end{enumerate}
\end{conclusion}
\newpage

\section {Guessing clubs by countable $C$-s} \label{C4}

Recently\footnote{added Fall 2002} Zapletal \cite{Za0x} proved a beautiful theorem
\begin{theorem}
\label{n2}
If $I$ is a ``nice" (definition) of a $\sigma$-complete ideal 
on ${\clP}(\bbR)$ for suitable $\mathrm{LC}$ if $\ZFC + \mathrm{LC} \vdash \cov(I) = 
2^{\aleph_0}$ then $\ZFC + \mathrm{LC} \models \mathrm{Unif}(I) < \aleph_4$.
\end{theorem}

\noindent
He also showed that $\aleph_4$ cannot be replaced by $\aleph_2$.  The
$< \aleph_4$ comes from quoting guessing clubs.  The following shows 
we can replace $\aleph_4$ by $\aleph_3$ (other continuation see
\cite{Sh:561}, \cite{Sh:610}).
\begin{claim}
\label{n5}
Assume $\delta^* < \omega_1$ is a limit ordinal and 
$S \subseteq S^{\aleph_2}_{\aleph_0}$ is stationary.
\underline{then}  we can find $\bar C = \langle C_\alpha:\alpha \in S \rangle$
such that
\mn
\begin{enumerate}
\item[$(a)$]  $C_\alpha \subseteq S$
\sn
\item[$(b)$]  $C_\alpha \subseteq \alpha$
\sn
\item[$(c)$]   $\beta \in C_\alpha \Rightarrow \beta \in S \text{ and }
C_\beta = C_\alpha \cap \beta$ 
\sn
\item[$(d)$]  $\otp(C_\alpha) \le \delta^*$
\sn
\item[$(e)$]   for every club $E$ of $\omega_2$ the set 
$$\{\delta \in S : \delta = \sup(C_\delta),\ \delta^* = 
\otp(C_\delta), \text{and } C_\delta \subseteq E\}$$ is a 
stationary subset of $\omega_2$.
\end{enumerate}
\end{claim}

\begin{PROOF}{\ref{n5}} 
For each $\alpha < \omega_2$ choose $\langle a^\alpha_i : i < \omega_1 \rangle$; 
it will be an increasing continuous sequence
of countable subsets of $\alpha$ with union $\alpha$.  For each
$\alpha < \omega_2$ let 

\begin{align*}
    C^0_\alpha = \big\{i < \omega_1 :\ & i \text{ is a limit ordinal such that }\\
    &\forall \beta \in a^\alpha_i\ [a^\beta_i = a^\alpha_i \cap \beta],\\ 
    &\alpha < \omega_1 \Rightarrow \alpha \subseteq a^\alpha_i, \text{and } \\
    & j < i \Rightarrow \text{the closure of $a^\alpha_j$ is }\subseteq a^\alpha_i
    \cup \{\alpha\}\big\}.
\end{align*}

Clearly
\mn
\begin{enumerate}
\item[$(*)_1(a)$]   each $C^0_\alpha$ is a club of $\omega_1$
\sn
\item[$(b)$]   if $i \in C^0_\alpha$ and $\beta \in a^\alpha_i$ then
$i \in C^0_\beta$.
\end{enumerate}
\mn
Now
\mn
\begin{enumerate}
\item[$(*)_2$]   for some $\zeta = \zeta^* < \omega_1$, for every
club $E$ of $\omega_1$ the following set $F_{\zeta^*}(E) \cap S$ is non empty
where
\sn
\item[$\boxtimes$]   $F_\zeta(E)$ is the set of $\delta < \omega_2$ such
that:
\sn
\begin{enumerate}
\item[$(a)$]   $\delta = \otp(\delta \cap E \cap S)$
\sn
\item[$(b)$]   $\delta = \sup(E \cap \delta \cap S) =
\sup(a^\delta_\zeta \cap E \cap S)$
\sn
\item[$(c)$]  $\otp(a^\delta_\zeta \cap E \cap S)$ is divisible
by $\delta^*$
\sn
\item[$(d)$]   $\zeta \in C^0_\delta$ hence $\beta \in
a^\delta_\zeta \Rightarrow \zeta \in C^0_\beta$.
\end{enumerate}
\end{enumerate}
\mn
[Why does $(*)_2$ hold?  Otherwise for each $\zeta < \omega_1$ there is a
club $E_\zeta$ of $\omega_2$ such that $F_\zeta(E_\zeta) = \varnothing$.
Let $E^* = \bigcap\{E_\zeta:\zeta < \omega_1\} \setminus \omega_1$,
clearly $E^*$ is a club of $\omega_2$, 
and so is $E' = \{\delta < \omega_2:\delta = 
\otp(E^* \cap \delta \cap S)\}$ and choose $\delta \in E' \cap S$, 
exists as $E'$ is a club of $\omega_2$ and $S \subseteq S^2_0$ is 
stationary.  Easily the set 
$$C^* = \{\zeta < \omega_1 : \zeta \text{ limit, } \otp(a^\delta_\zeta \cap E \cap S)
\text{ is divisible by } \delta^*\}$$ is a club of $\omega_1$.  
So there is $\zeta^* \in C^* \cap C^0_\delta$, clearly
$\delta \in F_{\zeta^*}(E^*)$ hence $\delta \in F_{\zeta^*}(E_\zeta)$:
contradiction.]
\mn
\begin{enumerate}
\item[$(*)_3$]  if $E_1 \subseteq E_0$ are clubs of $\omega_2$ then
$F_{\zeta^*}(E_1) \subseteq F_{\zeta^*}(E_0)$. 
\end{enumerate}
\mn
[Why?  Note that $a \subseteq b \subseteq \delta =
\sup(a) \text{ and } \delta^*|\text{otp}(a)| \Rightarrow \delta^* \otp(b)$.]
\mn
\begin{enumerate}
\item[$(*)_4$]   for some club $E_0$ of $\omega_2$ for every club
$E_1$ of $\omega_2$ the set $F_{\zeta^*}(E_1,E_0) \ne \varnothing$ 
where $F_{\zeta^*}(E_1,E_0) = \{\delta:\delta \in F_{\zeta^*}(E_0)$ 
and $a^\delta_{\zeta^*} \cap E_0 \cap E_1 = a^\delta_{\zeta^*} \cap E_0\}$.
\end{enumerate}
\mn
[Why?  If not we choose by induction on $\eps < \omega_1$ a
club $E_\eps$ of $\omega_2$ such that $i <\eps
\Rightarrow E_\eps \subseteq E_i$ and $F_{\zeta^*}(E_{\eps
+1},E_\eps) = \varnothing$.  So $E^* =
\cap\{E_\eps:\eps < \omega_1\}$ is a club of
$\omega_2$ so we can find $\delta \in F_{\zeta^*}(E^*)$, hence $\delta
\in \cap\{F_{\zeta^*}(E_\eps):\eps < \omega_1\}$ by $(*)_3$.
Now trivially $\langle a^\delta_{\zeta^*} \cap
E_\eps:\eps < \omega_1 \rangle$ is a decreasing sequence
of subsets of $a^\delta_{\zeta^*}$ which is countable 
and $\eps < \omega_1 \Rightarrow a^\delta_{\zeta^*}
\cap E_\eps \ne a^\delta_{\zeta^*} \cap E_{\eps +1}$ 
as $F_{\zeta^*}(E_{\eps +1},E_\eps) = \varnothing$, contradiction.]

We fix $E_0$ as in $(*)_4$,
\mn
\begin{enumerate}
\item[$(*)_5$]   for some $\xi_\zeta < \omega_1$ we have: 
for every club $E_1$ of $\omega_2$ for some $\delta \in
F_{\zeta^*}(E_1,E_0)$ we have $\otp(a^\delta_{\zeta^*} \cap S \cap E_0) =\xi$. 
\end{enumerate}
\mn
[Why?  As in the proof of $(*)_4$.]

So necessarily $\xi$ is divisible by $\delta^*$.  Choose $b \subseteq
\xi = \sup(b)$, otp$(b) = \delta^*$.  Let

\[
S' = \{\alpha < \omega_1:\otp(a^\alpha_{\zeta^*} \cap S \cap
E_0) \in b \cup \{\xi\}\}.
\]

\mn
Now we define $\bar C = \langle C_\alpha:\alpha \in S \rangle$ as
follows: if $\alpha \in S \setminus S'$ we let $C_\alpha = \varnothing$
and if $\alpha \in S'$ we let $C_\alpha = 
\{\beta:\beta \in a^\alpha_{\zeta^*} \cap S \cap E_0$ and
$\otp(\beta \cap a^\alpha_{\zeta^*} \cap S \cap E_0) \in b\}$.

Now you can check that $\bar C = \langle C_\alpha:\alpha \in S
\rangle$ is as required.  (Noting that is clause (e), ``stationarily
many" ``at least one" are equivalent demands.)
\end{PROOF}

\begin{remark}
\label{n8}
Can we demand above that if $C_\alpha$ has no last
element then $C_\alpha$ is a closed subset of $\alpha$?

Not clear to me, but we can find
\mn
\begin{enumerate}
\item[$\circledast$]   there is $\langle {\cC}_\alpha:\alpha \in
S \rangle$ such that
\sn
\begin{enumerate}
\item[$(a)$]   ${\cC}_\alpha$ is a countable family of
countable subsets of $\alpha \cap S$, each of order type $\le \delta(*)$
\sn
\item[$(b)$]   if $C \in {\cC}_\alpha$ then $C$ is closed as a
subset of $\alpha$
\sn
\item[$(c)$]   if $\beta \in C \in {\cC}_\alpha$ then $C \cap
\beta \in {\cC}_\beta$
\sn
\item[$(d)$]  if $E$ is a club of $\omega_2$ then for
stationarily many $\alpha \in S$ for some $C \in {\cC}_\alpha$ we
have $\delta(*) = \otp(C) \text{ and } C \subseteq E$.
\end{enumerate}
\end{enumerate}
\end{remark}

\noindent
In some cases Zapletal \cite{Za0x}
uses ${\gd} \le {\gb}^{+n}$ we can
replace this by $\cf([{\gd}]^{\aleph_0},\subseteq) = {\gd}$ because
\begin{claim}
\label{n11}
Assume $\kappa$ is regular uncountable.

If $\lambda > \kappa$ and $\cf([\lambda]^{< \kappa},\subseteq) =
\lambda$, \underline{then}  for any $\alpha < \kappa$ there is $Y \subseteq
{}^\alpha([\lambda]^{< \kappa})$ which is
$E_{\kappa,\alpha}(\lambda)$-positive, i.e.,
\mn
\begin{enumerate}
\item[$\circledast$]   if $\chi > \lambda$ and $\xi \in 
\clH(\chi)$ then there is $\bar N = \langle N_i:i \le \alpha \rangle$
such that $x \in N_i \prec (\clH(\chi) \in <^*_\chi)$
\end{enumerate}

\[
\|N_i\| < \kappa
\]

\[
N_i \cap \kappa \in \kappa
\]

\[
N_i \text{ increasing continuous}
\]

\[
\bar N \restriction (i+1) \in N_{i+1}.
\]  
\end{claim}

\begin{PROOF}{\ref{n11}}
As in \cite[\S2]{Sh:420} \magenta{\textbf{}{(fill!)}}

Let $W_0 = \big\{\alpha \in E_0:\{\xi^*,\zeta^*\} \subseteq C^0_\alpha$
and $\otp(a^\alpha_{\xi^*} \cap E_0) \le \zeta \red{\big\}}$ and for $\alpha \in W_0$
let $b_\alpha = a^\alpha_{\xi^*} \cap E_0 \cap S$.  Clearly
\mn
\begin{enumerate}
\item[$\circledast)$]
\begin{enumerate}
\item[(a)]  $\alpha \in W_0 \Rightarrow b_\alpha
\subseteq W_0 \text{ and } \otp(b_\alpha) \le \zeta$
\sn
\item[(b)]  $\alpha \in b_\beta \cap \beta \in W \Rightarrow b^*_\alpha
= b_\beta \cap \alpha$
\sn
\item[(c)]   if $E_1$ is a club of $\omega_2$ then for
stationarily many $\alpha \in W_0 \cap S$ we have $b_\alpha \subseteq
E_1$ and $\otp(b_\alpha) = \xi$.
\end{enumerate}
\end{enumerate}
\end{PROOF}



\bibliographystyle{amsalpha}
\bibliography{shlhetal}

\end{document}